\documentclass{amsart}
\usepackage[all]{xy}
\usepackage{latexsym}
\usepackage{amssymb}
\usepackage{amsmath}
\usepackage{amsthm}
\usepackage{amscd}
\usepackage{fancyhdr}
\usepackage[dvips]{graphicx}
\usepackage{a4wide}
\usepackage{mathrsfs}
\input cyracc.def

 \newfam\cyrfam

  \font\tencyr=wncyr10

  \font\sevencyr=wncyr7

  \font\fivecyr=wncyr5

  \textfont\cyrfam=\tencyr \scriptfont\cyrfam=\sevencyr

    \scriptscriptfont\cyrfam=\fivecyr

  \newfam\cyifam

  \font\tencyi=wncyi10

  \font\sevencyi=wncyi7

  \font\fivecyi=wncyi5

  \textfont\cyifam=\tencyi \scriptfont\cyifam=\sevencyi

    \scriptscriptfont\cyifam=\fivecyi

\def\id{{\mbox{1 \hskip -7.2pt 1}}}
\newcommand{\sgn}{{\mathrm s  \mathrm g\mathrm  n}}
 \newcommand{\lon}{\longrightarrow}
 
 \newcommand{\rar}{\rightarrow}

\newcommand{\p}{{\partial}}

\newcommand{\Q}{{\mathbb Q}}
 \newcommand{\Z}{{\mathbb Z}}
 \newcommand{\bS}{{\mathbb S}}

 \newcommand{\R}{{\mathbb R}}
 \newcommand{\N}{{\mathbb N}}
 \newcommand{\K}{{\mathbb K}}
 \newcommand{\ot}{\otimes}

  \newcommand{\CoLie}{{\mathit C \mathit o \mathit L \mathit i\mathit e}}

\newcommand{\LB}{{\mathcal L} ie^1\hspace{-0.3mm} {\mathcal B}}

 \newcommand{\Beq}{\begin{equation}}
 \newcommand{\Eeq}{\end{equation}}
 \newcommand{\Beqr}{\begin{eqnarray}}
 \newcommand{\Eeqr}{\end{eqnarray}}
 \newcommand{\Beqrn}{\begin{eqnarray*}}
 \newcommand{\Eeqrn}{\end{eqnarray*}}
 \newcommand{\Ba}{\begin{array}}
 \newcommand{\Ea}{\end{array}}
 \newcommand{\Bi}{\begin{itemize}}
 \newcommand{\Ei}{\end{itemize}}
 \newcommand{\Bc}{\begin{center}}
 \newcommand{\Ec}{\end{center}}
 \newcommand{\fg}{{\mathfrak g}}

 \newcommand{\fG}{{\mathfrak G}}
\newcommand{\fs}{{\mathfrak s}}
\newcommand{\fS}{{\mathfrak S}}
 \newcommand{\f}{{\mathcal O}}
 \newcommand{\cA}{{\mathcal A}}
 \newcommand{\cB}{{\mathcal B}}
 \newcommand{\cC}{{\mathcal C}}
 
 \newcommand{\cE}{{\mathcal E}}
 \newcommand{\cF}{{\mathcal F}}
 \newcommand{\cG}{{\mathcal G}}
 
 \newcommand{\cI}{{\mathcal I}}
 
 \newcommand{\caL}{{\mathcal L}}
 \newcommand{\cM}{{\mathcal M}}
 \newcommand{\cN}{{\mathcal N}}
 \newcommand{\cP}{{\mathcal P}}
 \newcommand{\cQ}{{\mathcal Q}}
 \newcommand{\cR}{{\mathcal R}}
 \newcommand{\cS}{{\mathcal S}}
 \newcommand{\cT}{{\mathcal T}}

 \newcommand{\cU}{{\mathcal U}}

 \newcommand{\al}{\alpha}
 \newcommand{\be}{\beta}
 \newcommand{\ga}{\gamma}
 
 \newcommand{\Ga}{\Gamma}
 \newcommand{\bGa}{{\mathbf \Gamma}}

 \newcommand{\la}{\lambda}
 \newcommand{\om}{\omega}

 \newcommand{\Ker}{{\mathsf K \mathsf e \mathsf r}\, }
 \newcommand{\Img}{{\mathrm I\mathrm m}\, }
 \newcommand{\Hom}{{\mathrm H\mathrm o\mathrm m}}

 \def\hgw{{\mbox {\large $\circlearrowright$}}}
 \newcommand{\sip}{\smallskip}
 
 \newcommand{\mip}{\vspace{2.5mm}}

\theoremstyle{plain}
\swapnumbers

\newtheorem{prop-def}[theorem]{Proposition-definition}

\newtheorem{main-theorem}{Main~Theorem}[section]
\newtheorem{section-theorem}{Theorem}[section]
\newtheorem{section-corollary}{Corollary}[section]

\theoremstyle{definition}

\begin{document}

 \sloppy

\long\def\symbolfootnote[#1]#2{\begingroup%
\def\thefootnote{\fnsymbol{footnote}}\footnote[#1]{#2}\endgroup}

 \title{Wheeled props in  algebra, geometry and quantization}
 \author{ S.A.\ Merkulov}
\address{Sergei~A.~Merkulov: Department of Mathematics, Stockholm University, 10691 Stockholm, Sweden}
\email{sm@math.su.se}
 \date{}

\begin{abstract}
Wheeled props is one  the latest species found in the world of operads and props. We attempt to give
 an elementary introduction into the main ideas of the  theory of wheeled props for beginners, and
 also a survey
of its most recent major applications (ranging from algebra and geometry to deformation theory and
Batalin-Vilkovisky quantization) which might be of interest to experts.
\end{abstract}
 \maketitle

\section{Introduction}
The theory of operads and props undergoes a rapid development in recent years; its applications
can be seen nowadays almost everywhere --- in algebraic topology,
 in homological algebra, in differential geometry, in non-commutative geometry, in string topology, in deformation theory,
 in quantization theory etc.
 The theory demonstrates a remarkable unity of mathematics;
 for example, one and the same {\em operad of little 2-disks}\,   solves
the recognition problem for based 2-loop spaces in algebraic topology,
describes homotopy Gerstenhaber
structure on the Hochschild deformation complex in homological algebra, and also
controls diffeomorphism invariant Hertling-Manin's integrability equations
in differential
geometry!

 \sip

First examples of operads and props were constructed
 in the 1960s in the classical papers of Gerstenhaber on deformation theory of algebras
 and of Stasheff on homotopy theory
 of loop spaces. The notion of prop was introduced by MacLane already in 1963 as a
 useful way to code axioms for operations with many inputs and outputs. The notion of operad was
 ultimately coined 10 years later by
 P.May through axiomatization of  properties of earlier discovered {\em associahedra}\, polytopes and
 the
 associated $A_\infty$-{\em spaces}\, by Stasheff and of the {\em little cubes operad}\, by Boardman
 and Vogt.

\sip

In this paper we attempt to explain the main ideas and constructions of the theory
of wheeled operads and props and illustrate them with some of the most recent applications
[Gr1, Gr2, Me1-Me7, MMS, MeVa, Mn, Str1, Str2]
to geometry, deformation
theory and
Batalin-Vilkovisky quantization formalism of theoretical physics. In the heart of these applications
lies  the fact
 that some categories of local geometric and theoretical physics
structures can be identified with the derived categories of surprisingly simple algebraic structures.
The language of graphs is essential for the proof of this fact  and permits
us to reformulate  it as follows: solution spaces of several
 important highly non-linear differential equations in geometry and physics
are controlled by (wheeled) props which are resolutions of very compact graphical data, a kind of
``genome"; for example, the ``genome" of the
 species {\em local Poisson
structures}\, is the prop of Lie 1-bialgebras built from two ``genes",
$
 \begin{xy}
 <0mm,-0.55mm>*{};<0mm,-3.5mm>*{}**@{-},
 <0.5mm,0.5mm>*{};<2.2mm,3.2mm>*{}**@{-},
 <-0.48mm,0.48mm>*{};<-2.2mm,3.2mm>*{}**@{-},
 <0mm,0mm>*{\circ};<0mm,0mm>*{}**@{},
 \end{xy}$ and
 $\begin{xy}
 <0mm,0.66mm>*{};<0mm,4mm>*{}**@{-},
 <0.39mm,-0.39mm>*{};<2.2mm,-3.2mm>*{}**@{-},
 <-0.35mm,-0.35mm>*{};<-2.2mm,-3.2mm>*{}**@{-},
 <0mm,0mm>*{\bullet};<0mm,0mm>*{}**@{},
 \end{xy}$,
subject to the following engineering  rules (see \S 2 for precise details),
$$
\Ba{c}
\begin{xy}
 <0mm,0mm>*{\circ};<0mm,0mm>*{}**@{},
 <0mm,-0.49mm>*{};<0mm,-3.0mm>*{}**@{-},
 <0.49mm,0.49mm>*{};<3.9mm,4.9mm>*{}**@{-},
 <-0.5mm,0.5mm>*{};<-1.9mm,1.9mm>*{}**@{-},
 <-2.3mm,2.3mm>*{\circ};<-2.3mm,2.3mm>*{}**@{},
 <-1.8mm,2.8mm>*{};<0mm,4.9mm>*{}**@{-},
 <-2.8mm,2.9mm>*{};<-4.6mm,4.9mm>*{}**@{-},
 \end{xy}
-
\begin{xy}
 <0mm,0mm>*{\circ};<0mm,0mm>*{}**@{},
 <0mm,-0.49mm>*{};<0mm,-3.0mm>*{}**@{-},
 <0.49mm,0.49mm>*{};<1.8mm,1.8mm>*{}**@{-},
 <1.8mm,1.8mm>*{};<-2.1mm,4.9mm>*{}**@{-},
 <-0.5mm,0.5mm>*{};<-1.9mm,1.9mm>*{}**@{-},
 <-2.3mm,2.3mm>*{\circ};<-2.3mm,2.3mm>*{}**@{},
 <-1.8mm,2.8mm>*{};<0mm,4.9mm>*{}**@{-},
 <-2.8mm,2.9mm>*{};<-4.6mm,4.9mm>*{}**@{-},
 \end{xy}
-
\begin{xy}
 <0mm,0mm>*{\circ};<0mm,0mm>*{}**@{},
 <0mm,-0.49mm>*{};<0mm,-3.0mm>*{}**@{-},
 <0.49mm,0.49mm>*{};<1.9mm,1.9mm>*{}**@{-},
 <-0.4mm,0.4mm>*{};<-3.2mm,4.9mm>*{}**@{-},
 <2.3mm,2.3mm>*{\circ};<2.3mm,2.3mm>*{}**@{},
 <1.9mm,2.8mm>*{};<0.5mm,4.9mm>*{}**@{-},
 <2.8mm,2.9mm>*{};<4.6mm,4.9mm>*{}**@{-},
 \end{xy}
 \Ea
=0, \ \ \
\Ba{c}
 \begin{xy}
 <0mm,0mm>*{\bullet};<0mm,0mm>*{}**@{},
 <0mm,0.69mm>*{};<0mm,3.0mm>*{}**@{-},
 <0.39mm,-0.39mm>*{};<1.9mm,-1.9mm>*{}**@{-},
 <-0.35mm,-0.35mm>*{};<-4.2mm,-4.9mm>*{}**@{-},
 <2.4mm,-2.4mm>*{\bullet};<2.4mm,-2.4mm>*{}**@{},
 <2.0mm,-2.8mm>*{};<0mm,-4.9mm>*{}**@{-},
 <2.8mm,-2.9mm>*{};<4.6mm,-4.9mm>*{}**@{-},
 \end{xy}
+ \,
\begin{xy}
 <0mm,0mm>*{\bullet};<0mm,0mm>*{}**@{},
 <0mm,0.69mm>*{};<0mm,3.0mm>*{}**@{-},
 <0.39mm,-0.39mm>*{};<1.9mm,-1.9mm>*{}**@{-},
 <0mm,0mm>*{};<-1.5mm,-1.5mm>*{}**@{-},
 <-1.5mm,-1.5mm>*{};<2.4mm,-4.9mm>*{}**@{-},
 <2.4mm,-2.4mm>*{\bullet};<2.4mm,-2.4mm>*{}**@{},
 <2.0mm,-2.8mm>*{};<0mm,-4.9mm>*{}**@{-},
 <2.8mm,-2.9mm>*{};<4.8mm,-4.9mm>*{}**@{-},
 \end{xy}
+
 \begin{xy}
 <0mm,0mm>*{\bullet};<0mm,0mm>*{}**@{},
 <0mm,0.69mm>*{};<0mm,3.0mm>*{}**@{-},
 <0.39mm,-0.39mm>*{};<4.3mm,-4.9mm>*{}**@{-},
 <-0.35mm,-0.35mm>*{};<-1.9mm,-1.9mm>*{}**@{-},
 <-2.4mm,-2.4mm>*{\bullet};<-2.4mm,-2.4mm>*{}**@{},
 <-2.0mm,-2.8mm>*{};<0mm,-4.9mm>*{}**@{-},
 <-2.8mm,-2.9mm>*{};<-4.7mm,-4.9mm>*{}**@{-},
 \end{xy}
\Ea
=0, \ \ \
\Ba{c}
 \begin{xy}
 <0mm,2.47mm>*{};<0mm,-0.5mm>*{}**@{-},
 <0.5mm,3.5mm>*{};<2.2mm,5.2mm>*{}**@{-},
 <-0.48mm,3.48mm>*{};<-2.2mm,5.2mm>*{}**@{-},
 <0mm,3mm>*{\circ};<0mm,3mm>*{}**@{},
  <0mm,-0.8mm>*{\bullet};<0mm,-0.8mm>*{}**@{},
<0mm,-0.8mm>*{};<-2.2mm,-3.5mm>*{}**@{-},
 <0mm,-0.8mm>*{};<2.2mm,-3.5mm>*{}**@{-},
\end{xy}
\Ea
 =
\Ba{c}
\begin{xy}
 <0mm,-1.3mm>*{};<0mm,-3.5mm>*{}**@{-},
 <0.38mm,-0.2mm>*{};<2.2mm,2.2mm>*{}**@{-},
 <-0.38mm,-0.2mm>*{};<-2.2mm,2.2mm>*{}**@{-},
<0mm,-0.8mm>*{\circ};-<0mm,0.8mm>*{}**@{},
 <-2.25mm,2.2mm>*{};<-2.2mm,5.2mm>*{}**@{-},
 <2.4mm,2.4mm>*{\bullet};<2.4mm,2.4mm>*{}**@{},
 <2.5mm,2.3mm>*{};<4.4mm,-0.8mm>*{}**@{-},
 <4.4mm,-0.8mm>*{};<4.4mm,-3.5mm>*{}**@{-},
 <2.4mm,2.5mm>*{};<2.4mm,5.2mm>*{}**@{-},
 \end{xy}
 +
\begin{xy}
<3mm,-1.3mm>*{};<3mm,-3.5mm>*{}**@{-},
<3.4mm,-0.2mm>*{};<5.2mm,2.2mm>*{}**@{-},
<2.6mm,-0.2mm>*{};<0.4mm,2.2mm>*{}**@{-},
<3mm,-0.8mm>*{\circ};<3mm,-0.8mm>*{}**@{},
<5.2mm,2.2mm>*{};<5.2mm,5.2mm>*{}**@{-},
<0.3mm,2.4mm>*{\bullet};<0.3mm,2.4mm>*{}**@{},
<0.1mm,2.3mm>*{};<-2.2mm,-0.8mm>*{}**@{-},
<-2.2mm,-0.8mm>*{};<-2.2mm,-3.5mm>*{}**@{-},
<0.3mm,2.5mm>*{};<0.3mm,5.2mm>*{}**@{-},
\end{xy}
 +
\begin{xy}
<3mm,-0.8mm>*{\circ};<3mm,-0.8mm>*{}**@{},
<3mm,-1.3mm>*{};<3mm,-3.5mm>*{}**@{-},
<3.7mm,-0.5mm>*{};<5.7mm,1.2mm>*{}**@{-},
<5.7mm,1.2mm>*{};<3.2mm,2.4mm>*{}**@{-},
<2.5mm,-0.5mm>*{};<-0.5mm,2.4mm>*{}**@{-},
<-0.5mm,2.4mm>*{};<-0.5mm,5.2mm>*{}**@{-},
<3mm,2.4mm>*{\bullet};<3mm,2.4mm>*{}**@{},
<3.3mm,2.9mm>*{};<1.4mm,1.0mm>*{}**@{-},
<0.9mm,0.5mm>*{};<-0.5mm,-0.8mm>*{}**@{-},
<-0.5mm,-0.8mm>*{};<-0.5mm,-3.5mm>*{}**@{-},
<3mm,2.4mm>*{};<3mm,5.2mm>*{}**@{-},
\end{xy}
 +
\begin{xy}
<3mm,-0.8mm>*{\circ};<3mm,-0.8mm>*{}**@{},
<3mm,-1.3mm>*{};<3mm,-3.5mm>*{}**@{-},
<2.3mm,-0.5mm>*{};<0.3mm,1.2mm>*{}**@{-},
<0.3mm,1.2mm>*{};<2.4mm,2.4mm>*{}**@{-},
<3.5mm,-0.5mm>*{};<7mm,2.4mm>*{}**@{-},
<7mm,2.4mm>*{};<7mm,5.2mm>*{}**@{-},
<3mm,2.4mm>*{\bullet};<3mm,2.4mm>*{}**@{},
<2.7mm,2.9mm>*{};<4.8mm,1.0mm>*{}**@{-},
<5.4mm,0.5mm>*{};<7mm,-0.8mm>*{}**@{-},
<7mm,-0.8mm>*{};<7mm,-3.5mm>*{}**@{-},
<3mm,2.4mm>*{};<3mm,5.2mm>*{}**@{-},
\end{xy}
\Ea
$$
We shall explain how a slight modification of the above rules by addition of two extra conditions,
$
\Ba{c}
\begin{xy}
 <0mm,-0.55mm>*{};<0mm,-2.5mm>*{}**@{-},
 <0.5mm,0.5mm>*{};<2.2mm,2.2mm>*{}**@{-},
 <-0.48mm,0.48mm>*{};<-2.9mm,3.2mm>*{}**@{-},
 <0mm,0mm>*{\circ};<0mm,0mm>*{}**@{},
(2.2,2.2)*{}
   \ar@{->}@(ur,d) (0,-1)*{}
 \end{xy}=0$ and
$ \begin{xy}
 <0mm,0.66mm>*{};<0mm,3mm>*{}**@{-},
 <0.39mm,-0.39mm>*{};<2.2mm,-2.2mm>*{}**@{-},
 <-0.35mm,-0.35mm>*{};<-2.8mm,-3.2mm>*{}**@{-},
 <0mm,0mm>*{\bullet};<0mm,0mm>*{}**@{},
(0.0,1.0)*{}
   \ar@{->}@(u,dr) (2.2,-2.2)*{}
\end{xy}=0,
\Ea$,
changes the resulting ``species" dramatically: instead of the category of local Poisson structures
 one gets the category of {\em quantum BV manifolds with split quasi-classical limit
 }  which, for example, naturally emerges in the study of  quantum master equations \cite{BV, Sc}
 for $BF$-type quantum field theories (see \S5 for precise details).
 Moreover, in the homotopy theory sense, this category
 is as perfect as, for example, the nowadays famous category of $\caL ie_\infty$-algebras:
 quasi-isomorphisms of quantum BV manifolds turn out to be equivalence relations.

\sip

It is yet to see how non-trivial  topology can be incorporated into the current pro(p)file
of {\em local}\,
differential geometry, but it is worth stressing already now that this approach
to geometry and physics turns  space-time ---
``the  background of everything" ---
into an ordinary observable, a certain function (representation) on a prop and hence
 unveils a possibility for a new architecture:
\vspace{1mm}
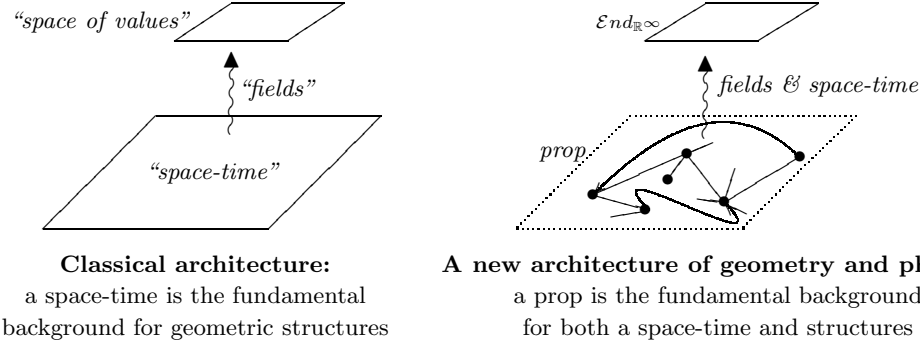
$$
\underset{\mbox{\tabular{c}\small\bf Classical architecture:\\
 \small a space-time is the fundamental \\
\small background for geometric structures  \endtabular }}
{
\Ba{c}
\begin{xy}
 <0mm,0mm>*{};<30mm,0mm>*{}**@{-},
 <0mm,0mm>*{};<15mm,15mm>*{}**@{-},
 <15mm,15mm>*{};<45mm,15mm>*{}**@{-},
 <0mm,0mm>*{};<20mm,0mm>*{}**@{-},
 <30mm,0mm>*{};<45mm,15mm>*{}**@{-},
  <22.5mm,7.5mm>*{\mbox{\small\em ``space-time"}};
 <25mm,12.5mm>*{};<25mm,22.5mm>*{}**@{~},
 <25mm,22.5mm>*{\blacktriangle};
<31mm,18.5mm>*{\mbox{\small\em ``fields"}};
<17.5mm,25mm>*{};<32.5mm,25mm>*{}**@{-},
 <17.5mm,25mm>*{};<25mm,30mm>*{}**@{-},
 <25mm,30mm>*{};<40mm,30mm>*{}**@{-},
 <32.5mm,25mm>*{};<40mm,30mm>*{}**@{-},
<7.5mm,27.5mm>*{\mbox{\small\em ``space of values"}};
\end{xy}
\Ea
}
\ \ \
\underset{\mbox{\tabular{c}\small\bf A new architecture of geometry and physics:\\
\small a prop is the fundamental background\\
\small  for both a space-time and structures\endtabular }}
{
\Ba{c}
\begin{xy}
 <0mm,0mm>*{};<30mm,0mm>*{}**@{.},
 <0mm,0mm>*{};<15mm,15mm>*{}**@{.},
 <15mm,15mm>*{};<45mm,15mm>*{}**@{.},
 <30mm,0mm>*{};<45mm,15mm>*{}**@{.},
  <6mm,10mm>*{\mbox{\small \em prop}};
 <25mm,12mm>*{};<25mm,22mm>*{}**@{~},
 <25mm,22mm>*{\blacktriangle};
<40mm,19mm>*{\mbox{{\small\em  fields \& space-time}}};
<17mm,25mm>*{};<32mm,25mm>*{}**@{-},
 <17mm,25mm>*{};<25mm,30mm>*{}**@{-},
 <25mm,30mm>*{};<40mm,30mm>*{}**@{-},
 <32mm,25mm>*{};<40mm,30mm>*{}**@{-},
<15mm,27mm>*{^{\cE nd_{\R^\infty}}};
<17mm,2.5mm>*{\bullet};
<27.5mm,3.5mm>*{\bullet};
<20mm,6.5mm>*{\bullet};
<37.5mm,9.5mm>*{\bullet};
<10mm,4.5mm>*{\bullet};
<22.5mm,10mm>*{\bullet};
<37.5mm,9.5mm>*{};<27.5mm,3.5mm>*{}**@{-},
<30.5mm,2mm>*{};<27.5mm,3.5mm>*{}**@{-},
<22.5mm,10mm>*{};<26mm,11.5mm>*{}**@{-},
<27.5mm,3.5mm>*{};<24mm,4mm>*{}**@{-},
<27.5mm,3.5mm>*{};<29mm,8mm>*{}**@{-},
<27.5mm,3.5mm>*{};<22.5mm,10mm>*{}**@{-},
<17.5mm,2.5mm>*{};<10mm,4.5mm>*{}**@{-},
<17.5mm,2.5mm>*{};<12.5mm,1.5mm>*{}**@{-},
<22.5mm,10mm>*{};<10mm,4.5mm>*{}**@{-},
<20mm,6.5mm>*{};<22.5mm,10mm>*{}**@{-},
 (17.5,2.5)*{}
   \ar@{->}@(ul,ur) (37.5,9.5)*{};
 (10.0,4.5)*{}
   \ar@{->}@(ul,dr) (27.5,3.5)*{}
\end{xy}
\Ea
}
$$
In fact, some elements of this architecture have been envisaged long ago by Roger Penrose
\cite{Pe} in his ``abstract index calculus".

\sip

The paper is organized as follows. In section 2 we  give a short but self-contained
introduction into the theory of (wheeled) operads and props. Sections 3 and 4 aim to give
an account of  most recent applications of that theory to geometry and, respectively, deformation
theory. In Section 5 we explain some ideas of Koszul duality theory and its relation
to the homotopy transfer formulae and Batalin-Vilkovisky formalism.

\sip

A few words about notations.
The symbol $\bS_n$ stands for the permutation group, i.e.\ for the group of all bijections,
$[n]\rar [n]$, where $[n]$ denotes (here and everywhere) the set
$\{1,2,\ldots,n\}$.
Given a partition,  $[n]=I_1\sqcup\ldots \sqcup I_k$,
the symbol $\sigma(I_1,\ldots,I_k)$ denotes the sign of the permutation $[n]\rar\{I_1,\ldots, I_k\}$.
 If $V=\oplus_{i\in \Z} V^i$ is a graded vector space, then
$V[k]$ is a graded vector space with $V[k]^i:=V^{i+k}$.
We work throughout over a field $\K$ of characteristic 0.

\section{An introduction to operads, dioperads, properads and props}

\subsection{Directed graphs} Let $m$ and $n$ be arbitrary non-negative integers.
 A {\em directed  $(m,n)$-graph} is a triple $(G,f_{in},f_{out})$, where $G$ is a finite
 $1$-dimensional
$CW$ complex all of whose every 1-dimensional cells (``edges") are oriented (``directed"), and
$$
f_{in}: [m] \rar \left\{
\Ba{c}
\mbox{\small the set of all $0$-cells, $v$, of $G$}\\
\mbox{\small  which have precisely one}\\
\mbox{\small  adjacent edge directed from $v$}
\Ea
\right\},\ \
f_{out}: [n] \rar \left\{
\Ba{c}
\mbox{\small the set of all $0$-cells, $v$, of $G$}\\
\mbox{\small  which have precisely one}\\
\mbox{\small  adjacent edge directed towards $v$}
\Ea
\right\}
$$
are injective maps of finite sets (called {\em labelling maps}\, or simply {\em labellings}) such that $\Img f_{in}\cap \Img f_{out}=\emptyset$.
The set, $\fG^\circlearrowright(m,n)$, of all
possible directed $(m,n)$-graphs carries an action,
$
(G, f_{in}, f_{out}) \rar (G, f_{in}\circ \sigma^{-1} , f_{out}\circ \tau)$,
of the group  $\bS_m\times \bS_n$ (more precisely, the {\em right}\, action of $\bS_m^{op}\times \bS_n$
but we declare this detail implicit from now). We often abbreviate a triple  $(G,f_{in},f_{out})$
to $G$.
For any $G\in \fG^\circlearrowright(m,n)$
the set,
$$
V(G):=\{\mbox{all 0-cells of}\ G\}\setminus\{\Img f_{in}\cup \Img f_{out}\},
$$
of all unlabelled
$0$-cells is called the set
of {\em vertices}\, of $G$. The edges attached to labelled $0$-cells, i.e.\ the ones
 lying in $\Img f_{in}$ or in $\Img f_{out}$ are called {\em incoming}\,
  or, respectively, {\em outgoing legs}\,
of the graph $G$. The set
$$
E(G):=\{\mbox{all 1-cells of}\ G\}\setminus\{\mbox{legs}\},
$$
is called the set of {\em (internal) edges}\, of $G$. Legs and edges of $G$ incident to a vertex
$v\in V(G)$ are often called {\em half-edges}\, of $v$; the set of half-edges of $v$ splits naturally
 into two disjoint
sets, $In_v$ and ${\mathit O\mathit u\mathit t}_v$, consisting of incoming and, respectively,
 outgoing half-edges.
 In all our pictures the vertices of a
graph will be denoted by bullets, the edges by intervals (or sometimes curves) connecting the vertices,
and legs by intervals attached from one side to vertices. A choice of orientation on an edge or a leg will
be visualized by the choice of a particular direction (arrow) on the associated interval/curve;
unless otherwise explicitly shown the direction of each edge in all our pictures is assumed to go {\em
from bottom to the top}. For example, the graph
$
\Ba{c}
\begin{xy}
 <0mm,0mm>*{\bullet};
<0.39mm,0.39mm>*{};<3.4mm,3.4mm>*{}**@{-},
<0mm,0.39mm>*{};<0mm,3.4mm>*{}**@{-},
 <0.39mm,-0.39mm>*{};<3.4mm,-3.7mm>*{}**@{-},
 <-0.35mm,-0.35mm>*{};<-2.9mm,-2.9mm>*{}**@{-},
 <-3.4mm,-3.4mm>*{\bullet};
 <-3.4mm,-3.4mm>*{};<0mm,-6.8mm>*{}**@{-},
 <3.4mm,-3.7mm>*{};<0mm,-6.8mm>*{}**@{-},
  <0mm,-6.8mm>*{\bullet};
  <0mm,-6.8mm>*{};<0mm,-10mm>*{\bullet}**@{-},
  <0mm,-6.8mm>*{};<3mm,-10mm>*{}**@{-},
<0mm,-6.8mm>*{};<-3mm,-10mm>*{}**@{-},
  <0.4mm,-8.5mm>*{};<3.8mm,-12mm>*{^2}**@{},
 <0.4mm,-6.5mm>*{};<-3.8mm,-12mm>*{^1}**@{},
 <0mm,0mm>*{};<0mm,4.5mm>*{^1}**@{},
 <0mm,0mm>*{};<3.6mm,4.5mm>*{^2}**@{},
(-0.39,0.39)*{}
   \ar@{->}@(ul,dl) (-3.6,-3.6)*{}
 \end{xy}
\Ea\in \fG^\circlearrowright(2,2)
$
has four vertices, four legs and five edges; the orientation of all legs and of four
internal edges is {\em not}\, shown explicitly and hence, by default, flows {\em upwards}.
Sometimes we skip showing explicitly
 labellings of legs (as in Table 1, for example).
We set $\fG^\circlearrowright:=\sqcup_{m,n\geq 0} \fG^\circlearrowright(m,n)$.
Note that elements of $\fG^\circlearrowright$ are not necessarily connected,
e.g.\
$\Ba{c}
\begin{xy}
 <0mm,0mm>*{\bullet};
 <0.39mm,-0.39mm>*{};<2.4mm,-2.4mm>*{}**@{-},
 <-0.35mm,-0.35mm>*{};<-1.9mm,-1.9mm>*{}**@{-},
 <-2.4mm,-2.4mm>*{\bullet};
 <-2.0mm,-2.8mm>*{};<-0.4mm,-4.5mm>*{}**@{-},
 <2.4mm,-2.4mm>*{};<0.4mm,-4.5mm>*{}**@{-},
  <0mm,-5.1mm>*{\bullet};
  <0.4mm,-5.5mm>*{};<2mm,-7.7mm>*{}**@{-},
<-0.4mm,-5.5mm>*{};<-2mm,-7.7mm>*{}**@{-},
  <0.4mm,-5.5mm>*{};<2.9mm,-9.7mm>*{^2}**@{},
 <0.4mm,-5.5mm>*{};<-2.9mm,-9.7mm>*{^3}**@{},
(0,0)*{}
   \ar@{->}@(ul,dl) (-2.4,-2.4)*{}
 \end{xy}
\Ea
\begin{xy}
 <0mm,-1.3mm>*{};<0mm,-3.5mm>*{}**@{-},
 <0.38mm,-0.2mm>*{};<2.2mm,2.2mm>*{}**@{-},
 <-0.38mm,-0.2mm>*{};<-2.2mm,2.2mm>*{}**@{-},
<0mm,-0.8mm>*{\bullet};
 <2.4mm,2.4mm>*{\bullet};
 <2.5mm,2.3mm>*{};<4.4mm,-0.8mm>*{}**@{-},
 <2.4mm,2.8mm>*{};<2.4mm,5.2mm>*{}**@{-},
     <0mm,-1.3mm>*{};<0mm,-5.3mm>*{^1}**@{},
     <2.5mm,2.3mm>*{};<5.1mm,-2.7mm>*{^4}**@{},
    <2.4mm,2.5mm>*{};<2.4mm,5.7mm>*{^2}**@{},
    <-0.38mm,-0.2mm>*{};<-2.8mm,2.5mm>*{^1}**@{},
    \end{xy}\in  \fG^\circlearrowright(2,4)
$.

\vspace{-2mm}

\subsection{Decorated directed graphs}
Let $E$ be an $\bS$-{\em bimodule}, that is, a family, $\{E(p,q)\}_{p,q\geq 0}$, of vector spaces on
which the group
$\bS_p$ acts on the left and the group $\bS_q$ acts on the right, and both actions commute with each other.
We shall use elements of $E$ to decorate vertices of an arbitrary graph
$G\in \fG^\circlearrowright$  as follows. First, for each vertex $v\in V(G)$  we construct
a vector space
$$
E({\mathit O\mathit u\mathit t}_v, In_v):=  \langle {\mathit O\mathit u\mathit t}_v
\rangle \ot_{\bS_p} E(p,q) \ot_{\bS_q} \langle In_v\rangle,
$$
where  $\langle {\mathit O\mathit u\mathit t}_v\rangle$ (resp.,
$\langle In_v\rangle$) is the vector space spanned by all bijections
$[\#   {\mathit O\mathit u\mathit t}_v]\rar {\mathit O\mathit u\mathit t}_v$
(resp., $In_v\rar [\# In_v])$.
It is (non-canonically) isomorphic to $E(p,q)$ as a vector space and carries  natural actions
of the automorphism groups of the sets ${\mathit O\mathit u\mathit t}_v$ and $In_v$. These
actions  make the following
{\em unordered tensor product}\, over the set $V(G)$ (of cardinality, say, $k$),
$$
\bigotimes_{v\in V(G)} E(Out_v, In_v):= \left(\bigoplus_{i:[k]\rar V(G) }
 E({\mathit O\mathit u\mathit t}_{i(1)}, In_{i(1)})
\ot\ldots \ot
 E({\mathit O\mathit u\mathit t}_{i(k)}, In_{i(k)})\right)_{\bS_k},
$$
into a representation space of the automorphism group, $Aut(G)$, of the graph $G$ which,
by definition, is the subgroup of the
 symmetry group of the 1-dimensional
$CW$-complex underlying the graph $G$ which fixes its legs. Hence with an arbitrary
graph $G\in \fG^\circlearrowright$
and an arbitrary $\bS$-bimodule $E$ one can associate a vector space,
$$
G\langle E\rangle:= \left(
\otimes_{v\in V(G)} E(Out_v, In_v)\right)_{Aut G},
$$
whose elements are called {\em decorated (by $E$) graphs}. For example, the automorphism
group of the graph
$G=\Ba{c}
\begin{xy}
 <0mm,0mm>*{\bullet};
<0mm,0.41mm>*{};<0mm,2.9mm>*{}**@{-},
 <0.39mm,-0.39mm>*{};<2.4mm,-2.4mm>*{}**@{-},
 <-0.35mm,-0.35mm>*{};<-2.4mm,-2.4mm>*{}**@{-},
 <-2.4mm,-2.4mm>*{};<-0.4mm,-4.5mm>*{}**@{-},
 <2.4mm,-2.4mm>*{};<0.4mm,-4.5mm>*{}**@{-},
  <0mm,-5.1mm>*{\bullet};
  <0.4mm,-5.5mm>*{};<2mm,-7.7mm>*{}**@{-},
<-0.4mm,-5.5mm>*{};<-2mm,-7.7mm>*{}**@{-},
  <0.4mm,-5.5mm>*{};<2.9mm,-9.7mm>*{^2}**@{},
 <0.4mm,-5.5mm>*{};<-2.9mm,-9.7mm>*{^1}**@{},
 \end{xy}
\Ea$
is $\Z_2$ so that $G\langle E \rangle= E(1,2)\ot_{\Z_2} E(2,2)$. It is useful to think
of an element in $G\langle E\rangle$  as of the graph $G$
whose vertices are literarily decorated by some
elements $a\in E(1,2)$ and $b\in E(2,1)$ and subject to the following
relations,
$$
\Ba{c}
\begin{xy}
 <0mm,0mm>*{\bullet};
<2.5mm,0mm>*{^a};
<0mm,0.41mm>*{};<0mm,2.9mm>*{}**@{-},
 <0.39mm,-0.39mm>*{};<2.4mm,-2.4mm>*{}**@{-},
 <-0.35mm,-0.35mm>*{};<-2.4mm,-2.4mm>*{}**@{-},
 <-2.4mm,-2.4mm>*{};<-0.4mm,-4.5mm>*{}**@{-},
 <2.4mm,-2.4mm>*{};<0.4mm,-4.5mm>*{}**@{-},
  <0mm,-5.1mm>*{\bullet};
 <2.5mm,-5.1mm>*{_b};
  <0.4mm,-5.5mm>*{};<2mm,-7.7mm>*{}**@{-},
<-0.4mm,-5.5mm>*{};<-2mm,-7.7mm>*{}**@{-},
  <0.4mm,-5.5mm>*{};<2.9mm,-9.7mm>*{^2}**@{},
 <0.4mm,-5.5mm>*{};<-2.9mm,-9.7mm>*{^1}**@{},
 \end{xy}
\Ea
=
\Ba{c}
\begin{xy}
 <0mm,0mm>*{\bullet};
<5.2mm,0mm>*{^{a\sigma^{-1}}};
<0mm,0.41mm>*{};<0mm,2.9mm>*{}**@{-},
 <0.39mm,-0.39mm>*{};<2.4mm,-2.4mm>*{}**@{-},
 <-0.35mm,-0.35mm>*{};<-2.4mm,-2.4mm>*{}**@{-},
 <-2.4mm,-2.4mm>*{};<-0.4mm,-4.5mm>*{}**@{-},
 <2.4mm,-2.4mm>*{};<0.4mm,-4.5mm>*{}**@{-},
  <0mm,-5.1mm>*{\bullet};
 <3mm,-5.1mm>*{_{\sigma b}};
  <0.4mm,-5.5mm>*{};<2mm,-7.7mm>*{}**@{-},
<-0.4mm,-5.5mm>*{};<-2mm,-7.7mm>*{}**@{-},
  <0.4mm,-5.5mm>*{};<2.9mm,-9.7mm>*{^2}**@{},
 <0.4mm,-5.5mm>*{};<-2.9mm,-9.7mm>*{^1}**@{},
 \end{xy}
 \Ea  \ \sigma\in \Z_2,
\ \ \ \
\lambda\left(
\Ba{c}
\begin{xy}
 <0mm,0mm>*{\bullet};
<3mm,0mm>*{a};
<0mm,0.41mm>*{};<0mm,2.9mm>*{}**@{-},
 <0.39mm,-0.39mm>*{};<2.4mm,-2.4mm>*{}**@{-},
 <-0.35mm,-0.35mm>*{};<-2.4mm,-2.4mm>*{}**@{-},
 <-2.4mm,-2.4mm>*{};<-0.4mm,-4.5mm>*{}**@{-},
 <2.4mm,-2.4mm>*{};<0.4mm,-4.5mm>*{}**@{-},
  <0mm,-5.1mm>*{\bullet};
 <3mm,-5.1mm>*{b};
  <0.4mm,-5.5mm>*{};<2mm,-7.7mm>*{}**@{-},
<-0.4mm,-5.5mm>*{};<-2mm,-7.7mm>*{}**@{-},
  <0.4mm,-5.5mm>*{};<2.9mm,-9.7mm>*{^2}**@{},
 <0.4mm,-5.5mm>*{};<-2.9mm,-9.7mm>*{^1}**@{},
 \end{xy}
\Ea
\right)
=
\Ba{c}
\begin{xy}
 <0mm,0mm>*{\bullet};
<5mm,0mm>*{\lambda a};
<0mm,0.41mm>*{};<0mm,2.9mm>*{}**@{-},
 <0.39mm,-0.39mm>*{};<2.4mm,-2.4mm>*{}**@{-},
 <-0.35mm,-0.35mm>*{};<-2.4mm,-2.4mm>*{}**@{-},
 <-2.4mm,-2.4mm>*{};<-0.4mm,-4.5mm>*{}**@{-},
 <2.4mm,-2.4mm>*{};<0.4mm,-4.5mm>*{}**@{-},
  <0mm,-5.1mm>*{\bullet};
 <3mm,-5.1mm>*{b};
  <0.4mm,-5.5mm>*{};<2mm,-7.7mm>*{}**@{-},
<-0.4mm,-5.5mm>*{};<-2mm,-7.7mm>*{}**@{-},
  <0.4mm,-5.5mm>*{};<2.9mm,-9.7mm>*{^2}**@{},
 <0.4mm,-5.5mm>*{};<-2.9mm,-9.7mm>*{^1}**@{},
 \end{xy}
\Ea
=
\Ba{c}
\begin{xy}
 <0mm,0mm>*{\bullet};
<3mm,0mm>*{a};
<0mm,0.41mm>*{};<0mm,2.9mm>*{}**@{-},
 <0.39mm,-0.39mm>*{};<2.4mm,-2.4mm>*{}**@{-},
 <-0.35mm,-0.35mm>*{};<-2.4mm,-2.4mm>*{}**@{-},
 <-2.4mm,-2.4mm>*{};<-0.4mm,-4.5mm>*{}**@{-},
 <2.4mm,-2.4mm>*{};<0.4mm,-4.5mm>*{}**@{-},
  <0mm,-5.1mm>*{\bullet};
 <5mm,-5.1mm>*{\lambda b};
  <0.4mm,-5.5mm>*{};<2mm,-7.7mm>*{}**@{-},
<-0.4mm,-5.5mm>*{};<-2mm,-7.7mm>*{}**@{-},
  <0.4mm,-5.5mm>*{};<2.9mm,-9.7mm>*{^2}**@{},
 <0.4mm,-5.5mm>*{};<-2.9mm,-9.7mm>*{^1}**@{},
 \end{xy}
\Ea
\ \ \ \forall \lambda \in \K,
$$
$$
\hspace{-40mm}
\Ba{c}
\begin{xy}
 <0mm,0mm>*{\bullet};
<7mm,0mm>*{a_1\hspace{-0.5mm}+\hspace{-0.5mm} a_2};
<0mm,0.41mm>*{};<0mm,2.9mm>*{}**@{-},
 <0.39mm,-0.39mm>*{};<2.4mm,-2.4mm>*{}**@{-},
 <-0.35mm,-0.35mm>*{};<-2.4mm,-2.4mm>*{}**@{-},
 <-2.4mm,-2.4mm>*{};<-0.4mm,-4.5mm>*{}**@{-},
 <2.4mm,-2.4mm>*{};<0.4mm,-4.5mm>*{}**@{-},
  <0mm,-5.1mm>*{\bullet};
 <3mm,-5.1mm>*{b};
  <0.4mm,-5.5mm>*{};<2mm,-7.7mm>*{}**@{-},
<-0.4mm,-5.5mm>*{};<-2mm,-7.7mm>*{}**@{-},
  <0.4mm,-5.5mm>*{};<2.9mm,-9.7mm>*{^2}**@{},
 <0.4mm,-5.5mm>*{};<-2.9mm,-9.7mm>*{^1}**@{},
 \end{xy}
\Ea=
\Ba{c}
\begin{xy}
 <0mm,0mm>*{\bullet};
<3mm,0mm>*{a_1};
<0mm,0.41mm>*{};<0mm,2.9mm>*{}**@{-},
 <0.39mm,-0.39mm>*{};<2.4mm,-2.4mm>*{}**@{-},
 <-0.35mm,-0.35mm>*{};<-2.4mm,-2.4mm>*{}**@{-},
 <-2.4mm,-2.4mm>*{};<-0.4mm,-4.5mm>*{}**@{-},
 <2.4mm,-2.4mm>*{};<0.4mm,-4.5mm>*{}**@{-},
  <0mm,-5.1mm>*{\bullet};
 <3mm,-5.1mm>*{b};
  <0.4mm,-5.5mm>*{};<2mm,-7.7mm>*{}**@{-},
<-0.4mm,-5.5mm>*{};<-2mm,-7.7mm>*{}**@{-},
  <0.4mm,-5.5mm>*{};<2.9mm,-9.7mm>*{^2}**@{},
 <0.4mm,-5.5mm>*{};<-2.9mm,-9.7mm>*{^1}**@{},
 \end{xy}
\Ea
+
\Ba{c}
\begin{xy}
 <0mm,0mm>*{\bullet};
<3mm,0mm>*{a_2};
<0mm,0.41mm>*{};<0mm,2.9mm>*{}**@{-},
 <0.39mm,-0.39mm>*{};<2.4mm,-2.4mm>*{}**@{-},
 <-0.35mm,-0.35mm>*{};<-2.4mm,-2.4mm>*{}**@{-},
 <-2.4mm,-2.4mm>*{};<-0.4mm,-4.5mm>*{}**@{-},
 <2.4mm,-2.4mm>*{};<0.4mm,-4.5mm>*{}**@{-},
  <0mm,-5.1mm>*{\bullet};
 <3mm,-5.1mm>*{b};
  <0.4mm,-5.5mm>*{};<2mm,-7.7mm>*{}**@{-},
<-0.4mm,-5.5mm>*{};<-2mm,-7.7mm>*{}**@{-},
  <0.4mm,-5.5mm>*{};<2.9mm,-9.7mm>*{^2}**@{},
 <0.4mm,-5.5mm>*{};<-2.9mm,-9.7mm>*{^1}**@{},
 \end{xy}
\Ea
\ \ \ \mbox{and similarly for $b$}.
$$
It also follows from the definition that
$
\Ba{c}
\begin{xy}
 <0mm,0mm>*{\bullet};
<2.5mm,0mm>*{^a};
<0mm,0.41mm>*{};<0mm,2.9mm>*{}**@{-},
 <0.39mm,-0.39mm>*{};<2.4mm,-2.4mm>*{}**@{-},
 <-0.35mm,-0.35mm>*{};<-2.4mm,-2.4mm>*{}**@{-},
 <-2.4mm,-2.4mm>*{};<-0.4mm,-4.5mm>*{}**@{-},
 <2.4mm,-2.4mm>*{};<0.4mm,-4.5mm>*{}**@{-},
  <0mm,-5.1mm>*{\bullet};
 <2.5mm,-5.1mm>*{_b};
  <0.4mm,-5.5mm>*{};<2mm,-7.7mm>*{}**@{-},
<-0.4mm,-5.5mm>*{};<-2mm,-7.7mm>*{}**@{-},
  <0.4mm,-5.5mm>*{};<2.9mm,-9.7mm>*{^2}**@{},
 <0.4mm,-5.5mm>*{};<-2.9mm,-9.7mm>*{^1}**@{},
 \end{xy}
\Ea =
\Ba{c}
\begin{xy}
 <0mm,0mm>*{\bullet};
<2.5mm,0mm>*{^a};
<0mm,0.41mm>*{};<0mm,2.9mm>*{}**@{-},
 <0.39mm,-0.39mm>*{};<2.4mm,-2.4mm>*{}**@{-},
 <-0.35mm,-0.35mm>*{};<-2.4mm,-2.4mm>*{}**@{-},
 <-2.4mm,-2.4mm>*{};<-0.4mm,-4.5mm>*{}**@{-},
 <2.4mm,-2.4mm>*{};<0.4mm,-4.5mm>*{}**@{-},
  <0mm,-5.1mm>*{\bullet};
 <5mm,-5.1mm>*{_{b(12)}};
  <0.4mm,-5.5mm>*{};<2mm,-7.7mm>*{}**@{-},
<-0.4mm,-5.5mm>*{};<-2mm,-7.7mm>*{}**@{-},
  <0.4mm,-5.5mm>*{};<2.9mm,-9.7mm>*{^1}**@{},
 <0.4mm,-5.5mm>*{};<-2.9mm,-9.7mm>*{^2}**@{},
 \end{xy}
\Ea$,  $(12)\in \Z_2$.

\subsubsection{Remark} If $E=\{E(p,q)\}$ is a {\em differential
graded}\, (dg, for short) $\bS$-bimodule, i.e. if each vector space $E(p,q)$ is a complex equipped with an $\bS_p\times \bS_q$-equivariant
differential $\delta$, then, for any graph $G\in \fG^\circlearrowright(m,n)$,
the associated  graded
vector space
$G\langle E \rangle$ comes equipped with an induced  $\bS_m\times \bS_n$-equivariant differential $\delta_G$ so
that the collection, $\{\bigoplus_{G\in \fG^\circlearrowright(m,n)} G\langle E \rangle\}_{m,n\geq 0}$,
 is again  a {\em dg}\, $\bS$-bimodule. We sometimes abbreviate $\delta_G$ to $\delta$.

\subsubsection{Remark}\label{1-corolla} The one vertex graph
$
{\mathfrak C}_{m,n}:=
\Ba{c}
\begin{xy}
 <0mm,0mm>*{\bullet};
 <-0.5mm,0.2mm>*{};<-8mm,3mm>*{}**@{-},
 <-0.4mm,0.3mm>*{};<-4.5mm,3mm>*{}**@{-},
 <0mm,0mm>*{};<0mm,2.6mm>*{\ldots}**@{},
 <0.4mm,0.3mm>*{};<4.5mm,3mm>*{}**@{-},
 <0.5mm,0.2mm>*{};<8mm,3mm>*{}**@{-},
<-0.4mm,-0.2mm>*{};<-8mm,-3mm>*{}**@{-},
 <-0.5mm,-0.3mm>*{};<-4.5mm,-3mm>*{}**@{-},
 <0mm,0mm>*{};<0mm,-2.6mm>*{\ldots}**@{},
 <0.5mm,-0.3mm>*{};<4.5mm,-3mm>*{}**@{-},
 <0.4mm,-0.2mm>*{};<8mm,-3mm>*{}**@{-};
<0mm,5mm>*{\overbrace{\ \ \ \ \ \ \ \ \ \ \ \ \ \  }};
<0mm,-5mm>*{\underbrace{\ \ \ \ \ \ \ \ \ \ \ \ \ \ }};
<0mm,7mm>*{^{m\ \ output\ legs}};
<0mm,-7mm>*{_{n\ \ input\ legs}};
 \end{xy}\in \fG^\circlearrowright(m,n)\Ea
$
is often called {\em the $(m,n)$-corolla}. It is clear  that for any $\bS$-bimodule $E$ one has
$G\langle E \rangle = E(m,n)$.

\subsection{Wheeled props}\label{1: subsect Wheeled props}
A {\em wheeled prop}\,  is an $\bS$-bimodule $\cP=\{\cP(m,n)\}$ together with
a family of linear $\bS_m\times \bS_n$-equivariant maps,
$$
\left\{\mu_G: G\langle \cP\rangle\rar \cP(m,n)\right\}_{G\in \fG^\circlearrowright(m,n)},\ \ m,n\geq 0,
$$
parameterized by elements $G\in \fG^\circlearrowright$, which
satisfy a ``three-dimensional" associativity condition,
\Beq\label{graph-associativity}
\mu_G=\mu_{G/H}\circ \mu_H',
\Eeq
 for any subgraph $H\subset G$. Here $G/H$ is the graph obtained from $G$ by shrinking
 the whole subgraph $H$ into a single internal vertex, and
  $\mu_H': G\langle E \rangle \rar (G/H)\langle E\rangle$ stands for the map
which equals $\mu_H$ on the decorated vertices lying in $H$ and which is identity on all other vertices of $G$.

\sip

If the $\bS$-bimodule $\cP$ underlying a wheeled prop has a differential $\delta$ satisfying,
for any $G\in \fG^\circlearrowright$, the condition $\delta\circ \mu_G=\mu_G\circ \delta_G$, then the
wheeled prop $\cP$
is called {\em differential}.

\mip

By Remark~\ref{1-corolla}, the values of the maps $\mu_G$ can be identified with decorated corollas, and hence
the maps themselves can  be visually understood as  {\em contraction}\, maps,
$\mu_{G\in \fG^\circlearrowright(m,n)}: G\langle \cP\rangle\rar {\mathfrak C}_{m,n}\langle \cP \rangle$,
contracting all the  edges and vertices of $G$ into a single vertex.

\subsubsection{Remark}
 Strictly speaking, the notion introduced in \S~\ref{1: subsect Wheeled props} should be called a wheeled
prop {\em without unit}. A wheeled prop {\em with unit}\, can be defined as in \S 2.1.1
provided one  enlarges
$\fG^\circlearrowright$ by adding a family of graphs,
$\{
\uparrow \ \uparrow  \cdots \uparrow
\hgw\hgw\cdots \hgw\}
$,  {\em without vertices}
 \cite{MMS}.

\subsection{Props, properads, operads, etc.\ as $\fG$-algebras} Let $\fG=\sqcup_{m,n}\fG(m,n)$
be a subset of the set $\fG^\circlearrowright$, say, one of the subsets defined in Table 1 below.
 A subgraph $H$ of a graph $G\in \fG$ is called
{\em admissible}\, if $H\in \fG$ and $G/H\in \fG$, i.e.\ a contraction of a graph from $\fG$ by a subgraph belonging to $\fG$
gives a graph which again belongs to $\fG$.
\sip

A $\fG$-{\em algebra}\,  is, by definition (cf.\ \S~\ref{1: subsect Wheeled props}),
an $\bS$-bimodule $\cP=\{\cP(m,n)\}$ together with
a family of linear $\bS_m\times \bS_n$-equivariant maps,
$\left\{\mu_G: G\langle \cP\rangle\rar \cP(m,n)\right\}_{G\in \fG^\circlearrowright(m,n)},$
parameterized by elements $G\in \fG$, which
satisfy condition (\ref{graph-associativity}) for any admissible subgraph $H\subset G$.
Applying this idea to the subfamilies
$\fG\subset \fG^\circlearrowright$ from Table 1 gives us, in the chronological order,
the notions of {\em prop, operad, dioperad, properad,
$\frac{1}{2}$-prop}\,
and their {\em wheeled}\, versions which have been introduced, respectively,
in the papers \cite{Mc, May, G, V, Ko0, Me3, MMS}.
\sip

 We leave it to the reader as an exercise to check that $\fG^\mid$-algebra structures on
an $\bS$-bimodule $E$ with only $E(1,1)$ non-zero are precisely
 associative algebra structures on $E(1,1)$. This fact implies that, for any $\fG$-algebra
 $E=\{E(m,n)\}_{m,n\geq 0}$,
 the space $E(1,1)$ is an associative algebra.

\subsection{Basic examples of $\fG$-algebras}
{ ({\it i}\,)} For any $\fG$ and any finite-dimensional vector space $V$ the $\bS$-bimodule
$\cE nd_V=\{ \Hom(V^{\ot n}, V^{\ot m})\}$ is naturally a
 $\fG$-algebra with contraction maps $\mu_{G\in \fG}$ being ordinary compositions and, possibly, traces of
 linear maps; it is
 called the {\em endomorphism
$\fG$-algebra
of $V$}. In the cases $\fG\neq \fG^\circlearrowright,\fG^\circlearrowright_c$ the assumption of finite-dimensionality
of $V$ can be dropped (as the defining operations $\mu_G$ do not employ traces).

\sip

({\it ii}\,) With any  $\bS$-bimodule,
 $E=\{E(m,n)\}$, there is associated  another $\bS$-bimodule,
$\cF^\fG\langle E\rangle=\{\cF^\fG\langle E\rangle (m,n)\}$ with
$
\cF^\fG\langle E\rangle (m,n):= \bigoplus_{G\in \fG(m,n)} G\langle E\rangle$,
which has a natural $\fG$-algebra structure with the contraction maps
$\mu_G$ being tautological. The $\fG$-algebra  $\cF^\fG \langle E\rangle$ is called {\em the free $\fG$-algebra
generated
by the $\bS$-bimodule $E$}. We often abbreviate notations by replacing $\cF^{\fG^\circlearrowright}$
by $\cF^\circlearrowright$, $\cF^{\fG^\curlywedge}$ by $\cF^\curlywedge$, etc.

\sip

({\it iii}\,) Definitions of $\fG$-{\em sub}algebras, $\mathcal Q\subset \cP$, of $\fG$-algebras, of their
ideals, $\cI\subset \cP$, and the associated quotient $\fG$-algebras, $\cP/\cI$, are
straightforward. We omit the details.

\begin{center}
\begin{tabular}{|c|c||c|c|c|}
\multicolumn{4}{c}{{\bf Table 1}\vspace{2 mm}: A list of $\fG$-algebras}\\
\hline
&&&\vspace{-2mm}\\
 $\fG$ & {\bf Definition}
  &   $\fG${\bf-algebra}  & {\bf Typical examples} \\
  &&& \vspace{-2mm} \\
\hline
\hline
&&&\vspace{-6mm} \\
 $\fG^\circlearrowright$ & All possible directed graphs  &
  $\Ba{c}\mbox{Wheeled}\\\mbox{prop}\Ea $ &
  \begin{xy}
 <0.4mm,0.0mm>*{};<2.4mm,2.1mm>*{}**@{-},
 <-0.38mm,-0.2mm>*{};<-2.8mm,2.5mm>*{}**@{-},
<0mm,-0.8mm>*{\bullet};
<0mm,-1.0mm>*{};<0mm,-3.6mm>*{}**@{-},
 <2.96mm,2.4mm>*{\bullet};
 <2.4mm,2.8mm>*{};<0mm,5mm>*{}**@{-},
  <3.4mm,3.1mm>*{};<5.1mm,5mm>*{}**@{-},
%
<0mm,5mm>*{\bullet}**@{},
<-2.8mm,2.5mm>*{};<0mm,5mm>*{}**@{-},
<0mm,5mm>*{};<0mm,8.6mm>*{}**@{-},
\end{xy}
\ \
\begin{xy}
 <0mm,2.47mm>*{};<0mm,-0.5mm>*{}**@{-},
 <0.5mm,3.5mm>*{};<2.2mm,5.2mm>*{}**@{-},
 <-0.48mm,3.48mm>*{};<-2.2mm,5.2mm>*{}**@{-},
 <0mm,3mm>*{\bullet};
  <0mm,-0.8mm>*{\bullet};
<0mm,-0.8mm>*{};<-2.2mm,-3.5mm>*{}**@{-},
<0mm,-0.8mm>*{};<2.2mm,-3.5mm>*{}**@{-},
 <-2.5mm,5.7mm>*{\bullet};
<-2.5mm,5.7mm>*{};<-2.5mm,8.4mm>*{}**@{-},
<-2.5mm,5.7mm>*{};<-5mm,3mm>*{}**@{-},
<-5mm,3mm>*{};<-5mm,-0.8mm>*{}**@{-},
 <-2.5mm,-4.2mm>*{\bullet};
 <-2.8mm,-3.6mm>*{};<-5mm,-0.8mm>*{}**@{-},
 <-2.5mm,-4.6mm>*{};<-2.5mm,-6.3mm>*{}**@{-},
 (-2.2,8.4)*{}
   \ar@{->}@(ur,dr) (-2.2,-6.4)*{}
\end{xy}
\vspace{-6mm}
 \\
&&& \\
\hline
&&&
\vspace{-6mm}
 \\
 $\fG^\circlearrowright_c$ &$\Ba{c} \mbox{A subset}\ \fG^\circlearrowright_c\subset \fG^\circlearrowright
 \ \mbox{consisting}
\\ \mbox{ of all {\em connected}\, graphs}\Ea $ &
  $\Ba{c}\mbox{Wheeled}\\\mbox{properad}\Ea $ &
\begin{xy}
 <0mm,2.47mm>*{};<0mm,-0.5mm>*{}**@{-},
 <0.5mm,3.5mm>*{};<2.2mm,5.2mm>*{}**@{-},
 <-0.48mm,3.48mm>*{};<-2.2mm,5.2mm>*{}**@{-},
 <0mm,3mm>*{\bullet};
  <0mm,-0.8mm>*{\bullet};
<0mm,-0.8mm>*{};<-2.2mm,-3.5mm>*{}**@{-},
<0mm,-0.8mm>*{};<2.2mm,-3.5mm>*{}**@{-},
 <-2.5mm,5.7mm>*{\bullet};
<-2.5mm,5.7mm>*{};<-2.5mm,8.4mm>*{}**@{-},
<-2.5mm,5.7mm>*{};<-5mm,3mm>*{}**@{-},
<-5mm,3mm>*{};<-5mm,-0.8mm>*{}**@{-},
 <-2.5mm,-4.2mm>*{\bullet};
 <-2.8mm,-3.6mm>*{};<-5mm,-0.8mm>*{}**@{-},
 <-2.5mm,-4.6mm>*{};<-2.5mm,-6.3mm>*{}**@{-},
 (-2.2,8.4)*{}
   \ar@{->}@(ur,dr) (-2.2,-6.4)*{}
\end{xy}
\vspace{-3mm} \\
&&& \vspace{-3mm}\\
\hline
&&& \vspace{-5mm}\\
$\fG^\circlearrowright_{oper}$ & $\Ba{c} \mbox{A subset}\ \fG^\circlearrowright_{oper}\subset
\fG^\circlearrowright_{c}\
\mbox{consisting of graphs}
\\
 \mbox{ whose vertices have at most {one} output leg
}\Ea $ &
  $\Ba{c}\mbox{Wheeled}\\\mbox{operad}\Ea $ &
   \begin{xy}
 <0mm,0mm>*{\bullet};<0mm,0mm>*{}**@{},
 <0mm,0.69mm>*{};<0mm,4.0mm>*{}**@{-},
 <0.39mm,-0.39mm>*{};<2.4mm,-2.4mm>*{}**@{-},
 <-0.35mm,-0.35mm>*{};<-1.9mm,-1.9mm>*{}**@{-},
 <-2.4mm,-2.4mm>*{\bullet};<-2.4mm,-2.4mm>*{}**@{},
 <-2.0mm,-2.8mm>*{};<0mm,-4.9mm>*{}**@{-},
 <-2.8mm,-2.9mm>*{};<-4.7mm,-4.9mm>*{}**@{-},
  (0.0,4.0)*{}
   \ar@{->}@(ur,dr) (0.0,-4.9)*{}
 \end{xy}
  \\
&&& \vspace{-5mm}\\
\hline
&&& \vspace{-3mm}\\

$\fG^\uparrow$ & $\Ba{c} \mbox{A subset}\ \fG^\uparrow\subset \fG^\circlearrowright
\ \mbox{consisting}
\\ \mbox{of graphs with no {\em wheels}, i.e.}\\
\mbox{with no directed closed paths of edges}\Ea $ &
  $\Ba{c}\mbox{Prop}\Ea $ &
    \begin{xy}
 <0.4mm,0.0mm>*{};<2.4mm,2.1mm>*{}**@{-},
 <-0.38mm,-0.2mm>*{};<-2.8mm,2.5mm>*{}**@{-},
<0mm,-0.8mm>*{\bullet};
<0mm,-1.0mm>*{};<0mm,-3.6mm>*{}**@{-},
 <2.96mm,2.4mm>*{\bullet};
 <2.4mm,2.8mm>*{};<0mm,5mm>*{}**@{-},
  <3.4mm,3.1mm>*{};<5.1mm,5mm>*{}**@{-},
%
<-2.8mm,2.5mm>*{};<0mm,5mm>*{\bullet}**@{},
<-2.8mm,2.5mm>*{};<0mm,5mm>*{}**@{-},
<0mm,5mm>*{};<0mm,8.6mm>*{}**@{-},
\end{xy}
  \\
&&& \vspace{-3mm}\\
\hline
&&&\vspace{-3mm} \\
$\fG^\uparrow_c$ & $\Ba{c} \mbox{A subset}\ \fG^\uparrow_c\subset \fG^\uparrow\
\mbox{consisting}
\\ \mbox{ of all {\em connected}\, graphs}\Ea $ &
  $\Ba{c}\mbox{Properad}\Ea $ &
   \begin{xy}
 <0.4mm,0.0mm>*{};<2.4mm,2.1mm>*{}**@{-},
 <-0.38mm,-0.2mm>*{};<-2.8mm,2.5mm>*{}**@{-},
<0mm,-0.8mm>*{\bullet};
<0mm,-1.0mm>*{};<0mm,-3.6mm>*{}**@{-},
 <2.96mm,2.4mm>*{\bullet};
 <2.4mm,2.8mm>*{};<0mm,5mm>*{}**@{-},
  <3.4mm,3.1mm>*{};<5.1mm,5mm>*{}**@{-},
%
<-2.8mm,2.5mm>*{};<0mm,5mm>*{\bullet}**@{},
<-2.8mm,2.5mm>*{};<0mm,5mm>*{}**@{-},
<0mm,5mm>*{};<0mm,8.6mm>*{}**@{-},
\end{xy}
  \\
&&& \vspace{-3mm}\\
\hline
&&& \vspace{-3mm}\\
$\fG^\uparrow_{c,0}$ & $\Ba{c} \mbox{A subset}\ \fG^\uparrow_{c,0}\subset \fG^\uparrow_c\
\mbox{consisting}
\\ \mbox{ of graphs of genus zero
}\Ea $ &
  $\Ba{c}\mbox{Dioperad}\Ea $ &
\begin{xy}
 <0mm,2.47mm>*{};<0mm,-0.5mm>*{}**@{-},
 <0.5mm,3.5mm>*{};<2.2mm,5.2mm>*{}**@{-},
 <-0.48mm,3.48mm>*{};<-2.2mm,5.2mm>*{}**@{-},
 <0mm,3mm>*{\bullet};<0mm,3mm>*{}**@{},
  <0mm,-0.8mm>*{\bullet};<0mm,-0.8mm>*{}**@{},
<0mm,-0.8mm>*{};<-2.2mm,-3.5mm>*{}**@{-},
<0mm,-0.8mm>*{};<2.2mm,-3.5mm>*{}**@{-},
 <-2.5mm,5.7mm>*{\bullet};<0mm,0mm>*{}**@{},
<-2.5mm,5.7mm>*{};<-2.5mm,9.4mm>*{}**@{-},
<-2.5mm,5.7mm>*{};<-5mm,3mm>*{}**@{-},
%
 <-2.5mm,-4.2mm>*{\bullet};<0mm,3mm>*{}**@{},
 <-2.8mm,-3.6mm>*{};<-5mm,-0.8mm>*{}**@{-},
 <-2.5mm,-4.6mm>*{};<-2.5mm,-7.3mm>*{}**@{-},
\end{xy}
  \\
&&& \vspace{-3mm}\\
\hline
&&& \vspace{-3mm}\\
$\fG^{\frac{1}{2}}$ & $\Ba{c} \mbox{A subset}\ \fG^{\frac{1}{2}}\subset\fG^\uparrow_{c,0}\
\mbox{consisting of}
\\ \mbox{all $(m,n)$-graphs with the number }\\ \mbox{of directed paths from input legs}\\
 \mbox{to the output legs equal to $mn$}\\
\Ea $ &
  $\Ba{c}\mbox{$\frac{1}{2}$-Prop}\Ea $ &
   \begin{xy}
 <0mm,0mm>*{\bullet};
 <0mm,0.69mm>*{};<0mm,3.5mm>*{}**@{-},
 <0.39mm,-0.39mm>*{};<2.4mm,-2.4mm>*{}**@{-},
 <-0.35mm,-0.35mm>*{};<-1.9mm,-1.9mm>*{}**@{-},
 <-2.4mm,-2.4mm>*{\bullet};
 <-2.0mm,-2.8mm>*{};<0mm,-4.9mm>*{}**@{-},
 <-2.8mm,-2.9mm>*{};<-4.7mm,-4.9mm>*{}**@{-},
 <0mm,3.5mm>*{\bullet};
 <0mm,3.5mm>*{};<-2.4mm,5.9mm>*{}**@{-},
 <0mm,3.5mm>*{};<2.4mm,5.9mm>*{}**@{-},
 \end{xy}
  \\
&&& \vspace{-3mm}\\
\hline
&&& \vspace{-3mm}\\
$\fG^\curlywedge$ & $\Ba{c} \mbox{A subset}\ \fG^\curlywedge\subset \fG^\uparrow_{c,0}\
\mbox{consisting of graphs}
\\
 \mbox{ whose vertices have precisely {one} output leg
}\Ea $ &
  $\Ba{c}\mbox{Operad}\Ea $ &

   \begin{xy}
 <0mm,0mm>*{\bullet};<0mm,0mm>*{}**@{},
 <0mm,0.69mm>*{};<0mm,4.0mm>*{}**@{-},
 <0.39mm,-0.39mm>*{};<2.4mm,-2.4mm>*{}**@{-},
 <-0.35mm,-0.35mm>*{};<-1.9mm,-1.9mm>*{}**@{-},
 <-2.4mm,-2.4mm>*{\bullet};<-2.4mm,-2.4mm>*{}**@{},
 <-2.0mm,-2.8mm>*{};<0mm,-4.9mm>*{}**@{-},
 <-2.8mm,-2.9mm>*{};<-4.7mm,-4.9mm>*{}**@{-},
 \end{xy}
  \\
&&& \vspace{-3mm}\\
\hline
&&& \vspace{-3mm}\\

$\fG^\mid$ & $\Ba{c} \mbox{A subset}\ \fG^\mid\subset \fG^\curlywedge
\\
 \mbox{consisting of graphs whose}\\
\mbox{vertices have precisely}\\ \mbox{ {one} input leg
}\Ea $ &
  $\Ba{c}\mbox{Associative}\\
  \mbox{algebra}\Ea $ &

   \begin{xy}
 <0mm,0mm>*{\bullet};
  <0mm,3.5mm>*{\bullet};
  <0mm,-3.5mm>*{\bullet};
 <0mm,-7mm>*{};<0mm,7.0mm>*{}**@{-},
 \end{xy}
\vspace{-3mm}  \\
&&&\\
\hline
\end{tabular}
\end{center}

\subsection{Morphisms and resolutions of $\fG$-algebras}\label{1: section morphisms and resolutions}
 A morphisms of $\fG$-algebras, $\rho: \cP_1\rar \cP_2$,
is a morphism of the underlying $\bS$-bimodules such that, for any graph $G$,
 one has $\rho\circ \mu_G= \mu_G\circ (\rho^{\ot G})$,
where $\rho^{\ot G}$ is a map, $G\langle \cP_1\rangle \rar G\langle \cP_2\rangle$, which changes
decorations of each vertex
 in $G$ in accordance with $\rho$.
A morphism of $\fG$-algebras, $\cP\rar \cE nd \langle V\rangle$, is called a {\em representation}\, of the
$\fG$-algebra $\cP$ in a graded vector space $V$.

\sip

 A {\em free resolution}\, of a dg $\fG$-algebra
$\cP$ is, by definition, a dg free $\fG$-algebra, $(\cF^\fG\langle E \rangle, \delta)$,
together with a morphism,
$\pi: (\cF\langle E \rangle, \delta) \rar \cP$, which induces a cohomology isomorphism.
If the differential $\delta$ in $\cF\langle \cE \rangle$ is
decomposable with respect to compositions $\mu_G$, then it is
 called a {\em minimal model}\, of $\cP$ and  is often denoted by
$\cP_\infty$.

 \section{Applications to algebra and geometry}
\subsection{The operad of associative algebras}
 Let $A_0=\{A_0(m,n)\}$ be an $\bS$-bimodule with all $A_0(m,n)=0$ except
$A_0(1,2):= \K[\bS_2]$. The associated free operad $\cF^\curlywedge\langle A_0\rangle$ can be identified with
the vector space spanned by all connected planar graphs of the form
$\Ba{c}\begin{xy}
 <0mm,0mm>*{\bullet};
 <0mm,0.69mm>*{};<0mm,3.0mm>*{}**@{-},
 <0.39mm,-0.39mm>*{};<2.4mm,-2.4mm>*{}**@{-},
 <-0.35mm,-0.35mm>*{};<-1.9mm,-1.9mm>*{}**@{-},
<-2.4mm,-2.4mm>*{\bullet};
 <-2.0mm,-2.8mm>*{};<0mm,-4.9mm>*{}**@{-},
 <-2.8mm,-2.9mm>*{};<-4.7mm,-4.9mm>*{}**@{-},
    <0.39mm,-0.39mm>*{};<3.3mm,-4.0mm>*{^3}**@{},
    <-2.0mm,-2.8mm>*{};<2.8mm,-9.6mm>*{^2}**@{},
    <-2.0mm,-2.8mm>*{};<-2.8mm,-9.6mm>*{^4}**@{},
    <-2.8mm,-2.9mm>*{};<-5.2mm,-6.7mm>*{^1}**@{},
 <0mm,-5mm>*{\bullet};
 <0mm,-5mm>*{};<2.4mm,-7.4mm>*{}**@{-},
 <0mm,-5mm>*{};<-2.4mm,-7.4mm>*{}**@{-},
 \end{xy}\Ea$. In particular,
 $
 \cF^\curlywedge\langle A_0\rangle(1,2)\cong A_0(1,2)\cong
\mbox{span}\langle
\begin{xy}
 <0mm,0.66mm>*{};<0mm,3mm>*{}**@{-},
 <0.39mm,-0.39mm>*{};<2.2mm,-2.2mm>*{}**@{-},
 <-0.35mm,-0.35mm>*{};<-2.2mm,-2.2mm>*{}**@{-},
 <0mm,0mm>*{\bullet};<0mm,0mm>*{}**@{},
   <0.39mm,-0.39mm>*{};<2.9mm,-4mm>*{^2}**@{},
   <-0.35mm,-0.35mm>*{};<-2.8mm,-4mm>*{^1}**@{},
\end{xy}
\,
,\,
\begin{xy}
 <0mm,0.66mm>*{};<0mm,3mm>*{}**@{-},
 <0.39mm,-0.39mm>*{};<2.2mm,-2.2mm>*{}**@{-},
 <-0.35mm,-0.35mm>*{};<-2.2mm,-2.2mm>*{}**@{-},
 <0mm,0mm>*{\bullet};
   <0.39mm,-0.39mm>*{};<2.9mm,-4mm>*{^1}**@{},
   <-0.35mm,-0.35mm>*{};<-2.8mm,-4mm>*{^2}**@{},
\end{xy}
\rangle
$. Let $I_0$ be an ideal of $\cF^\curlywedge\langle A_0\rangle$ generated by the following 6
 planar graphs,
\begin{equation}\label{2: assoc_relations}
\Ba{c}
\begin{xy}
 <0mm,0mm>*{\bullet};<0mm,0mm>*{}**@{},
 <0mm,0.69mm>*{};<0mm,3.0mm>*{}**@{-},
 <0.39mm,-0.39mm>*{};<2.4mm,-2.4mm>*{}**@{-},
 <-0.35mm,-0.35mm>*{};<-1.9mm,-1.9mm>*{}**@{-},
 <-2.4mm,-2.4mm>*{\bullet};<-2.4mm,-2.4mm>*{}**@{},
 <-2.0mm,-2.8mm>*{};<0mm,-4.9mm>*{}**@{-},
 <-2.8mm,-2.9mm>*{};<-4.7mm,-4.9mm>*{}**@{-},
    <0.39mm,-0.39mm>*{};<5mm,-4.0mm>*{^{\sigma(3)}}**@{},
    <-2.0mm,-2.8mm>*{};<1.5mm,-6.7mm>*{^{\sigma(2)}}**@{},
    <-2.8mm,-2.9mm>*{};<-6.2mm,-6.7mm>*{^{\sigma(1)}}**@{},
 \end{xy}
\ - \
 \begin{xy}
 <0mm,0mm>*{\bullet};<0mm,0mm>*{}**@{},
 <0mm,0.69mm>*{};<0mm,3.0mm>*{}**@{-},
 <0.39mm,-0.39mm>*{};<2.4mm,-2.4mm>*{}**@{-},
 <-0.35mm,-0.35mm>*{};<-1.9mm,-1.9mm>*{}**@{-},
 <2.4mm,-2.4mm>*{\bullet};<-2.4mm,-2.4mm>*{}**@{},
 <2.0mm,-2.8mm>*{};<0mm,-4.9mm>*{}**@{-},
 <2.8mm,-2.9mm>*{};<4.7mm,-4.9mm>*{}**@{-},
    <0.39mm,-0.39mm>*{};<-4mm,-4.0mm>*{^{\sigma(1)}}**@{},
    <-2.0mm,-2.8mm>*{};<-0.7mm,-6.7mm>*{^{\sigma(2)}}**@{},
    <-2.8mm,-2.9mm>*{};<7.2mm,-6.7mm>*{^{\sigma(3)}}**@{},
 \end{xy}
 \Ea
 \in \cF^\curlywedge\langle A_0\rangle(1,3) , \ \forall\ \sigma\in \bS_3.
\end{equation}

\subsubsection{\bf Claim}{\em
There is a 1-1 correspondence between representations,
$\rho:\cA ss \rar \cE nd_V$, of the quotient operad,
$\cA ss:= \cF^\curlywedge\langle A_0\rangle/\langle I_0\rangle$,
in a space $V$
and  associative algebra structures on $V$.
}

\begin{proof}
The values of $\rho$ on arbitrary (equivalence classes of) planar graphs is uniquely determined
by its value, $\rho(\begin{xy}
 <0mm,0.66mm>*{};<0mm,3mm>*{}**@{-},
 <0.39mm,-0.39mm>*{};<2.2mm,-2.2mm>*{}**@{-},
 <-0.35mm,-0.35mm>*{};<-2.2mm,-2.2mm>*{}**@{-},
 <0mm,0mm>*{\bullet};
   <0.39mm,-0.39mm>*{};<2.9mm,-4mm>*{^2}**@{},
   <-0.35mm,-0.35mm>*{};<-2.8mm,-4mm>*{^1}**@{},
\end{xy})\in \Hom(V^{\ot 2}, V)$, on one of the two generators. Denote this value by $\mu$. As $\rho$
sends any of the graphs (\ref{2: assoc_relations}) to zero, the multiplication in $V$
given by $\mu$ must be associative.
\end{proof}

\sip

Thus the operad $\cA ss$ can be called the {\em operad of associative algebras}.
What could be a (minimal)  free resolution of $\cA ss$? By definition in
\S~\ref{1: section morphisms and resolutions}, this must be a {\em free}\, operad,
$\cF^\curlywedge\langle A\rangle$, generated by some $\bS$-bimodule
${A}=\{{A}(1,n)\}_{n\geq 2}$ equipped with a differential $\delta$ and a projection $\pi:
\cF^\curlywedge\langle {A}\rangle\rar \cA ss$ inducing an isomorphism,
$H(\cF^\curlywedge\langle {A}\rangle,
\delta)=\cA ss$, at the cohomology level. The latter condition suggests that we can choose
${A}(1,2)$ to be identical to $A_0(1,2)$ and set a differential $\delta$ to satisfy
$\delta\begin{xy}
 <0mm,0.66mm>*{};<0mm,3mm>*{}**@{-},
 <0.39mm,-0.39mm>*{};<2.2mm,-2.2mm>*{}**@{-},
 <-0.35mm,-0.35mm>*{};<-2.2mm,-2.2mm>*{}**@{-},
 <0mm,0mm>*{\bullet};
   <0.39mm,-0.39mm>*{};<2.9mm,-4mm>*{^2}**@{},
   <-0.35mm,-0.35mm>*{};<-2.8mm,-4mm>*{^1}**@{},
\end{xy}=0$. Then the graphs (\ref{2: assoc_relations})  are  cocycles in
$\cF^\curlywedge\langle {A}\rangle(1,3)$. In view of the cohomology isomorphism
$\cF^\curlywedge\langle {A}\rangle\rar \cA ss$, we have to make them coboundaries, and hence
are forced
to introduce an $\bS_3$-module,
$$
{A}(1,3):= \K[\bS_3][1]=\mbox{span}\left\langle
\begin{xy}
 <0mm,-4mm>*{};<0mm,4mm>*{}**@{-},
 <0.39mm,-0.39mm>*{};<3.2mm,-4mm>*{}**@{-},
 <-0.35mm,-0.35mm>*{};<-3.2mm,-4mm>*{}**@{-},
<0mm,0mm>*{\bullet};
<6mm,-6.7mm>*{^{\sigma(3)}}**@{},
<-6mm,-6.7mm>*{^{\sigma(1)}}**@{},
<0mm,-6.7mm>*{^{\sigma(2)}}**@{},
\end{xy}
\right\rangle_{\sigma\in \bS_3},\vspace{-4mm}
$$
and set
\begin{equation}\label{2: delta_ass(3)}
\delta
\begin{xy}
 <0mm,-4mm>*{};<0mm,4mm>*{}**@{-},
 <0.39mm,-0.39mm>*{};<3.2mm,-4mm>*{}**@{-},
 <-0.35mm,-0.35mm>*{};<-3.2mm,-4mm>*{}**@{-},
<0mm,0mm>*{\bullet};
   <0.39mm,-0.39mm>*{};<4mm,-6.9mm>*{^3}**@{},
   <-0.35mm,-0.35mm>*{};<-4mm,-6.9mm>*{^1}**@{},
    <-0.35mm,-0.35mm>*{};<0mm,-6.9mm>*{^2}**@{},
\end{xy}=
\begin{xy}
 <0mm,0mm>*{\bullet};
 <0mm,0.69mm>*{};<0mm,3.0mm>*{}**@{-},
 <0.39mm,-0.39mm>*{};<2.4mm,-2.4mm>*{}**@{-},
 <-0.35mm,-0.35mm>*{};<-1.9mm,-1.9mm>*{}**@{-},
 <-2.4mm,-2.4mm>*{\bullet};
 <-2.0mm,-2.8mm>*{};<0mm,-4.9mm>*{}**@{-},
 <-2.8mm,-2.9mm>*{};<-4.7mm,-4.9mm>*{}**@{-},
    <0.39mm,-0.39mm>*{};<3.3mm,-4.0mm>*{^3}**@{},
    <-2.0mm,-2.8mm>*{};<0.5mm,-6.7mm>*{^2}**@{},
    <-2.8mm,-2.9mm>*{};<-5.2mm,-6.7mm>*{^1}**@{},
 \end{xy}
\ - \
 \begin{xy}
 <0mm,0mm>*{\bullet};
 <0mm,0.69mm>*{};<0mm,3.0mm>*{}**@{-},
 <0.39mm,-0.39mm>*{};<2.4mm,-2.4mm>*{}**@{-},
 <-0.35mm,-0.35mm>*{};<-1.9mm,-1.9mm>*{}**@{-},
 <2.4mm,-2.4mm>*{\bullet};
 <2.0mm,-2.8mm>*{};<0mm,-4.9mm>*{}**@{-},
 <2.8mm,-2.9mm>*{};<4.7mm,-4.9mm>*{}**@{-},
    <0.39mm,-0.39mm>*{};<-3mm,-4.0mm>*{^1}**@{},
    <-2.0mm,-2.8mm>*{};<0mm,-6.7mm>*{^2}**@{},
    <-2.8mm,-2.9mm>*{};<5.2mm,-6.7mm>*{^3}**@{},
 \end{xy}.
\end{equation}
We get in this way a well-defined dg {\em free}\, operad together with a well-defined epimorphism,
$
\left( \cF^\curlywedge(
\begin{xy}
 <0mm,0.66mm>*{};<0mm,3mm>*{}**@{-},
 <0.39mm,-0.39mm>*{};<2.2mm,-2.2mm>*{}**@{-},
 <-0.35mm,-0.35mm>*{};<-2.2mm,-2.2mm>*{}**@{-},
 <0mm,0mm>*{\bullet};
\end{xy}
\ , \,
\begin{xy}
 <0mm,-2.2mm>*{};<0mm,3mm>*{}**@{-},
 <0.39mm,-0.39mm>*{};<2.2mm,-2.2mm>*{}**@{-},
 <-0.35mm,-0.35mm>*{};<-2.2mm,-2.2mm>*{}**@{-},
 <0mm,0mm>*{\bullet};
\end{xy})\ , \  \delta
\right)
\lon (\cA ss, 0),
$
sending $(1,3)$-corollas to zero. However, this epimorphisms fails to be a quasi-isomorphism as
$
\delta
\underbrace{
\left(
\Ba{c}
\begin{xy}
 <0mm,0.66mm>*{};<0mm,3mm>*{}**@{-},
 <0.39mm,-0.39mm>*{};<2.2mm,-2.2mm>*{}**@{-},
 <-0.35mm,-0.35mm>*{};<-2.2mm,-2.2mm>*{}**@{-},
 <0mm,0mm>*{\bullet};
 <-2.2mm,-2.2mm>*{\bullet};
<-2.2mm,-2.2mm>*{};<-5.2mm,-6.2mm>*{}**@{-},
<-2.2mm,-2.2mm>*{};<-2.2mm,-6.2mm>*{}**@{-},
<-2.2mm,-2.2mm>*{};<0.8mm,-6.2mm>*{}**@{-},
   <0.39mm,-0.39mm>*{};<2.9mm,-4.1mm>*{^4}**@{},
   <-0.35mm,-0.35mm>*{};<-5.3mm,-8.5mm>*{^1}**@{},
<-0.35mm,-0.35mm>*{};<-2.2mm,-8.5mm>*{^2}**@{},
<-0.35mm,-0.35mm>*{};<0.9mm,-8.5mm>*{^3}**@{},
\end{xy}
+
\begin{xy}
 <0mm,0.66mm>*{};<0mm,3mm>*{}**@{-},
 <0.39mm,-0.39mm>*{};<2.2mm,-2.2mm>*{}**@{-},
 <-0.35mm,-0.35mm>*{};<-2.2mm,-2.2mm>*{}**@{-},
 <0mm,0mm>*{\bullet};
 <2.2mm,-2.2mm>*{\bullet};
<2.2mm,-2.2mm>*{};<5.2mm,-6.2mm>*{}**@{-},
<2.2mm,-2.2mm>*{};<2.2mm,-6.2mm>*{}**@{-},
<2.2mm,-2.2mm>*{};<-0.8mm,-6.2mm>*{}**@{-},
   <0.39mm,-0.39mm>*{};<-2.9mm,-4.1mm>*{^1}**@{},
   <-0.35mm,-0.35mm>*{};<5.3mm,-8.5mm>*{^4}**@{},
<-0.35mm,-0.35mm>*{};<2.2mm,-8.5mm>*{^3}**@{},
<-0.35mm,-0.35mm>*{};<-0.9mm,-8.5mm>*{^2}**@{},
\end{xy}
-
\begin{xy}
 <0mm,-4mm>*{};<0mm,4mm>*{}**@{-},
 <0.39mm,-0.39mm>*{};<3.2mm,-4mm>*{}**@{-},
 <-0.35mm,-0.35mm>*{};<-3.2mm,-4mm>*{}**@{-},
<0mm,0mm>*{\bullet};
<-3.2mm,-4mm>*{\bullet};
<-3.3mm,-4mm>*{};<-5mm,-7mm>*{}**@{-},
<-3.3mm,-4mm>*{};<-1.6mm,-7mm>*{}**@{-},
   <0.39mm,-0.39mm>*{};<4mm,-6.6mm>*{^4}**@{},
   <-0.35mm,-0.35mm>*{};<-5mm,-9mm>*{^1}**@{},
<-0.35mm,-0.35mm>*{};<-1.6mm,-9mm>*{^2}**@{},
    <-0.35mm,-0.35mm>*{};<0.6mm,-6.6mm>*{^3}**@{},
\end{xy}
+
\begin{xy}
 <0mm,-4mm>*{};<0mm,4mm>*{}**@{-},
 <0.39mm,-0.39mm>*{};<3.2mm,-4mm>*{}**@{-},
 <-0.35mm,-0.35mm>*{};<-3.2mm,-4mm>*{}**@{-},
<0mm,0mm>*{\bullet};
<0mm,-4mm>*{\bullet};
<0mm,-4mm>*{};<-1.7mm,-7mm>*{}**@{-},
<0mm,-4mm>*{};<1.7mm,-7mm>*{}**@{-},
   <0.39mm,-0.39mm>*{};<4mm,-6.6mm>*{^4}**@{},
   <-0.35mm,-0.35mm>*{};<-1.7mm,-9mm>*{^2}**@{},
<-0.35mm,-0.35mm>*{};<1.7mm,-9mm>*{^3}**@{},
    <-0.35mm,-0.35mm>*{};<-3.4mm,-6.6mm>*{^1}**@{},
\end{xy}
-
\begin{xy}
 <0mm,-4mm>*{};<0mm,4mm>*{}**@{-},
 <0.39mm,-0.39mm>*{};<3.2mm,-4mm>*{}**@{-},
 <-0.35mm,-0.35mm>*{};<-3.2mm,-4mm>*{}**@{-},
<0mm,0mm>*{\bullet};
<3.2mm,-4mm>*{\bullet};
<3.3mm,-4mm>*{};<5mm,-7mm>*{}**@{-},
<3.3mm,-4mm>*{};<1.6mm,-7mm>*{}**@{-},
   <0.39mm,-0.39mm>*{};<-4mm,-6.6mm>*{^1}**@{},
   <-0.35mm,-0.35mm>*{};<5mm,-9mm>*{^3}**@{},
<-0.35mm,-0.35mm>*{};<1.6mm,-9mm>*{^4}**@{},
    <-0.35mm,-0.35mm>*{};<-0.6mm,-6.6mm>*{^2}**@{},
\end{xy}
\Ea
\right)
}_{a\ nontrivial\ cohomology\ class\ in\
H^{-1}( \cF^\curlywedge(
\begin{xy}
 <0mm,0.66mm>*{};<0mm,2.2mm>*{}**@{-},
 <0.39mm,-0.39mm>*{};<2.2mm,-2.2mm>*{}**@{-},
 <-0.35mm,-0.35mm>*{};<-2.2mm,-2.2mm>*{}**@{-},
 <0mm,0mm>*{\bullet};<0mm,0mm>*{}**@{},
\end{xy}
\ , \,
\begin{xy}
 <0mm,-2.2mm>*{};<0mm,2.2mm>*{}**@{-},
 <0.39mm,-0.39mm>*{};<2.2mm,-2.2mm>*{}**@{-},
 <-0.35mm,-0.35mm>*{};<-2.2mm,-2.2mm>*{}**@{-},
 <0mm,0mm>*{\bullet};<0mm,0mm>*{}**@{},
\end{xy}\, )\ ,  \delta)}
=0
$.
To kill this cohomology  class we have to introduce a new generating $(1,4)$-corolla,
$\Ba{c}\begin{xy}
 <0mm,0mm>*{};<0mm,3mm>*{}**@{-},
 <0mm,0mm>*{};<4mm,-3mm>*{}**@{-},
 <0mm,0mm>*{};<1.5mm,-3mm>*{}**@{-},
 <0mm,0mm>*{};<-1.5mm,-3mm>*{}**@{-},
 <0mm,0mm>*{};<-4mm,-3mm>*{}**@{-},
<0mm,0mm>*{\bullet};
   <0.39mm,-0.39mm>*{};<5mm,-5.1mm>*{^4}**@{},
   <-0.35mm,-0.35mm>*{};<-5mm,-5.1mm>*{^1}**@{},
    <-0.35mm,-0.35mm>*{};<2mm,-5.1mm>*{^3}**@{},
<-0.35mm,-0.35mm>*{};<-2mm,-5.1mm>*{^2}**@{},
\end{xy}\Ea$,
of degree $-2$ and set the value of the differential on it to be equal to the underbraced expression above.
Again we get a well-defined dg free operad together with a natural homomorphism,
$
\left( \cF(
\begin{xy}
 <0mm,0.66mm>*{};<0mm,3mm>*{}**@{-},
 <0.39mm,-0.39mm>*{};<2.2mm,-2.2mm>*{}**@{-},
 <-0.35mm,-0.35mm>*{};<-2.2mm,-2.2mm>*{}**@{-},
 <0mm,0mm>*{\bullet};<0mm,0mm>*{}**@{},
\end{xy}
\ , \,
\begin{xy}
 <0mm,-2.2mm>*{};<0mm,3mm>*{}**@{-},
 <0.39mm,-0.39mm>*{};<2.2mm,-2.2mm>*{}**@{-},
 <-0.35mm,-0.35mm>*{};<-2.2mm,-2.2mm>*{}**@{-},
 <0mm,0mm>*{\bullet};<0mm,0mm>*{}**@{},
\end{xy}
\ , \,
\begin{xy}
 <0mm,0mm>*{};<0mm,3mm>*{}**@{-},
 <0mm,0mm>*{};<3mm,-2.2mm>*{}**@{-},
<0mm,0mm>*{};<1mm,-2.2mm>*{}**@{-},
<0mm,0mm>*{};<-1mm,-2.2mm>*{}**@{-},
 <0mm,0mm>*{};<-3mm,-2.2mm>*{}**@{-},
 <0mm,0mm>*{\bullet};
\end{xy}
),\delta
\right)
\rar \cA ss,
$
which, again, fails to be a quasi-isomorphism. To treat the new problem one has to introduce a new
generating corolla of degree $-3$ with 5 input legs and so on.

\sip
\subsubsection{\bf Theorem~\cite{St}}\label{2: Stasheff_theorem}
 {\em The minimal resolution of $\cA ss$ is a dg free operad,
$\cA ss_\infty:=(\cF^\curlywedge\langle {A} \rangle, \delta)$, generated by the $\bS$-bimodule
${A}=\{{A}(1,n)\}$,
\Beq\label{2: Ass infty generators}
{A}(1,n):=\K[\bS_n][n-2]= \mbox{\em span}
\left\langle
\Ba{c}
\begin{xy}
 <0mm,0mm>*{\bullet};
 <0mm,0mm>*{};<-8mm,-5mm>*{}**@{-},
 <0mm,0mm>*{};<-4.5mm,-5mm>*{}**@{-},
 <0mm,0mm>*{};<0mm,-4mm>*{\ldots}**@{},
 <0mm,0mm>*{};<4.5mm,-5mm>*{}**@{-},
 <0mm,0mm>*{};<8mm,-5mm>*{}**@{-},
   <0mm,0mm>*{};<-11mm,-7.9mm>*{^{\sigma(1)}}**@{},
   <0mm,0mm>*{};<-4mm,-7.9mm>*{^{\sigma(2)}}**@{},
   <0mm,0mm>*{};<10.0mm,-7.9mm>*{^{\sigma(n)}}**@{},
 <0mm,0mm>*{};<0mm,5mm>*{}**@{-},
 \end{xy}\Ea
\right\rangle_{\sigma\in \bS_n},
\Eeq
and with the differential given on the generators by
\begin{equation}\label{2: delta_ass(n)}
\delta
\begin{xy}
<0mm,0mm>*{\bullet},
<0mm,5mm>*{}**@{-},
<-5mm,-5mm>*{}**@{-},
<-2mm,-5mm>*{}**@{-},
<2mm,-5mm>*{}**@{-},
<5mm,-5mm>*{}**@{-},
<0mm,-7mm>*{_{\sigma(1)\ \ \ \ldots\ \ \ \sigma(n)}},
\end{xy}
=\sum_{k=0}^{n-2}\sum_{l=2}^{n-k}
(-1)^{k+l(n-k-l)+1}
\begin{xy}
<0mm,0mm>*{\bullet},
<0mm,5mm>*{}**@{-},
<4mm,-7mm>*{^{\sigma(1)\dots\sigma(k)\qquad\sigma(k+l+1)\dots\sigma(n)}},
<-14mm,-5mm>*{}**@{-},
<-6mm,-5mm>*{}**@{-},
<20mm,-5mm>*{}**@{-},
<8mm,-5mm>*{}**@{-},
<0mm,-5mm>*{}**@{-},
<0mm,-5mm>*{\bullet};
<-5mm,-10mm>*{}**@{-},
<-2mm,-10mm>*{}**@{-},
<2mm,-10mm>*{}**@{-},
<5mm,-10mm>*{}**@{-},
<0mm,-12mm>*{_{\sigma(k+1)\dots\sigma(k+l)}},
\end{xy}.
\end{equation}
}
\vspace{-2mm}
\subsubsection{\bf Definition}
Representations, $\cA ss_\infty\rar \cE nd_V$, of the dg operad $(\cA ss_\infty,\delta)$ in a
dg vector space $V$ are called {\em $A_\infty$-structures}\, in $V$.

\subsubsection{\bf Remark} We now suggest the reader to re-read Stasheff's Theorem~\ref{2: Stasheff_theorem}
from the end to the beginning:
{\em given an infinite dimensional graph complex, $(\cA ss_\infty,\delta)$,
spanned by all possible planar graphs (without wheels) built from $(1,n)$-corollas with $n\geq 2$
and equipped with differential (\ref{2: delta_ass(n)}), then
its cohomology, $H(\cA ss_\infty,\delta)$, is
generated by only  $(1,2)$-corollas}, i.e.\  it is surprisingly small.
It is often impossible to obtain such a result by a direct computation.
One of the main theorem-proving technique
in the theory of operads and props
is called the {\em Koszul duality theory}, and a result of type (\ref{2: Stasheff_theorem})
often requires a combination of ideas from homological algebra,
algebraic topology, the theory of Cohen-Macaulay posets \cite{V2} and so on.
Stasheff \cite{St} proved Theorem~\ref{2: Stasheff_theorem} by constructing a remarkable family
of polytopes  called nowadays {\em associahedra}; in his approach the
surprising  smallness of $H(\cA ss_\infty,\delta)$ gets nicely explained by the obvious
contractibility of Stasheff's polytopes as topological spaces.
We shall review some  theorem-proving techniques in \S 5 and continue this section with a list
of examples which are most relevant to differential geometry.

\subsection{The wheeled operad of finite-dimensional associative algebras}
Theorem~\ref{2: Stasheff_theorem} has been obtained in the category of algebras over the
family of graphs, $\fG^\curlywedge$,
which contain {\em no}\, closed directed paths of internal edges. What happens if we keep the {\em same}\,
family of generators as in the case of $\cA ss$,
$$
A_0(m,n)=\left\{\Ba{cl}
\K[\bS_2]=\mbox{span}\langle
\begin{xy}
 <0mm,0.66mm>*{};<0mm,3mm>*{}**@{-},
 <0.39mm,-0.39mm>*{};<2.2mm,-2.2mm>*{}**@{-},
 <-0.35mm,-0.35mm>*{};<-2.2mm,-2.2mm>*{}**@{-},
 <0mm,0mm>*{\bullet};<0mm,0mm>*{}**@{},
   <0.39mm,-0.39mm>*{};<2.9mm,-4mm>*{^2}**@{},
   <-0.35mm,-0.35mm>*{};<-2.8mm,-4mm>*{^1}**@{},
\end{xy}
\,
,\,
\begin{xy}
 <0mm,0.66mm>*{};<0mm,3mm>*{}**@{-},
 <0.39mm,-0.39mm>*{};<2.2mm,-2.2mm>*{}**@{-},
 <-0.35mm,-0.35mm>*{};<-2.2mm,-2.2mm>*{}**@{-},
 <0mm,0mm>*{\bullet};<0mm,0mm>*{}**@{},
   <0.39mm,-0.39mm>*{};<2.9mm,-4mm>*{^1}**@{},
   <-0.35mm,-0.35mm>*{};<-2.8mm,-4mm>*{^2}**@{},
\end{xy}
\rangle & \mbox{for}\ m=1,n=2\\
0 & \mbox{otherwise}
\Ea
\right.
$$
the {\em same}\, family of relations (\ref{graph-associativity}), but enlarge the family of graphs
we work over from $\fG^\curlywedge$ to
$\fG^\circlearrowright$? The associated  quotient wheeled operad,
$
\cA ss^\circlearrowright:= \cF^\circlearrowright\langle A_0 \rangle/ \langle I_0 \rangle,
$
can be called
 the {\em operad of finite-dimensional associative algebras}. Indeed, one has the following

 \subsubsection{\bf Claim}{\em
There is a 1-1 correspondence between representations,
$\rho:\cA ss^\circlearrowright \rar \cE nd_V$, of $\cA ss^\circlearrowright$
in a finite-dimensional vector space $V$
and  associative algebra structures on $V$.
}

\begin{proof}
We need to explain only the subjective {\em finite-dimensional}, and that follows from the fact
that representations
of graphs $\in \cA ss^\circlearrowright$ which have wheels involve traces. For example, the element
$\begin{xy}
 <0mm,0.66mm>*{};<0mm,2mm>*{}**@{-},
 <0.39mm,-0.39mm>*{};<2mm,-2mm>*{}**@{-},
 <-0.35mm,-0.35mm>*{};<-2.6mm,-2.6mm>*{}**@{-},
 <0mm,0mm>*{\bullet};<0mm,0mm>*{}**@{},
(0.0,0.0)*{}
   \ar@{->}@(u,dr) (1.0,-1.0)*{}
\end{xy}\in \cA ss^\circlearrowright(0,1)$ gets represented in $V$
as the image of the multiplication map $\rho(\begin{xy}
 <0mm,0.66mm>*{};<0mm,3mm>*{}**@{-},
 <0.39mm,-0.39mm>*{};<2.2mm,-2.2mm>*{}**@{-},
 <-0.35mm,-0.35mm>*{};<-2.2mm,-2.2mm>*{}**@{-},
 <0mm,0mm>*{\bullet};<0mm,0mm>*{}**@{},\end{xy})
 \in \Hom(V\ot V, V)$ under a natural  trace map
 $ \Hom(V\ot V, V)\rar \Hom(V,\K)$.
\end{proof}

It is easy to see that the straightforward analogue of Theorem~\ref{2: Stasheff_theorem} can {\em not}\,
hold true
for the operad of {\em finite-dimensional}\, associative algebras
as, for example, formula (\ref{2: delta_ass(3)}) implies
\vspace{-4mm}
$$
\delta
\begin{xy}
 <0mm,-4mm>*{};<0mm,2mm>*{}**@{-},
 <0mm,-4mm>*{};<1mm,-6mm>*{}**@{-},
 <0.39mm,-0.39mm>*{};<3.2mm,-4mm>*{}**@{-},
 <-0.35mm,-0.35mm>*{};<-3.2mm,-4mm>*{}**@{-},
<0mm,0mm>*{\bullet};
   <0.39mm,-0.39mm>*{};<4mm,-6mm>*{^2}**@{},
   <-0.35mm,-0.35mm>*{};<-4mm,-6mm>*{^1}**@{},
(0,2)*{}
   \ar@{->}@(u,dr) (1,-6)*{}
\end{xy}=
\begin{xy}
 <0mm,0mm>*{\bullet};
 <0mm,0.69mm>*{};<0mm,3.0mm>*{}**@{-},
 <0.39mm,-0.39mm>*{};<2.4mm,-2.4mm>*{}**@{-},
 <-0.35mm,-0.35mm>*{};<-1.9mm,-1.9mm>*{}**@{-},
 <-2.4mm,-2.4mm>*{\bullet};
 <-2.0mm,-2.8mm>*{};<0mm,-4.9mm>*{}**@{-},
 <-2.8mm,-2.9mm>*{};<-4.7mm,-4.9mm>*{}**@{-},
    <0.39mm,-0.39mm>*{};<3mm,-4.3mm>*{^2}**@{},
    <-2.8mm,-2.9mm>*{};<-5.2mm,-6.7mm>*{^1}**@{},
(0,2)*{}
   \ar@{->}@(u,dr) (0,-5)*{}
 \end{xy}
\ - \
 \begin{xy}
 <0mm,0mm>*{\bullet};
 <0mm,0.69mm>*{};<0mm,3.0mm>*{}**@{-},
 <0.39mm,-0.39mm>*{};<2.4mm,-2.4mm>*{}**@{-},
 <-0.35mm,-0.35mm>*{};<-1.9mm,-1.9mm>*{}**@{-},
 <2.4mm,-2.4mm>*{\bullet};
 <2.0mm,-2.8mm>*{};<0mm,-4.9mm>*{}**@{-},
 <2.8mm,-2.9mm>*{};<4.7mm,-4.9mm>*{}**@{-},
    <0.39mm,-0.39mm>*{};<-2.3mm,-4.0mm>*{^1}**@{},
    <-2.8mm,-2.9mm>*{};<5.2mm,-6.7mm>*{^2}**@{},
(0,2)*{}
   \ar@{->}@(u,dl) (0,-5)*{}
 \end{xy}
 =
\begin{xy}
 <0mm,0mm>*{\bullet};
 <0mm,0.69mm>*{};<0mm,3.0mm>*{}**@{-},
 <0.39mm,-0.39mm>*{};<2.4mm,-2.4mm>*{}**@{-},
 <-0.35mm,-0.35mm>*{};<-1.9mm,-1.9mm>*{}**@{-},
 <-2.4mm,-2.4mm>*{\bullet};
 <-2.0mm,-2.8mm>*{};<0mm,-4.9mm>*{}**@{-},
 <-2.8mm,-2.9mm>*{};<-4.7mm,-4.9mm>*{}**@{-},
    <0.39mm,-0.39mm>*{};<3mm,-4.3mm>*{^2}**@{},
    <-2.8mm,-2.9mm>*{};<-5.2mm,-6.7mm>*{^1}**@{},
(0,2)*{}
   \ar@{->}@(u,dr) (0,-5)*{}
 \end{xy}
 \ - \
\begin{xy}
 <0mm,0mm>*{\bullet};
 <0mm,0.69mm>*{};<0mm,3.0mm>*{}**@{-},
 <0.39mm,-0.39mm>*{};<2.4mm,-2.4mm>*{}**@{-},
 <-0.35mm,-0.35mm>*{};<-1.9mm,-1.9mm>*{}**@{-},
 <-2.4mm,-2.4mm>*{\bullet};
 <-2.0mm,-2.8mm>*{};<0mm,-4.9mm>*{}**@{-},
 <-2.8mm,-2.9mm>*{};<-4.7mm,-4.9mm>*{}**@{-},
    <0.39mm,-0.39mm>*{};<3mm,-4.3mm>*{^2}**@{},
    <-2.8mm,-2.9mm>*{};<-5.2mm,-6.7mm>*{^1}**@{},
(0,2)*{}
   \ar@{->}@(u,dr) (0,-5)*{}
 \end{xy}
 =0
 \vspace{-3mm}
$$
and hence provides us with a non-trivial cohomology class in
$H^{-1}(\cF^\circlearrowright\langle{A}\rangle,\delta)$ which maps under the natural projection
$\cF^\circlearrowright\langle{A}\rangle\rar
\cA ss^\circlearrowright$ to zero. The correct analogue of Stasheff's result for finite-dimensional
associative algebras was found in \cite{MMS}.

\subsubsection{\bf Theorem}{\em
The minimal resolution of $\cA ss^\circlearrowright$ is a dg free wheeled operad,
$(\cA ss^\circlearrowright)_\infty:=(\cF^\circlearrowright\langle\hat{A}\rangle, \delta)$ generated by an
$\bS$-bimodule $\hat{A}=\{\hat{A}(m,n)\}$,
$$
\hat{A}(m,n):=
\left\{\Ba{cl}
\K[\bS_n][n-2]= \mbox{\em span}
\left\langle
\Ba{c}
\begin{xy}
 <0mm,0mm>*{\bullet};
 <0mm,0mm>*{};<-8mm,-5mm>*{}**@{-},
 <0mm,0mm>*{};<-4.5mm,-5mm>*{}**@{-},
 <0mm,0mm>*{};<0mm,-4mm>*{\ldots}**@{},
 <0mm,0mm>*{};<4.5mm,-5mm>*{}**@{-},
 <0mm,0mm>*{};<8mm,-5mm>*{}**@{-},
   <0mm,0mm>*{};<-11mm,-7.9mm>*{^{\sigma(1)}}**@{},
   <0mm,0mm>*{};<-4mm,-7.9mm>*{^{\sigma(2)}}**@{},
   <0mm,0mm>*{};<10.0mm,-7.9mm>*{^{\sigma(n)}}**@{},
 <0mm,0mm>*{};<0mm,5mm>*{}**@{-},
 \end{xy}\Ea
\right\rangle_{\sigma\in \bS_n} &\hspace{-2mm} \mbox{for}\ m=1, n\geq 2\\
\hspace{-2mm}
\bigoplus_{p=1}^{n-1}\K[\bS_n]_{C_p\times C_{n-p}}[n]
= \mbox{\em span}\hspace{-0.3mm}
\left\langle\hspace{-1mm}
\Ba{c}
\begin{xy}
 <0mm,-0.5mm>*{\blacktriangledown};
 <0mm,0mm>*{};<-16mm,-5mm>*{}**@{-},
 <0mm,0mm>*{};<-11mm,-5mm>*{}**@{-},
 <0mm,0mm>*{};<-3.5mm,-5mm>*{}**@{-},
 <0mm,0mm>*{};<-6mm,-5mm>*{...}**@{},
   <0mm,0mm>*{};<-16mm,-7.9mm>*{^{\sigma(1)}}**@{},
   <0mm,0mm>*{};<-10mm,-7.9mm>*{^{\sigma(2)}}**@{},
   <0mm,0mm>*{};<-3mm,-7.9mm>*{^{\sigma(p)}}**@{},
 <0mm,0mm>*{};<16mm,-5mm>*{}**@{-},
 <0mm,0mm>*{};<12mm,-5mm>*{}**@{-},
 <0mm,0mm>*{};<3.5mm,-5mm>*{}**@{-},
 <0mm,0mm>*{};<6.6mm,-5mm>*{...}**@{},
   <0mm,0mm>*{};<16mm,-7.9mm>*{^{\sigma(n)}}**@{},
   <0mm,0mm>*{};<6mm,-7.9mm>*{^{\sigma(p\hspace{-0.1mm}+\hspace{-0.5mm}1)}}**@{},
 \end{xy}
 \Ea
 \hspace{-2mm}
\right\rangle_{\hspace{-1mm}\sigma\in \bS_n} & \hspace{-2mm}\mbox{for}\ m=0, n\geq 2\\
0 &\hspace{-2mm} \mbox{otherwise}
\Ea
\right.
$$
where $C_p\times C_{n-p}$  is the  subgroup of\, $\bS_n$ generated by two commuting
cyclic permutations $\zeta:=(12\ldots p)$ and $\xi:=(p+1\ldots n)$, and
$k[\bS_n]_{C_p\times C_{n-p}}$ stands for coinvariants.

The differential is given on the generators of $\hat{A}(1,n)$  by (\ref{2: delta_ass(n)}) and
on the generators of  $\hat{A}(0,n)$ by\vspace{-3mm}
\Beqrn
\delta\hspace{-1mm}
\begin{xy}
 <0mm,-0.5mm>*{\blacktriangledown};<0mm,0mm>*{}**@{},
 <0mm,0mm>*{};<-16mm,-5mm>*{}**@{-},
 <0mm,0mm>*{};<-12mm,-5mm>*{}**@{-},
 <0mm,0mm>*{};<-3.5mm,-5mm>*{}**@{-},
 <0mm,0mm>*{};<-6mm,-5mm>*{...}**@{},
   <0mm,0mm>*{};<-16mm,-7.9mm>*{^{1}}**@{},
   <0mm,0mm>*{};<-11mm,-7.9mm>*{^{2}}**@{},
   <0mm,0mm>*{};<-4mm,-7.9mm>*{^{p}}**@{},
 <0mm,0mm>*{};<16mm,-5mm>*{}**@{-},
 <0mm,0mm>*{};<12mm,-5mm>*{}**@{-},
 <0mm,0mm>*{};<3.5mm,-5mm>*{}**@{-},
 <0mm,0mm>*{};<6.6mm,-5mm>*{...}**@{},
   <0mm,0mm>*{};<16.5mm,-7.9mm>*{^{n}}**@{},
   <0mm,0mm>*{};<5mm,-7.9mm>*{^{p\hspace{-0mm}+\hspace{-0.5mm}1}}**@{},
 \end{xy}
 \hspace{-2mm}
&=&
\hspace{-2mm}
\displaystyle\underset{(1\ldots p)}{\oint}\underset{(p+1\ldots n)}{\oint}
\left(\rule{0em}{2.3em}\right.
\begin{xy}
 <0mm,0mm>*{\bullet};<0mm,0mm>*{}**@{},
 <0mm,0mm>*{};<-11mm,-5mm>*{}**@{-},
 <0mm,0mm>*{};<-8.5mm,-5mm>*{}**@{-},
 <0mm,0mm>*{};<-5mm,-5mm>*{...}**@{},
 <0mm,0mm>*{};<-2.5mm,-5mm>*{}**@{-},
   <0mm,0mm>*{};<-11mm,-7.9mm>*{^1}**@{},
   <0mm,0mm>*{};<-7.6mm,-7.9mm>*{^2}**@{},
   <0mm,0mm>*{};<-3mm,-7.9mm>*{^p}**@{},
<0mm,0mm>*{};<11mm,-5mm>*{}**@{-},
 <0mm,0mm>*{};<8.5mm,-5mm>*{}**@{-},
 <0mm,0mm>*{};<5.5mm,-5mm>*{...}**@{},
 <0mm,0mm>*{};<2.5mm,-5mm>*{}**@{-},
   <0mm,0mm>*{};<11.8mm,-7.9mm>*{^{n}}**@{},
   <0mm,0mm>*{};<3mm,-7.9mm>*{^{p\hspace{-0mm}+\hspace{-0.5mm}1}}**@{},
 <0mm,-10mm>*{};<0mm,5mm>*{}**@{-},
(0,5)*{}
   \ar@{->}@(ur,dr) (0,-10)*{}
 \end{xy}
 +
 \sum_{k=2}^{p}
(-1)^{kn} \hskip -2em
\begin{xy}
 <0mm,-0.5mm>*{\blacktriangledown};<0mm,0mm>*{}**@{},
 <0mm,0mm>*{};<-16mm,-5mm>*{}**@{-},
 <0mm,0mm>*{};<-11mm,-5mm>*{}**@{-},
 <0mm,0mm>*{};<-3.5mm,-5mm>*{}**@{-},
 <0mm,0mm>*{};<-6mm,-5mm>*{...}**@{},
   <0mm,0mm>*{};<-10mm,-7.9mm>*{^{k\hspace{-0mm}+\hspace{-0.5mm}1}}**@{},
   <0mm,0mm>*{};<-4mm,-7.9mm>*{^{p}}**@{},
 <0mm,0mm>*{};<16mm,-5mm>*{}**@{-},
 <0mm,0mm>*{};<12mm,-5mm>*{}**@{-},
 <0mm,0mm>*{};<3.5mm,-5mm>*{}**@{-},
 <0mm,0mm>*{};<6.6mm,-5mm>*{...}**@{},
   <0mm,0mm>*{};<16.6mm,-7.9mm>*{^{n}}**@{},
   <0mm,0mm>*{};<5mm,-7.9mm>*{^{p\hspace{0mm}+\hspace{-0.5mm}1}}**@{},
 <-16mm,-5.5mm>*{\bullet};<0mm,0mm>*{}**@{},
 <-16mm,-5.5mm>*{};<-20mm,-11mm>*{}**@{-},
 <-16mm,-5.5mm>*{};<-12mm,-11mm>*{}**@{-},
 <-16mm,-5.5mm>*{};<-18mm,-11mm>*{}**@{-},
 <-16mm,-5.5mm>*{};<-14mm,-11mm>*{}**@{-},
 <-16mm,-13mm>*{_{1\ \, \dots\ \ k}},
 \end{xy}
\\
&&
+ \sum_{k=2}^{n-2}
(-1)^{p+k(1+n-p)+1}
\begin{xy}
<0mm,-0.5mm>*{\blacktriangledown};<0mm,0mm>*{}**@{},
 <0mm,0mm>*{};<-16mm,-5mm>*{}**@{-},
 <0mm,0mm>*{};<-12mm,-5mm>*{}**@{-},
 <0mm,0mm>*{};<-3.5mm,-5mm>*{}**@{-},
 <0mm,0mm>*{};<-6mm,-5mm>*{...}**@{},
   <0mm,0mm>*{};<-16.7mm,-7.9mm>*{^{1}}**@{},
   <0mm,0mm>*{};<-12mm,-7.9mm>*{^{2}}**@{},
   <0mm,0mm>*{};<-4mm,-7.9mm>*{^{p}}**@{},
 <0mm,0mm>*{};<17mm,-5mm>*{}**@{-},
 <0mm,0mm>*{};<7mm,-5mm>*{}**@{-},
 <0mm,0mm>*{};<3.5mm,-5mm>*{}**@{-},
 <0mm,0mm>*{};<11.6mm,-5mm>*{...}**@{},
   <0mm,0mm>*{};<19mm,-7.9mm>*{^{n}}**@{},
   <0mm,0mm>*{};<11mm,-7.9mm>*{^{p\hspace{-0mm}+\hspace{-0.2mm}k+\hspace{-0.2mm}1}}**@{},
 <3.5mm,-5.5mm>*{\bullet};<0mm,0mm>*{}**@{},
 <3.5mm,-5.5mm>*{};<-0.5mm,-11mm>*{}**@{-},
 <3.5mm,-5.5mm>*{};<1.5mm,-11mm>*{}**@{-},
 <3.5mm,-5.5mm>*{};<5.5mm,-11mm>*{}**@{-},
 <3.5mm,-5.5mm>*{};<7.5mm,-11mm>*{}**@{-},
 <3.5mm,-13mm>*{_{p+1\, \dots\ p+k}},
 \end{xy}\left. \rule{0em}{2.3em}\right)
\Eeqrn
where the symbol $\displaystyle\underset{(i_1\ldots i_k)}{\oint}$ stands for the cyclic skewsymmetrization
of the indices $(i_1\ldots i_k)$.
}

Thus the minimal resolution, $(\cA ss^\circlearrowright)_\infty$,
 of the operad of finite-dimensional associative algebras is different
from the naive ``wheelification", $(\cA ss_\infty)^\circlearrowright$, of the Stasheff's
minimal resolution of the operad, $\cA ss$, of arbitrary associative algebras.
A similar phenomenon occurs for
the operad of commutative algebras \cite{MMS}.
In contrast, the {operad}, $\caL ie$, of Lie algebras
 is rigid
with respect to the wheelification:

\subsubsection{\bf Fact \cite{Me3}}
$
(\caL ie^\circlearrowright)_\infty
 = \ (\caL ie_\infty)^\circlearrowright$,
{ i.e.\ wheeled $L_\infty$-algebras are exactly the same
as ordinary finite-dimensional $L_\infty$-algebras.}

\subsubsection{\bf Reminder on $L_\infty$-algebras and their
 homotopy classification}\label{2: subsub on Lie_infty}
For future reference we recall here a few useful facts about Lie and $L_\infty$-algebras \cite{Ko}.
The operad, $\caL ie$, of Lie algebras is the quotient operad,
$\caL ie:= \cF^\curlywedge \langle L_0 \rangle/I$, of the free operad generated by an $\bS$-bimodule
$L_0=\{L_0(m,n)\}$,
\Beq\label{2: Lie generators}
L_0(m,n)=\left\{\Ba{cl}
\sgn_2 =\mbox{span}\langle
\begin{xy}
 <0mm,0.66mm>*{};<0mm,3mm>*{}**@{-},
 <0.39mm,-0.39mm>*{};<2.2mm,-2.2mm>*{}**@{-},
 <-0.35mm,-0.35mm>*{};<-2.2mm,-2.2mm>*{}**@{-},
 <0mm,0mm>*{\bullet};
   <0.39mm,-0.39mm>*{};<2.9mm,-4mm>*{^2}**@{},
   <-0.35mm,-0.35mm>*{};<-2.8mm,-4mm>*{^1}**@{},
\end{xy}
=-
\begin{xy}
 <0mm,0.66mm>*{};<0mm,3mm>*{}**@{-},
 <0.39mm,-0.39mm>*{};<2.2mm,-2.2mm>*{}**@{-},
 <-0.35mm,-0.35mm>*{};<-2.2mm,-2.2mm>*{}**@{-},
 <0mm,0mm>*{\bullet};<0mm,0mm>*{}**@{},
   <0.39mm,-0.39mm>*{};<2.9mm,-4mm>*{^1}**@{},
   <-0.35mm,-0.35mm>*{};<-2.8mm,-4mm>*{^2}**@{},
\end{xy}
\rangle & \mbox{for}\ m=1,n=2\\
0 & \mbox{otherwise}
\Ea
\right.
\Eeq
modulo the ideal $I$ generated by the following relations
\Beq\label{2: jacobi equations}
 \begin{xy}
 <0mm,0mm>*{\bullet};<0mm,0mm>*{}**@{},
 <0mm,0.69mm>*{};<0mm,3.0mm>*{}**@{-},
 <0.39mm,-0.39mm>*{};<2.4mm,-2.4mm>*{}**@{-},
 <-0.35mm,-0.35mm>*{};<-1.9mm,-1.9mm>*{}**@{-},
 <-2.4mm,-2.4mm>*{\bullet};<-2.4mm,-2.4mm>*{}**@{},
 <-2.0mm,-2.8mm>*{};<0mm,-4.9mm>*{}**@{-},
 <-2.8mm,-2.9mm>*{};<-4.7mm,-4.9mm>*{}**@{-},
    <0.39mm,-0.39mm>*{};<3.3mm,-4.0mm>*{^3}**@{},
    <-2.0mm,-2.8mm>*{};<0.5mm,-6.7mm>*{^2}**@{},
    <-2.8mm,-2.9mm>*{};<-5.2mm,-6.7mm>*{^1}**@{},
 \end{xy}
\ + \
 \begin{xy}
 <0mm,0mm>*{\bullet};<0mm,0mm>*{}**@{},
 <0mm,0.69mm>*{};<0mm,3.0mm>*{}**@{-},
 <0.39mm,-0.39mm>*{};<2.4mm,-2.4mm>*{}**@{-},
 <-0.35mm,-0.35mm>*{};<-1.9mm,-1.9mm>*{}**@{-},
 <-2.4mm,-2.4mm>*{\bullet};<-2.4mm,-2.4mm>*{}**@{},
 <-2.0mm,-2.8mm>*{};<0mm,-4.9mm>*{}**@{-},
 <-2.8mm,-2.9mm>*{};<-4.7mm,-4.9mm>*{}**@{-},
    <0.39mm,-0.39mm>*{};<3.3mm,-4.0mm>*{^2}**@{},
    <-2.0mm,-2.8mm>*{};<0.5mm,-6.7mm>*{^1}**@{},
    <-2.8mm,-2.9mm>*{};<-5.2mm,-6.7mm>*{^3}**@{},
 \end{xy}
\ + \
 \begin{xy}
 <0mm,0mm>*{\bullet};<0mm,0mm>*{}**@{},
 <0mm,0.69mm>*{};<0mm,3.0mm>*{}**@{-},
 <0.39mm,-0.39mm>*{};<2.4mm,-2.4mm>*{}**@{-},
 <-0.35mm,-0.35mm>*{};<-1.9mm,-1.9mm>*{}**@{-},
 <-2.4mm,-2.4mm>*{\bullet};<-2.4mm,-2.4mm>*{}**@{},
 <-2.0mm,-2.8mm>*{};<0mm,-4.9mm>*{}**@{-},
 <-2.8mm,-2.9mm>*{};<-4.7mm,-4.9mm>*{}**@{-},
    <0.39mm,-0.39mm>*{};<3.3mm,-4.0mm>*{^1}**@{},
    <-2.0mm,-2.8mm>*{};<0.5mm,-6.7mm>*{^3}**@{},
    <-2.8mm,-2.9mm>*{};<-5.2mm,-6.7mm>*{^2}**@{},
 \end{xy}=0.
 \Eeq
Its minimal resolution, $\caL ie_\infty$, is a dg free operad, $\cF^\curlywedge\langle L\rangle$
generated by an $\bS_n$-bimodule
$$\label{2: Lie infty generators}
L(m,n):=\left\{\Ba{cl}
\sgn_n[n-2]=\mbox{span}\left\langle
\begin{xy}
 <0mm,0mm>*{\bullet};<0mm,0mm>*{}**@{},
 <0mm,0mm>*{};<0mm,5mm>*{}**@{-},
   %
<0mm,0mm>*{\bullet};<0mm,0mm>*{}**@{},
 <0mm,0mm>*{};<-6mm,-5mm>*{}**@{-},
 <0mm,0mm>*{};<-3.1mm,-5mm>*{}**@{-},
 <0mm,0mm>*{};<0mm,-4.6mm>*{...}**@{},
 <0mm,0mm>*{};<3.1mm,-5mm>*{}**@{-},
 <0mm,0mm>*{};<6mm,-5mm>*{}**@{-},
   <0mm,0mm>*{};<-6.7mm,-6.4mm>*{_1}**@{},
   <0mm,0mm>*{};<-3.2mm,-6.4mm>*{_2}**@{},
   <0mm,0mm>*{};<3.1mm,-6.4mm>*{_{n\mbox{-}1}}**@{},
   <0mm,0mm>*{};<6.9mm,-6.4mm>*{_{n}}**@{},
 \end{xy}
\right\rangle & \mbox{for}\ m=1, n\geq 2,\\
0             & \mbox{otherwise}\\
\Ea
\right.
$$
with the differential given by
\Beq\label{2: Lie infty differential}
\delta\ \begin{xy}
 <0mm,0mm>*{\bullet};
 <0mm,0mm>*{};<0mm,4mm>*{}**@{-},
<0mm,0mm>*{\bullet};
 <0mm,0mm>*{};<-6mm,-5mm>*{}**@{-},
 <0mm,0mm>*{};<-3.1mm,-5mm>*{}**@{-},
 <0mm,0mm>*{};<0mm,-4.6mm>*{...}**@{},
 <0mm,0mm>*{};<3.1mm,-5mm>*{}**@{-},
 <0mm,0mm>*{};<6mm,-5mm>*{}**@{-},
   <0mm,0mm>*{};<-6.7mm,-6.4mm>*{_1}**@{},
   <0mm,0mm>*{};<-3.2mm,-6.4mm>*{_2}**@{},
   <0mm,0mm>*{};<3.1mm,-6.4mm>*{_{n\mbox{-}1}}**@{},
   <0mm,0mm>*{};<6.9mm,-6.4mm>*{_{n}}**@{},
 \end{xy}
 =
 \sum_{[n]=I_1\sqcup I_2\atop {\atop
 {\# I_1\geq 2, \# I_2\geq 1}}
}\hspace{0mm}
(-1)^{\sigma(I_1,I_2) + (\# I_1+1)\# I_2}
\begin{xy}
 <0mm,0mm>*{\bullet};<0mm,0mm>*{}**@{},
 <0mm,0mm>*{};<0mm,4mm>*{}**@{-},
 <0mm,0mm>*{};<-7mm,-5mm>*{}**@{-},
<-7mm,-5mm>*{\bullet};
 <-7mm,-5mm>*{};<-12mm,-10mm>*{}**@{-},
 <-7mm,-5mm>*{};<-2.2mm,-10mm>*{}**@{-},
 <-7mm,-5mm>*{};<-9.5mm,-10mm>*{}**@{-},
 <-7mm,-5mm>*{};<-4.55mm,-10mm>*{}**@{-},
  <-7mm,-5mm>*{};<-7mm,-9.5mm>*{...}**@{},
 <0mm,0mm>*{};<-7mm,-11.8mm>*{\underbrace{\ \ \ \ \ \ \ \ \
      }}**@{},
      <0mm,0mm>*{};<-7mm,-15.2mm>*{^{I_1}}**@{},
 <0mm,0mm>*{};<2mm,-6.4mm>*{\underbrace{\ \ \ \ \ \ \ \ \ \
      }}**@{},
    <0mm,0mm>*{};<2.8mm,-9.8mm>*{^{I_2}}**@{},
 <0mm,0mm>*{};<-3.5mm,-5mm>*{}**@{-},
 <0mm,0mm>*{};<-0mm,-4.6mm>*{...}**@{},
 <0mm,0mm>*{};<3.6mm,-5mm>*{}**@{-},
 <0mm,0mm>*{};<7mm,-5mm>*{}**@{-},
 \end{xy}.
 \vspace{-2mm}
\Eeq
Here (and elsewhere) $\sgn_n$ stands for the 1-dimensional sign representation of $\bS_n$.

\sip

With an arbitrary graded vector space $V$ one can associate a formal graded manifold, $\cM_V$,
whose structure sheaf, $\f_{\cM_V}$, is, by definition, the completed graded cocommutative coalgebra
 $\widehat{\odot}(V[1])$; if $V$ is finite dimensional, then one can equivalently view $\cM_V$ as
a small neighbourhood  of zero in the space $V[1]$ equipped with the algebra (rather than coalgebra),
 $\widehat{\odot}(V^*[-1])$, of ordinary smooth formal functions.
It is well known (see, e.g., \cite{Ko}) that  $L_\infty$-structures in a dg space $V$,
that is, representations $\caL_\infty\rar \cE nd_V$, are in one-to-one correspondence with degree
1 vector fields, $\eth$, on  $\cM_V$ which vanish at the distinguished point, $\eth\mid_{0\in \cM_V} =0$,
and satisfy the condition $[\eth,\eth]=0$ (such vector fields are called {\em cohomological}).
The pairs $(\cM_V, \eth)$ are often called {\em dg manifolds}. This interpretation   of
$L_\infty$-structures permits us to use
simple and concise
geometric instruments to describe notions which, in the pure algebraic translation,
look awkwardly large. For example, a {\em morphism}\, of $L_\infty$-algebras $V\rar W$ is nothing
but a smooth map, $f: \cM_{V}\rar \cM_{W}$, of the associated formal manifolds such that
$f_*(\eth_V)=\eth_W$.
\sip

A $L_\infty$ algebra $(\cM_V, \eth)$, is called {\em minimal}\,
if the first Taylor coefficient, $\eth_{(1)}$, of the homological
vector field $\eth$ at the distinguished point $0\in \cM_V$ vanishes. It is called {\em linear
contractible}\, if the higher Taylor coefficients $\eth_{(\geq 2)}$
vanish and the first one $\eth_{(1)}$ has trivial cohomology when
viewed as a differential in $V$.
According to Kontsevich \cite{Ko}, {\em any $L_\infty$-algebra (or, better, the associated dg manifold) is
isomorphic to the direct product of a minimal and of a
linear contractible one}.  This fact implies that quasi-isomorphisms in the category of $L_\infty$-algebras
 are equivalence relations. A dg manifold is called {\em
contractible}\, if it is isomorphic  to
 a linear contractible one.

\subsection{Unimodular Lie algebras}\label{3: subsect unimod Lie} Many important Lie algebras $\fg$ (e.g., all semisimple Lie algebras)
have the additional property that, for any $g\in \fg$,  the trace of the associated adjoint
action
$$
\Ba{rccc}
\mbox{Ad}_g: & \fg &\lon & \fg \\
         & e  & \lon & [g,e]
\Ea
$$
vanishes. Lie algebras with this property are called {\em unimodular}. The wheeled operad, $\cU\caL ie$,
controlling unimodular Lie algebras is the quotient of the free wheeled operad,
$\cF^\circlearrowright \langle L_0 \rangle$, generated by the $\bS$-bimodule (\ref{2: Lie generators})
modulo the ideal generated by the Jacobi relations (\ref{2: jacobi equations}) and the unimodularity relation,
$
\begin{xy}
 <0mm,0.66mm>*{};<0mm,3mm>*{}**@{-},
 <0.39mm,-0.39mm>*{};<2.2mm,-2.2mm>*{}**@{-},
 <-0.35mm,-0.35mm>*{};<-3mm,-3mm>*{}**@{-},
 <0mm,0mm>*{\bullet};
(0,2)*{}
   \ar@{->}@(u,dr) (2.2,-2.2)*{}
\end{xy}=0.
$
Its minimal resolution has been found in \cite{Gr1}:

\subsubsection{\bf Theorem} {\em The operad $\cU\caL ie_\infty$ is a dg free operad,
$\cF^\circlearrowright\langle\hat{L}\rangle$
generated by the $\bS$-bimodule},
$$
\hat{L}(m,n):=\left\{\Ba{ll}
\sgn_n[n-2]  & \mbox{for}\ m=1, n\geq 2,\\
\sgn_n[n]=
\mbox{span}\left\langle
\Ba{c}
\begin{xy}
 <0mm,0mm>*{\bullet};
 <0mm,0mm>*{};<-6mm,-5mm>*{}**@{-},
 <0mm,0mm>*{};<-3.1mm,-5mm>*{}**@{-},
 <0mm,0mm>*{};<0mm,-4.6mm>*{...}**@{},
 <0mm,0mm>*{};<3.1mm,-5mm>*{}**@{-},
 <0mm,0mm>*{};<6mm,-5mm>*{}**@{-},
   <0mm,0mm>*{};<-6.7mm,-6.4mm>*{_1}**@{},
   <0mm,0mm>*{};<-3.2mm,-6.4mm>*{_2}**@{},
   <0mm,0mm>*{};<3.1mm,-6.4mm>*{_{n\mbox{-}1}}**@{},
   <0mm,0mm>*{};<6.9mm,-6.4mm>*{_{n}}**@{},
 \end{xy}
 \Ea
\right\rangle & \mbox{for}\ m=0, n\geq 1,\\
0             & \mbox{otherwise}\\
\Ea
\right.
$$
{\em
with the differential on the generators of $\hat{L}(1,n)$ given by (\ref{2: Lie infty differential})
and on the generators of  $\hat{L}(1,n)$ by}
$$
\delta
\begin{xy}
 <0mm,0mm>*{\bullet};
<0mm,0mm>*{\bullet};
 <0mm,0mm>*{};<-6mm,-5mm>*{}**@{-},
 <0mm,0mm>*{};<-3.1mm,-5mm>*{}**@{-},
 <0mm,0mm>*{};<0mm,-4.6mm>*{...}**@{},
 <0mm,0mm>*{};<3.1mm,-5mm>*{}**@{-},
 <0mm,0mm>*{};<6mm,-5mm>*{}**@{-},
   <0mm,0mm>*{};<-6.7mm,-6.4mm>*{_1}**@{},
   <0mm,0mm>*{};<-3.2mm,-6.4mm>*{_2}**@{},
   <0mm,0mm>*{};<3.1mm,-6.4mm>*{_{n\mbox{-}1}}**@{},
   <0mm,0mm>*{};<6.9mm,-6.4mm>*{_{n}}**@{},
 \end{xy}
 =
 \sum_{[n]=I_1\sqcup I_2\atop {\atop
 {\# I_1\geq 2, \# I_2\geq 0}}
}
(-1)^{\sigma(I_1,I_2) + (\# I_1+1)\# I_2}\hspace{-6mm}
\begin{xy}
 <0mm,0mm>*{\bullet};<0mm,0mm>*{}**@{},
 <0mm,0mm>*{};<-7mm,-5mm>*{}**@{-},
<-7mm,-5mm>*{\bullet};
 <-7mm,-5mm>*{};<-12mm,-10mm>*{}**@{-},
 <-7mm,-5mm>*{};<-2.2mm,-10mm>*{}**@{-},
 <-7mm,-5mm>*{};<-9.5mm,-10mm>*{}**@{-},
 <-7mm,-5mm>*{};<-4.55mm,-10mm>*{}**@{-},
  <-7mm,-5mm>*{};<-7mm,-9.5mm>*{...}**@{},
 <0mm,0mm>*{};<-7mm,-11.8mm>*{\underbrace{\ \ \ \ \ \ \ \ \
      }}**@{},
      <0mm,0mm>*{};<-7mm,-15.2mm>*{^{I_1}}**@{},
 <0mm,0mm>*{};<2mm,-6.4mm>*{\underbrace{\ \ \ \ \ \ \ \ \ \
      }}**@{},
    <0mm,0mm>*{};<2.8mm,-9.8mm>*{^{I_2}}**@{},
 <0mm,0mm>*{};<-3.5mm,-5mm>*{}**@{-},
 <0mm,0mm>*{};<-0mm,-4.6mm>*{...}**@{},
 <0mm,0mm>*{};<3.6mm,-5mm>*{}**@{-},
 <0mm,0mm>*{};<7mm,-5mm>*{}**@{-},
 \end{xy}
 \
 +
 \
\begin{xy}
 <0mm,0mm>*{\bullet};
<0mm,0mm>*{};<0mm,3mm>*{}**@{-},
<0mm,0mm>*{\bullet};
 <0mm,0mm>*{};<-6mm,-5mm>*{}**@{-},
 <0mm,0mm>*{};<-3.1mm,-5mm>*{}**@{-},
 <0mm,0mm>*{};<0mm,-4.6mm>*{...}**@{},
 <0mm,0mm>*{};<3.1mm,-5mm>*{}**@{-},
 <0mm,0mm>*{};<6mm,-5mm>*{}**@{-},
   <0mm,0mm>*{};<-6.7mm,-6.4mm>*{_1}**@{},
   <0mm,0mm>*{};<-3.2mm,-6.4mm>*{_2}**@{},
   <0mm,0mm>*{};<3.1mm,-6.4mm>*{_{n}}**@{},
%
(0,2)*{}
   \ar@{->}@(u,dr) (5,-4.2)*{},
 \end{xy}
 .
$$

Geometrically, unimodular $L_\infty$-structures in $V$ can be interpreted
as pairs, $(\eth, \om)$, where $\eth$ is a cohomological vector field  and  $\om$
 a $Q$-invariant section of the Berezinian bundle on  $V^*[1]$ (see \cite{Gr1}).

\subsection{Lie $1$-bialgebras and Poisson geometry}\label{2: subsect LB and Poisson}
 A  {\em Lie n-bialgebra} on a graded vector space $V$ is a pair of linear maps,
$$
\Delta\simeq
\begin{xy}
 <0mm,-0.55mm>*{};<0mm,-2.5mm>*{}**@{-},
 <0.5mm,0.5mm>*{};<2.2mm,2.2mm>*{}**@{-},
 <-0.48mm,0.48mm>*{};<-2.2mm,2.2mm>*{}**@{-},
 <0mm,0mm>*{\circ};<0mm,0mm>*{}**@{},
   <0mm,-0.55mm>*{};<0mm,-4.6mm>*{_1}**@{},
   <0.5mm,0.5mm>*{};<2.7mm,3.6mm>*{^2}**@{},
   <-0.48mm,0.48mm>*{};<-2.7mm,3.6mm>*{^1}**@{},
 \end{xy}
: V\rightarrow V\wedge V, \ \ \
[\ \bullet\ ]\simeq \begin{xy}
 <0mm,0.66mm>*{};<0mm,3mm>*{}**@{-},
 <0.39mm,-0.39mm>*{};<2.2mm,-2.2mm>*{}**@{-},
 <-0.35mm,-0.35mm>*{};<-2.2mm,-2.2mm>*{}**@{-},
 <0mm,0mm>*{\bullet};<0mm,0mm>*{}**@{},
   <0mm,0.66mm>*{};<0mm,4mm>*{^1}**@{},
   <0.39mm,-0.39mm>*{};<2.9mm,-5mm>*{^2}**@{},
   <-0.35mm,-0.35mm>*{};<-2.8mm,-5mm>*{^1}**@{},
\end{xy}:
\wedge^2 (V[-n]) \rightarrow V[-n]
$$
making the space $V$ into a Lie coalgebra and the space  $V[-n]$ into a Lie algebra
and satisfying, for any $ a,b\in V$, the compatibility condition
$$
\Delta[a\bullet b] = \sum a_1\otimes [a_2\bullet b] +  [a\bullet
b_1]\otimes b_2 + (-1)^{|a||b|+n|a|+n|b|}( [b\bullet a_1]\otimes a_2
+ b_1\otimes [b_2\bullet a]),
$$
Here $\Delta a=:\sum a_1\otimes a_2$ and  $\Delta b=:\sum
b_1\otimes b_2$.
The case  $n=0$  gives  the notion of Lie bialgebra which was introduced by Drinfeld \cite{Dr}
in the context of quantum groups.
The case $n=1$, as we shall see below, is relevant to Poisson geometry. In this case one has
$\wedge^2 (V[-1])= (\odot^2V)[-2]$ so that the basic binary operations have the following symmetries,
$
\begin{xy}
 <0mm,-0.55mm>*{};<0mm,-2.5mm>*{}**@{-},
 <0.5mm,0.5mm>*{};<2.2mm,2.2mm>*{}**@{-},
 <-0.48mm,0.48mm>*{};<-2.2mm,2.2mm>*{}**@{-},
 <0mm,0mm>*{\circ};<0mm,0mm>*{}**@{},
   <0mm,-0.55mm>*{};<0mm,-4.6mm>*{_1}**@{},
   <0.5mm,0.5mm>*{};<2.7mm,3.6mm>*{^2}**@{},
   <-0.48mm,0.48mm>*{};<-2.7mm,3.6mm>*{^1}**@{},
 \end{xy}
 = -\,
 \begin{xy}
 <0mm,-0.55mm>*{};<0mm,-2.5mm>*{}**@{-},
 <0.5mm,0.5mm>*{};<2.2mm,2.2mm>*{}**@{-},
 <-0.48mm,0.48mm>*{};<-2.2mm,2.2mm>*{}**@{-},
 <0mm,0mm>*{\circ};<0mm,0mm>*{}**@{},
   <0mm,-0.55mm>*{};<0mm,-4.6mm>*{_1}**@{},
   <0.5mm,0.5mm>*{};<2.7mm,3.6mm>*{^1}**@{},
   <-0.48mm,0.48mm>*{};<-2.7mm,3.6mm>*{^2}**@{},
 \end{xy}
$
and
$
\begin{xy}
 <0mm,0.66mm>*{};<0mm,3mm>*{}**@{-},
 <0.39mm,-0.39mm>*{};<2.2mm,-2.2mm>*{}**@{-},
 <-0.35mm,-0.35mm>*{};<-2.2mm,-2.2mm>*{}**@{-},
 <0mm,0mm>*{\bullet};<0mm,0mm>*{}**@{},
   <0mm,0.66mm>*{};<0mm,4mm>*{^1}**@{},
   <0.39mm,-0.39mm>*{};<2.9mm,-5mm>*{^2}**@{},
   <-0.35mm,-0.35mm>*{};<-2.8mm,-5mm>*{^1}**@{},
\end{xy}
 =
 \begin{xy}
 <0mm,0.66mm>*{};<0mm,3mm>*{}**@{-},
 <0.39mm,-0.39mm>*{};<2.2mm,-2.2mm>*{}**@{-},
 <-0.35mm,-0.35mm>*{};<-2.2mm,-2.2mm>*{}**@{-},
 <0mm,0mm>*{\bullet};<0mm,0mm>*{}**@{},
   <0mm,0.66mm>*{};<0mm,4mm>*{^1}**@{},
   <0.39mm,-0.39mm>*{};<2.9mm,-5mm>*{^1}**@{},
   <-0.35mm,-0.35mm>*{};<-2.8mm,-5mm>*{^2}**@{},
\end{xy}
$.
Thus the prop of Lie 1-bialgebras, $\LB$,  is the quotient of the free prop, $\cF^\uparrow\langle B \rangle$,
 generated by an $\bS$-bimodule,
\Beq\label{2: generators LB}
B(m,n):=\left\{
\Ba{ll}
\sgn_2\ot \id_1=\mbox{span}\left\langle
\begin{xy}
 <0mm,-0.55mm>*{};<0mm,-2.5mm>*{}**@{-},
 <0.5mm,0.5mm>*{};<2.2mm,2.2mm>*{}**@{-},
 <-0.48mm,0.48mm>*{};<-2.2mm,2.2mm>*{}**@{-},
 <0mm,0mm>*{\circ};<0mm,0mm>*{}**@{},
 <0mm,-0.55mm>*{};<0mm,-3.8mm>*{_1}**@{},
 <0.5mm,0.5mm>*{};<2.7mm,2.8mm>*{^2}**@{},
 <-0.48mm,0.48mm>*{};<-2.7mm,2.8mm>*{^1}**@{},
 \end{xy}
   \right\rangle  & \mbox{if}\ m=2, n=1,\vspace{3mm}\\
\id_1\ot \id_2[-1]\equiv
\mbox{span}\left\langle
\begin{xy}
 <0mm,0.66mm>*{};<0mm,3mm>*{}**@{-},
 <0.39mm,-0.39mm>*{};<2.2mm,-2.2mm>*{}**@{-},
 <-0.35mm,-0.35mm>*{};<-2.2mm,-2.2mm>*{}**@{-},
 <0mm,0mm>*{\bullet};
   <0mm,0.66mm>*{};<0mm,3.4mm>*{^1}**@{},
   <0.39mm,-0.39mm>*{};<2.9mm,-4mm>*{^2}**@{},
   <-0.35mm,-0.35mm>*{};<-2.8mm,-4mm>*{^1}**@{},
\end{xy}
\right\rangle
\ & \mbox{if}\ m=1, n=2, \vspace{3mm}\\
0 & \mbox{otherwise}
\Ea
\right.
\Eeq
modulo the ideal generated by Jacobi relations (\ref{2: jacobi equations}) and the following ones
\Beq\label{2: Lie 1bi relations}
\begin{xy}
 <0mm,0mm>*{\circ};<0mm,0mm>*{}**@{},
 <0mm,-0.49mm>*{};<0mm,-3.0mm>*{}**@{-},
 <0.49mm,0.49mm>*{};<1.9mm,1.9mm>*{}**@{-},
 <-0.5mm,0.5mm>*{};<-1.9mm,1.9mm>*{}**@{-},
 <-2.3mm,2.3mm>*{\circ};<-2.3mm,2.3mm>*{}**@{},
 <-1.8mm,2.8mm>*{};<0mm,4.9mm>*{}**@{-},
 <-2.8mm,2.9mm>*{};<-4.6mm,4.9mm>*{}**@{-},
   <0.49mm,0.49mm>*{};<3.2mm,3.2mm>*{^3}**@{},
   <-1.8mm,2.8mm>*{};<0.4mm,6.3mm>*{^2}**@{},
   <-2.8mm,2.9mm>*{};<-5.1mm,6.3mm>*{^1}**@{},
 \end{xy}+
\begin{xy}
 <0mm,0mm>*{\circ};<0mm,0mm>*{}**@{},
 <0mm,-0.49mm>*{};<0mm,-3.0mm>*{}**@{-},
 <0.49mm,0.49mm>*{};<1.9mm,1.9mm>*{}**@{-},
 <-0.5mm,0.5mm>*{};<-1.9mm,1.9mm>*{}**@{-},
 <-2.3mm,2.3mm>*{\circ};<-2.3mm,2.3mm>*{}**@{},
 <-1.8mm,2.8mm>*{};<0mm,4.9mm>*{}**@{-},
 <-2.8mm,2.9mm>*{};<-4.6mm,4.9mm>*{}**@{-},
   <0.49mm,0.49mm>*{};<3.2mm,3.2mm>*{^2}**@{},
   <-1.8mm,2.8mm>*{};<0.4mm,6.3mm>*{^1}**@{},
   <-2.8mm,2.9mm>*{};<-5.1mm,6.3mm>*{^3}**@{},
 \end{xy}
 +
\begin{xy}
 <0mm,0mm>*{\circ};<0mm,0mm>*{}**@{},
 <0mm,-0.49mm>*{};<0mm,-3.0mm>*{}**@{-},
 <0.49mm,0.49mm>*{};<1.9mm,1.9mm>*{}**@{-},
 <-0.5mm,0.5mm>*{};<-1.9mm,1.9mm>*{}**@{-},
 <-2.3mm,2.3mm>*{\circ};<-2.3mm,2.3mm>*{}**@{},
 <-1.8mm,2.8mm>*{};<0mm,4.9mm>*{}**@{-},
 <-2.8mm,2.9mm>*{};<-4.6mm,4.9mm>*{}**@{-},
   <0.49mm,0.49mm>*{};<3.2mm,3.2mm>*{^1}**@{},
   <-1.8mm,2.8mm>*{};<0.4mm,6.3mm>*{^3}**@{},
   <-2.8mm,2.9mm>*{};<-5.1mm,6.3mm>*{^2}**@{},
 \end{xy}
 = 0, \ \ \ \ \ \ \
 \begin{xy}
 <0mm,2.47mm>*{};<0mm,-0.5mm>*{}**@{-},
 <0.5mm,3.5mm>*{};<2.2mm,5.2mm>*{}**@{-},
 <-0.48mm,3.48mm>*{};<-2.2mm,5.2mm>*{}**@{-},
 <0mm,3mm>*{\circ};
  <0mm,-0.8mm>*{\bullet};
<0mm,-0.8mm>*{};<-2.2mm,-3.5mm>*{}**@{-},
 <0mm,-0.8mm>*{};<2.2mm,-3.5mm>*{}**@{-},
     <0.5mm,3.5mm>*{};<2.8mm,6.7mm>*{^2}**@{},
     <-0.48mm,3.48mm>*{};<-2.8mm,6.7mm>*{^1}**@{},
   <0mm,-0.8mm>*{};<-2.7mm,-6.2mm>*{^1}**@{},
   <0mm,-0.8mm>*{};<2.7mm,-6.2mm>*{^2}**@{},
\end{xy}
 -
\begin{xy}
 <0mm,-1.3mm>*{};<0mm,-3.5mm>*{}**@{-},
 <0.38mm,-0.2mm>*{};<2.2mm,2.2mm>*{}**@{-},
 <-0.38mm,-0.2mm>*{};<-2.2mm,2.2mm>*{}**@{-},
<0mm,-0.8mm>*{\circ};
 <2.4mm,2.4mm>*{\bullet};
 <2.5mm,2.3mm>*{};<4.4mm,-0.8mm>*{}**@{-},
 <2.4mm,2.5mm>*{};<2.4mm,5.2mm>*{}**@{-},
     <0mm,-1.3mm>*{};<0mm,-6.3mm>*{^1}**@{},
     <2.5mm,2.3mm>*{};<5.1mm,-3.6mm>*{^2}**@{},
    <2.4mm,2.5mm>*{};<2.4mm,6.7mm>*{^2}**@{},
    <-0.38mm,-0.2mm>*{};<-2.8mm,3.5mm>*{^1}**@{},
    \end{xy}
  +
\begin{xy}
 <0mm,-1.3mm>*{};<0mm,-3.5mm>*{}**@{-},
 <0.38mm,-0.2mm>*{};<2.2mm,2.2mm>*{}**@{-},
 <-0.38mm,-0.2mm>*{};<-2.2mm,2.2mm>*{}**@{-},
<0mm,-0.8mm>*{\circ};
 <2.4mm,2.4mm>*{\bullet};
 <2.5mm,2.3mm>*{};<4.4mm,-0.8mm>*{}**@{-},
 <2.4mm,2.5mm>*{};<2.4mm,5.2mm>*{}**@{-},
     <0mm,-1.3mm>*{};<0mm,-6.3mm>*{^1}**@{},
     <2.5mm,2.3mm>*{};<5.1mm,-3.6mm>*{^2}**@{},
    <2.4mm,2.5mm>*{};<2.4mm,6.7mm>*{^1}**@{},
    <-0.38mm,-0.2mm>*{};<-2.8mm,3.5mm>*{^2}**@{},
    \end{xy}
 -
\begin{xy}
 <0mm,-1.3mm>*{};<0mm,-3.5mm>*{}**@{-},
 <0.38mm,-0.2mm>*{};<2.2mm,2.2mm>*{}**@{-},
 <-0.38mm,-0.2mm>*{};<-2.2mm,2.2mm>*{}**@{-},
<0mm,-0.8mm>*{\circ};
 <2.4mm,2.4mm>*{\bullet};
 <2.5mm,2.3mm>*{};<4.4mm,-0.8mm>*{}**@{-},
 <2.4mm,2.5mm>*{};<2.4mm,5.2mm>*{}**@{-},
     <0mm,-1.3mm>*{};<0mm,-6.3mm>*{^2}**@{},
     <2.5mm,2.3mm>*{};<5.1mm,-3.6mm>*{^1}**@{},
    <2.4mm,2.5mm>*{};<2.4mm,6.7mm>*{^2}**@{},
    <-0.38mm,-0.2mm>*{};<-2.8mm,3.5mm>*{^1}**@{},
    \end{xy}
 +
\begin{xy}
 <0mm,-1.3mm>*{};<0mm,-3.5mm>*{}**@{-},
 <0.38mm,-0.2mm>*{};<2.2mm,2.2mm>*{}**@{-},
 <-0.38mm,-0.2mm>*{};<-2.2mm,2.2mm>*{}**@{-},
<0mm,-0.8mm>*{\circ};
 <2.4mm,2.4mm>*{\bullet};
 <2.5mm,2.3mm>*{};<4.4mm,-0.8mm>*{}**@{-},
 <2.4mm,2.5mm>*{};<2.4mm,5.2mm>*{}**@{-},
     <0mm,-1.3mm>*{};<0mm,-6.3mm>*{^2}**@{},
     <2.5mm,2.3mm>*{};<5.1mm,-3.6mm>*{^1}**@{},
    <2.4mm,2.5mm>*{};<2.4mm,6.7mm>*{^1}**@{},
    <-0.38mm,-0.2mm>*{};<-2.8mm,3.5mm>*{^2}**@{},
    \end{xy}  = 0.
\Eeq
Its minimal resolution, $\LB_\infty$, has been computed in \cite{Me1}.

\subsubsection{\bf Theorem}\label{2: LB_infty}
(i)  {\em  $\LB_\infty$ is a dg free prop, $\cF^\uparrow\langle X\rangle$,
generated by an $\bS$-bimodule},
\Beq\label{2: LB_infty_generators}
X(m,n)[-1]= \sgn_m\ot \id_n[m-2]=\mbox{span}\left \langle
\begin{xy}
 <0mm,0mm>*{\bullet};<0mm,0mm>*{}**@{},
 <0mm,0mm>*{};<-8mm,4mm>*{}**@{-},
 <0mm,0mm>*{};<-4mm,4mm>*{}**@{-},
 <0mm,0mm>*{};<0mm,3.5mm>*{\ldots}**@{},
 <0mm,0mm>*{};<4mm,4mm>*{}**@{-},
 <0mm,0mm>*{};<8mm,4mm>*{}**@{-},
   <0mm,0mm>*{};<-8mm,5mm>*{^1}**@{},
   <0mm,0mm>*{};<-4.5mm,5mm>*{^2}**@{},
   <0mm,0mm>*{};<8.2mm,5mm>*{^m}**@{},
 <0mm,0mm>*{};<-8mm,-4mm>*{}**@{-},
 <0mm,0mm>*{};<-4mm,-4mm>*{}**@{-},
 <0mm,0mm>*{};<0mm,-4mm>*{\ldots}**@{},
 <0mm,0mm>*{};<4mm,-4mm>*{}**@{-},
 <0mm,0mm>*{};<8mm,-4mm>*{}**@{-},
   <0mm,0mm>*{};<-8mm,-5.5mm>*{_1}**@{},
   <0mm,0mm>*{};<-4.5mm,-5.5mm>*{_2}**@{},
   <0mm,0mm>*{};<8.5mm,-5.5mm>*{_n}**@{},
 \end{xy}
\right\rangle_{m\geq1,n\geq 1, m+n\geq3}
\vspace{-1mm}
\Eeq
{\em and with the differential given on the generators as follows},\vspace{-3mm}
$$
\delta
\begin{xy}
 <0mm,0mm>*{\bullet};<0mm,0mm>*{}**@{},
 <0mm,0mm>*{};<-8mm,4mm>*{}**@{-},
 <0mm,0mm>*{};<-4mm,4mm>*{}**@{-},
 <0mm,0mm>*{};<0mm,3.5mm>*{\ldots}**@{},
 <0mm,0mm>*{};<4mm,4mm>*{}**@{-},
 <0mm,0mm>*{};<8mm,4mm>*{}**@{-},
   <0mm,0mm>*{};<-8mm,5mm>*{^1}**@{},
   <0mm,0mm>*{};<-4.5mm,5mm>*{^2}**@{},
   <0mm,0mm>*{};<8.2mm,5mm>*{^m}**@{},
 <0mm,0mm>*{};<-8mm,-4mm>*{}**@{-},
 <0mm,0mm>*{};<-4mm,-4mm>*{}**@{-},
 <0mm,0mm>*{};<0mm,-4mm>*{\ldots}**@{},
 <0mm,0mm>*{};<4mm,-4mm>*{}**@{-},
 <0mm,0mm>*{};<8mm,-4mm>*{}**@{-},
   <0mm,0mm>*{};<-8mm,-5.5mm>*{_1}**@{},
   <0mm,0mm>*{};<-4.5mm,-5.5mm>*{_2}**@{},
   <0mm,0mm>*{};<8.5mm,-5.5mm>*{_n}**@{},
 \end{xy}
 =
 \sum_{[m]=I_1\sqcup I_2 , [n]=J_1\sqcup J_2\atop
 {|I_1|\geq 0, |I_2|\geq 1,
 |J_1|\geq 1, |J_2|\geq 0}
}\hspace{0mm} (-1)^{\sigma(I_1\sqcup I_2) + |I_1|(|I_2|+1)}
 \begin{xy}
 <0mm,0mm>*{\bullet};
 <0mm,0mm>*{};<-7mm,4mm>*{}**@{-},
 <0mm,0mm>*{};<-4mm,4mm>*{}**@{-},
 <0mm,0mm>*{};<0mm,3.5mm>*{\ldots}**@{},
 <0mm,0mm>*{};<4mm,4mm>*{}**@{-},
 <0mm,0mm>*{};<13mm,5mm>*{}**@{-},
     <0mm,0mm>*{};<-1mm,6mm>*{\overbrace{\ \ \ \ \ \ \ \ \ \ \ }}**@{},
     <0mm,0mm>*{};<-1mm,8mm>*{^{I_1}}**@{},
 <0mm,0mm>*{};<-7mm,-4mm>*{}**@{-},
 <0mm,0mm>*{};<-4mm,-4mm>*{}**@{-},
 <0mm,0mm>*{};<0mm,-3.5mm>*{\ldots}**@{},
 <0mm,0mm>*{};<4mm,-4mm>*{}**@{-},
 <0mm,0mm>*{};<7mm,-4mm>*{}**@{-},
      <0mm,0mm>*{};<0mm,-5.5mm>*{\underbrace{\ \ \ \ \ \ \ \ \ \ \ \ \ \
      }}**@{},
      <0mm,0mm>*{};<0mm,-8.6mm>*{_{J_1}}**@{},
 <13mm,5mm>*{\bullet}**@{},
 <13mm,5mm>*{};<6mm,9mm>*{}**@{-},
 <13mm,5mm>*{};<9mm,9mm>*{}**@{-},
 <13mm,5mm>*{};<13mm,8.5mm>*{\ldots}**@{},
 <13mm,5mm>*{};<17mm,9mm>*{}**@{-},
 <13mm,5mm>*{};<20mm,9mm>*{}**@{-},
      <13mm,5mm>*{};<13mm,10.7mm>*{\overbrace{\ \ \ \ \ \ \ \ \ \ \ \  }}**@{},
      <13mm,5mm>*{};<13mm,13mm>*{^{I_2}}**@{},
 <13mm,5mm>*{};<8mm,0mm>*{}**@{-},
 <13mm,5mm>*{};<12mm,0mm>*{\ldots}**@{},
 <13mm,5mm>*{};<16.5mm,0mm>*{}**@{-},
 <13mm,5mm>*{};<20mm,0mm>*{}**@{-},
     <13mm,5mm>*{};<14.3mm,-1.7mm>*{\underbrace{\ \ \ \ \ \ \ \ \ \ \ }}**@{},
     <13mm,5mm>*{};<14.3mm,-4.2mm>*{_{J_2}}**@{},
 \end{xy}
$$

(ii) {\em For any $d\in \N$, there is a one-to-one correspondence between representations
of the dg prop $\LB_\infty$ in $\R^d$ and formal Poisson structures, $\pi$,
on $\R^d$ vanishing at the origin}.

\begin{proof}
The proof of (i) is straightforward (see, e.g., \cite{MV,Me3,V}) once one uses rather
non-straightforward Koszul duality theory 
for dioperads, \cite{GK, G}, and Kontsevich's ideas of $\frac{1}{2}$-props and path filtrations
\cite{Ko0, MaVo}. We shall discuss some of these
ideas in \S 5 and show here now only the proof of (ii). Since $\R^p$ is concentrated in degree zero,
an arbitrary representation $\rho: \LB_\infty\rar \cE nd_{\R^p}$ can have non-zero values
only on $(m,n)$-corollas with $m=2$. Denote these values,
$
\rho\left(
\begin{xy}
 <0mm,0mm>*{\bullet};
 <0mm,0mm>*{};<-3mm,3mm>*{}**@{-},
 <0mm,0mm>*{};<3mm,3mm>*{}**@{-},
   <0mm,0mm>*{};<-3mm,4mm>*{^1}**@{},
   <0mm,0mm>*{};<3.5mm,4mm>*{^2}**@{},
 <0mm,0mm>*{};<-8mm,-4mm>*{}**@{-},
 <0mm,0mm>*{};<-4mm,-4mm>*{}**@{-},
 <0mm,0mm>*{};<0mm,-4mm>*{\ldots}**@{},
 <0mm,0mm>*{};<4mm,-4mm>*{}**@{-},
 <0mm,0mm>*{};<8mm,-4mm>*{}**@{-},
   <0mm,0mm>*{};<-8mm,-5.5mm>*{_1}**@{},
   <0mm,0mm>*{};<-4.5mm,-5.5mm>*{_2}**@{},
   <0mm,0mm>*{};<8.5mm,-5.5mm>*{_n}**@{},
 \end{xy}
 \right)\in \Hom(\odot^n \R^p, \wedge^2 \R^p),
$
by $\pi_n$. As the tangent space,  $\cT_{0}$, to $\R^p$ at zero can be identified with $\R^p$ itself,
we can identify
the total sum $\pi:= \sum_{n\geq 1}\pi_n \in \Hom(\odot^{\geq 1}\R^p, \wedge^2 \cT_0)$ with a formal
bi-vector field on $\R^p$. Then the equation $\rho\circ \delta=\delta \circ \rho$ becomes precisely
the Poisson equation, $[\pi,\pi]_S=0$, where $[\ ,\ ]_S$ is the Schouten bracket.
\end{proof}
It is worth pointing out that the vanishing condition $\pi|_{0\in \R^p}=0$ in Theorem~\ref{2: LB_infty}(ii)
is no serious restriction:
 given an arbitrary formal or analytic Poisson
structure $\pi$ on $\R^p$ (not necessary  vanishing at $0\in \R^p$), then, for any parameter $\lambda$
viewed as a coordinate
on $\R$, the product
$\lambda\pi$ is a Poisson structure on $\R^{p+1}=\R^p\times \R$ vanishing at zero $0\in \R^{n+1}$ and hence
is a representation of the prop $\LB_\infty$.

\subsubsection{\bf Bi-Hamiltonian geometry}
The prop profile of a {\em pair of compatible Poisson structures}\, (which is an important concept
in the theory of integrable systems) has been computed by Strohmayer in
\cite{Str2} with the help of an earlier
result of  Dotsenko and Khoroshkin \cite{DK}.

\subsubsection{\bf Wheeled Poisson structures?} Theorem~\ref{2: LB_infty} says that the minimal resolution,
$$
\LB_\infty=(\cF^\uparrow\langle X \rangle, \delta),
$$
of the prop, $\LB$, of arbitrary Lie 1-bialgebras controls  the category of local (formal)
smooth Poisson structures.
What can be said about a minimal resolution, $(\LB^\circlearrowright)_\infty$,
 of the wheeled prop, $\LB^\circlearrowright$,
 of {\em finite dimensional}\, Lie 1-bialgebras whose representations can, in view of
 Theorem~\ref{2: LB_infty}(ii), be called  {\em wheeled Poisson structures}?
 Note that
 $\LB^\circlearrowright$ has the same generators
 and relations as $\LB$, the only difference being that graphs now might have wheels. As in the case
 of associative algebra, the naive wheelification,
 $$
 (\LB_\infty)^\circlearrowright:=(\cF^\circlearrowright\langle X \rangle, \delta),
 $$
creates  new {\em non-trivial}\, cohomology classes, as e.g.\ this one \cite{Me3}
\vspace{-3mm}
\Beq\label{2: Three_wheels}
 \begin{xy}
<-5mm,5mm>*{\bullet};
<-5mm,5mm>*{};<-5mm,8mm>*{}**@{-},
<-5mm,5mm>*{};<-7mm,4mm>*{}**@{-},
 <0mm,0mm>*{\bullet};
<0mm,0mm>*{};<-5mm,-5mm>*{}**@{-},
 <0mm,0mm>*{};<-5mm,5mm>*{}**@{-},
<0mm,0mm>*{};<1.5mm,1.5mm>*{}**@{-},
 <0mm,0mm>*{};<1.5mm,-1.5mm>*{}**@{-};
<-5mm,-5mm>*{\bullet};
<-5mm,-5mm>*{};<-7mm,-2mm>*{}**@{-},
<-5mm,-5mm>*{};<-5mm,-8mm>*{}**@{-},
   \ar@{->}@(u,dl) (-5.0,8.0)*{};(-7.0,4.0)*{},
   \ar@{->}@(ur,dr) (1.5,1.5)*{};(1.5,-1.5)*{},
   \ar@{->}@(ul,d) (-7.0,-2.0)*{};(-5.0,-8.0)*{},
\end{xy}
\ - \
 \begin{xy}
<0mm,4mm>*{\bullet};
<0mm,4mm>*{};<0mm,8mm>*{}**@{-},
<0mm,4mm>*{};<3mm,2mm>*{}**@{-},
 <0mm,-5mm>*{\bullet};
<0mm,-5mm>*{};<0mm,4mm>*{}**@{-},
<0mm,-5mm>*{};<2mm,-3mm>*{}**@{-},
<0mm,-5mm>*{};<0mm,-8mm>*{}**@{-},
<-5mm,-5mm>*{};<0mm,4mm>*{}**@{-};
<-5mm,-5mm>*{\bullet};
<-5mm,-5mm>*{};<-7mm,-2mm>*{}**@{-},
<-5mm,-5mm>*{};<-5mm,-8mm>*{}**@{-},
   \ar@{->}@(u,dr) (0,8.0)*{};(3.0,2.0)*{},
   \ar@{->}@(ur,d) (2.0,-3.0)*{};(0.0,-8.0)*{},
   \ar@{->}@(ul,d) (-7.0,-2.0)*{};(-5.0,-8.0)*{},
\end{xy}
\ + \
 \begin{xy}
<0mm,-4mm>*{\bullet};
<0mm,-4mm>*{};<0mm,-8mm>*{}**@{-},
<0mm,-4mm>*{};<3mm,-2mm>*{}**@{-},
 <0mm,5mm>*{\bullet};
<0mm,5mm>*{};<0mm,-4mm>*{}**@{-},
<0mm,5mm>*{};<2mm,3mm>*{}**@{-},
<0mm,5mm>*{};<0mm,8mm>*{}**@{-},
<-5mm,5mm>*{};<0mm,-4mm>*{}**@{-};
<-5mm,5mm>*{\bullet};
<-5mm,5mm>*{};<-7mm,2mm>*{}**@{-},
<-5mm,5mm>*{};<-5mm,8mm>*{}**@{-},
   \ar@{->}@(ur,d) (3.0,-2.0)*{};(0.0,-7.0)*{},
   \ar@{->}@(u,dr) (0.0,8.0)*{};(2.0,3.0)*{},
   \ar@{->}@(u,dl) (-5.0,8.0)*{};(-7.0,2.0)*{},
\end{xy}
\ \ \ \in (\LB_\infty)^\circlearrowright\vspace{-2mm}
\Eeq
which map under the natural projection $(\LB_\infty)^\circlearrowright \rar \LB^\circlearrowright$
to zero. Thus the set of generators of a minimal resolution, $(\LB^\circlearrowright)_\infty$,
of $\LB^\circlearrowright$
must be larger than the set  (\ref{2: LB_infty_generators}), and at present its computation is beyond reach.
 All we can say now
about mysterious {\em wheeled}\, Poisson structures on a graded formal manifold $M$ is that
(i) they are
Maurer-Cartan elements of a certain $L_\infty$-algebra extension of the ordinary Schouten
bracket on $M$ which involves divergence operators (in fact, graph (\ref{2: Three_wheels}) gives us a glimpse
of the $\mu_3$ composition in that $L_\infty$-algebra), and
(ii) they can be deformation quantized in exactly the same sense as  ordinary Poisson
structures; moreover, it is proven in \cite{Me4} with the help of the theory of wheeled props that there exist universal formulae
for deformation quantization of wheeled Poisson structures which involve only rational numbers $\Q$.


\subsection{Pre-Lie algebras, Nijenhuis geometry and contractible dg manifolds}
A {\em pre-Lie}\, algebra
is a vector space together with a binary operation, $\circ: V^{\ot
2}\rar V$, satisfying the condition
$$
(a\circ b)\circ c - a\circ (b\circ c) - (-1)^{|b||c|}(a\circ
c)\circ b + (-1)^{|b||c|} a\circ (c\circ b) =0
$$
for any $a,b,c\in V$.
Any pre-Lie algebra is naturally a Lie algebra with the bracket,
$[a,b]:=a\circ b  - (-1)^{|a||b|}b\circ a$. Let us
consider the following extension of this notion:
a {\em pre-Lie}$^{\mathbf\frak 2}$ algebra is a pre-Lie algebra $(V,\circ)$
equipped with a compatible Lie bracket in  degree 1, i.e.\ with a linear map\footnote{Equivalently,
a linear map $[\, \bullet\, ]:  \odot^2 V \rar V[1]$.}
$[\, \bullet\, ]:  \wedge^2 (V[-1]) \rar V[-1]$  satisfying the Jacobi identities
and the following compatibility condition,
$$
[a\bullet b]\circ c + (-1)^{|b|}a\circ [b\bullet c] +
(-1)^{|b||a|+|b|} b\circ [a\bullet c] = \hspace{7cm}
$$
$$
 \hspace{4cm}
 (-1)^{|b||c|+|c|} [(a\circ c)\bullet b]
+ (-1)^{(|a|+1)(|b|+|c|)+|a|} [(b\circ c)\bullet a], \ \ \forall\, a,b,c\in V
$$

\sip
This compatibility condition can be understood as follows. The vector space
$V\oplus V[-1]$ is naturally a complex with trivial cohomology. If
we write elements of $V\oplus V[-1]$ as $a + \Pi b$, where $a,b\in
V$ and $\Pi$ is a formal symbol of degree 1,
 then the natural differential in $V\oplus V[-1]$ is given by
$d(a + \Pi b)= 0 + \Pi a$. Given two arbitrary binary operations,
$$
\circ: V\ot V \rar V, \ \ \
[\, \bullet\, ]:  \odot^2 V \rar V[1],
$$
define a degree zero map, $[\ ,\ ]: \wedge^2(V\oplus V[-1])\rar V\oplus V[-1]$ by setting,
$$
[a,b]:= a\circ b - (-1)^{|a||b|} b\circ a,\ \ \
\left[\Pi a, b\right] := - (-1)^{|a|}[a\bullet b] +\Pi a\circ b,\ \ \
\left[\Pi a, \Pi b\right] :=\Pi\left[a\bullet b\right].
$$
\subsubsection{\bf Proposition \cite{Me2}}
\label{2: contractible dg Lie}
 {\em The data $(V\oplus V[-1], d, [\ ,\ ])$ is a
(contractible) dg Lie algebra if and only if $(V,\circ, [\ \bullet\ ])$ is a pre-Lie$^2$ algebra.}

\sip
Rather surprisingly, the minimal resolution, $pre$-$\caL ie^2_\infty$, of the operad of pre-Lie$^2$-algebras
 has much to do with
the famous Nijenhuis integrability condition in differential geometry. The following result
is based on the works \cite{CL, Me2, Str2}.

\subsubsection{\bf Theorem} (i) {\em The operad $pre$-$\caL ie^2_\infty$ is a free operad, $\cF^\curlywedge\langle N\rangle$,
generated by an $\bS$-bimodule $N$ with all $N(m,n)=0$ except the following ones},
$$
N(1,n):=
\bigoplus_{p=1}^{n} {\mathrm
I\mathrm n\mathrm d}^{\bS_n}_{\bS_{p}\times \bS_{n-p}} {\id}_{p}\ot
\sgn_{n-p}[n-p-1]=\mbox{span}
\left\langle \hspace{-1mm}
\Ba{c}
\begin{xy}
 <0mm,-0.5mm>*{\blacktriangledown};
 <0mm,0mm>*{};<0mm,4mm>*{}**@{-},
 <0mm,0mm>*{};<-16mm,-5mm>*{}**@{-},
 <0mm,0mm>*{};<-11mm,-5mm>*{}**@{-},
 <0mm,0mm>*{};<-3.5mm,-5mm>*{}**@{-},
 <0mm,0mm>*{};<-6mm,-5mm>*{...}**@{},
   <0mm,0mm>*{};<-15mm,-7.5mm>*{^{i_1}}**@{},
   <0mm,0mm>*{};<-11mm,-7.5mm>*{^{i_2}}**@{},
   <0mm,0mm>*{};<-3mm,-7.5mm>*{^{i_p}}**@{},
 <-9.7mm,-9.4mm>*{\underbrace{\ \ \ \ \ \ \ \ \ \ \ \ \  }}**@{},
  <-10mm,-12mm>*{_{symmetric}}**@{},
 <0mm,0mm>*{};<16mm,-5mm>*{}**@{-},
 <0mm,0mm>*{};<12mm,-5mm>*{}**@{-},
 <0mm,0mm>*{};<3.5mm,-5mm>*{}**@{-},
 <0mm,0mm>*{};<6.6mm,-5mm>*{...}**@{},
   <0mm,0mm>*{};<17mm,-7.5mm>*{^{i_n}}**@{},
   <0mm,0mm>*{};<6mm,-7.5mm>*{^{i_{p+1}}}**@{},
 <9.7mm,-9.4mm>*{\underbrace{\ \ \ \ \ \ \ \ \ \ \ \ \ \ }}**@{},
  <10mm,-12mm>*{_{skewsymmetric}}**@{},
 \end{xy}
\Ea
\hspace{-1.5mm}
\right\rangle,  \ n\geq 2,
$$
{\em and equipped with a differential given on the generators by}
$$
d\Ba{c}
 \begin{xy}
<18mm,-0.5mm>*{{\blacktriangledown}};
<18mm,0cm>*{};<18mm,5mm>*{}**@{-},
<18mm,0cm>*{};<0mm,-7mm>*{}**@{-},
 <18mm,0cm>*{};<5mm,-7mm>*{}**@{-},
 <18mm,0cm>*{};<11mm,-6mm>*{}**@{-},
<11mm,-7mm>*{{\ldots}};
 <18mm,0cm>*{};<15mm,-7mm>*{}**@{-},
 <18mm,0cm>*{};<21mm,-7mm>*{}**@{-},
 <18mm,0cm>*{};<26mm,-7mm>*{}**@{-},
 <18mm,0cm>*{};<30mm,-6mm>*{}**@{-},
<31mm,-7mm>*{{\ldots}};
 <18mm,0cm>*{};<36mm,-7mm>*{}**@{-},
 <18mm,0cm>*{};<-0.5mm,-10mm>*{^1}**@{},
 <18mm,0cm>*{};<4.5mm,-10mm>*{^2}**@{},
 <18mm,0cm>*{};<13.5mm,-10mm>*{^{{p}}}**@{},
 <18mm,0cm>*{};<22.5mm,-10mm>*{^{{p\hspace{-0.4mm}
 +\hspace{-0.5mm}1}}}**@{},
 <18mm,0cm>*{};<36.5mm,-10mm>*{^{n}}**@{},
 \end{xy}
 \Ea
 =
 \sum_{I_1\sqcup I_2=(1,\ldots,p) \atop
{J_1\sqcup J_2=(p+1,\ldots,n) \atop
 {\#I_2\geq 1, \#I_1+\#J_2\geq 1\atop
 { \#I_2+\#J_1\geq 2}}}}
 \hspace{-5mm} (-1)^{\#J_2 + \sigma(J_1,J_2)}
\begin{xy}
<28.4mm,-0.5mm>*{{\blacktriangledown}};
 <28.4mm,0cm>*{};<28.4mm,5mm>*{}**@{-},
 <28mm,0cm>*{};<7mm,-7mm>*{}**@{-},
 <28mm,0cm>*{};<14mm,-7mm>*{\ldots}**@{-},
 <28mm,0cm>*{};<20mm,-7mm>*{}**@{-},
 <28mm,0cm>*{};<27mm,-7mm>*{}**@{-},
 <28mm,0cm>*{};<34mm,-7mm>*{}**@{-},
 <28mm,0cm>*{};<40mm,-7mm>*{\ldots}**@{-},
 <28mm,0cm>*{};<48mm,-7mm>*{}**@{-},
 <28mm,0cm>*{};<42mm,-11.5mm>*{^{J_2}}**@{},
 <28mm,0cm>*{};<41mm,-9mm>*{\underbrace{\ \ \ \ \ \  \ \ \ \ \ \  }}**@{},
<27mm,-7.5mm>*{{\blacktriangledown}};
 <27mm,-7mm>*{};<14mm,-14mm>*{}**@{-},
 <27mm,-7mm>*{};<19mm,-14mm>*{\cdots}**@{-},
 <27mm,-7mm>*{};<24mm,-14mm>*{}**@{-},
 <27mm,-7mm>*{};<30mm,-14mm>*{}**@{-},
 <27mm,-7mm>*{};<35mm,-14mm>*{\cdots}**@{-},
 <27mm,-7mm>*{};<40mm,-14mm>*{}**@{-},
 <38mm,-0cm>*{};<35mm,-19mm>*{^{J_1}}**@{},
 <38mm,0cm>*{};<35mm,-16mm>*{\underbrace{\ \ \ \ \ \ \ \ \   }}**@{},
<38mm,0cm>*{};<19mm,-19mm>*{^{I_2}}**@{},
 <38mm,0cm>*{};<19mm,-16mm>*{\underbrace{\ \ \  \ \ \ \ \ \ \   }}**@{},
 <38mm,0cm>*{};<13mm,-11.5mm>*{^{I_1}}**@{},
 <38mm,0cm>*{};<13mm,-9mm>*{\underbrace{\ \ \ \ \ \ \ \ \ \ \ \  }}**@{},
 \end{xy}
$$
$$
\hspace{30mm}
 -\sum_{I_1\sqcup I_2=(1,\ldots,n-p) \atop
{J_1\sqcup J_2\sqcup J_3=(n-p+1,\ldots,n) \atop
 {\#I_1\geq 1, \#I_2\geq 1\atop
 { \#I_1+\#J_3\geq 1 ,\#I_2+\#J_2\geq 1}}}}
 \hspace{-5mm} (-1)^{\#J_2 + \#J_3 +\sigma(J_1,J_2,J_3)}
  \hspace{-2mm}
 \begin{xy}
<28mm,-0.5mm>*{{\blacktriangledown}};
 <28mm,0cm>*{};<28mm,5mm>*{}**@{-},
 <28mm,0cm>*{};<11mm,-7mm>*{}**@{-},
 <28mm,0cm>*{};<18mm,-7mm>*{\ldots}**@{-},
 <28mm,0cm>*{};<24mm,-7mm>*{}**@{-},
 <28mm,0cm>*{};<31mm,-7mm>*{}**@{-},
 <28mm,0cm>*{};<38mm,-7mm>*{}**@{-},
 <28mm,0cm>*{};<44mm,-7mm>*{\ldots}**@{-},
 <28mm,0cm>*{};<51mm,-7mm>*{}**@{-},
 <28mm,0cm>*{};<46mm,-11.5mm>*{^{J_3}}**@{},
 <28mm,0cm>*{};<45mm,-9mm>*{\underbrace{\ \ \ \ \ \  \ \ \ \ \ \ \ }}**@{},
<31mm,-7.5mm>*{{\blacktriangledown}};
 <31mm,-7mm>*{};<17mm,-14mm>*{}**@{-},
 <31mm,-7mm>*{};<22mm,-14mm>*{\cdots}**@{-},
 <31mm,-7mm>*{};<27mm,-14mm>*{}**@{-},
 <31mm,-7mm>*{};<35mm,-14mm>*{}**@{-},
 <31mm,-7mm>*{};<40mm,-14mm>*{\cdots}**@{-},
 <31mm,-7mm>*{};<45mm,-14mm>*{}**@{-},
<31mm,-7mm>*{};<31.7mm,-16.5mm>*{^{J_1}}**@{},
<31mm,-7mm>*{};<31mm,-14mm>*{}**@{-},
 <42mm,-0cm>*{};<40mm,-19mm>*{^{J_2}}**@{},
 <42mm,0cm>*{};<40mm,-16mm>*{\underbrace{\ \ \ \ \ \ \ \ \  }}**@{},
<42mm,0cm>*{};<22mm,-19mm>*{^{I_2}}**@{},
 <42mm,0cm>*{};<22mm,-16mm>*{\underbrace{\ \ \  \ \ \ \ \ \  }}**@{},
 <42mm,0cm>*{};<17mm,-11.5mm>*{^{I_1}}**@{},
 <42mm,0cm>*{};<17mm,-9mm>*{\underbrace{\ \ \ \ \ \ \ \ \ \ \ \  }}**@{},
 \end{xy}
$$
(ii) {\em For any $d\in \N$,
there is a one-to-one correspondence between representations of $pre$-$\caL ie^2_\infty$ in $\R^d$
and endomorphisms, $J: \cT_{\R^d}\rar \cT_{\R^d}$, of the tangent bundle on the affine space $\R^d$
satisfying the Nijenhuis integrability condition, $N_J=0$, and the vanishing condition $J\mid_{0\in \R^d}=0$.
}

\mip

We recall that the Nijenhuis
tensor of an endomorphism, $J:\cT_M \rar \cT_M$, of the tangent bundle of an arbitrary  smooth manifold
$M$ (in particular, of $\R^m$) can be defined as a map
$$
\Ba{rccl} N_J: & \wedge^2 \cT_{\R^m} & \lon & \cT_{\R^m} \\
               & X\ot Y       & \lon & N_J(X,Y):=
               [JX,JY] + J^2[X,Y] - J[X,JY] - J[JX,Y],
\Ea
$$
and that its beauty is hidden in the far from being obvious fact that it is linear not only over $\R$
but also over arbitrary smooth functions, $f\in \f_M$, on $M$, that is, $N_J(fX,Y)=N_J(X,fY)=fN_J(X,Y)$.

\sip

A representation of this dg operad in an arbitrary {\em graded}\, vector space $V$
 might be called a
{\em graded}\, or
{\em extended}\, Nijenhuis structure on $V$ (viewed as a formal manifold).
Interestingly, the category of these {\em extended Nijenhuis manifolds}\, is almost identical
(see \cite{Me2})
to  the category of {\em contractible dg manifolds}\, which we first met
in \S\ref{2: subsub on Lie_infty} when discussing Kontsevich's homotopy classification of dg manifolds.
 Proposition~\ref{2: contractible dg Lie} above is in fact one of the simplest manifestations of this
 more general phenomenon.

\subsection{Gerstenhaber algebras and Hertling-Manin geometry}
\label{3: Gerst algebras and HM}
 We conclude this section with an example which was actually the first one
to reveal strong interconnections between derived (via minimal resolutions)
categories of rather simple algebraic structures and solution sets of highly non-linear
{\em diffeomorphism covariant}\, differential equations on ordinary smooth manifolds.
\sip

A {\em Gerstenhaber}\,
algebra  is, by definition, a graded vector space $V$ together with two
linear maps, $\circ: \odot^2 V  \rar V$ and $[\, \bullet\, ]: \odot^2 V \rar V[1]$
such that $(V,\circ)$ is a graded commutative algebra, $(V[-1], [\ \bullet\ ])$ is a graded
Lie algebra, and the compatibility equation,
 $$
 [(a\circ b)\bullet c]= a\circ [b\bullet c] +    (-1)^{|b|(|c|+1)}[a\bullet c]\circ b,
 \ \ \forall\ a,b,c\in V,
 $$
holds. The operad of Gerstenhaber algebras is often denoted by $\cG$. Its minimal resolution, $\cG_\infty$,
has been computed in \cite{GJ}; it is one of the most important operads in mathematics which found
many applications in homological  algebra, algebraic topology and deformation quantization.
It was shown in \cite{Me-F} that $\cG_\infty$
has also a differential geometric dimension:

\subsubsection{\bf Theorem}{\em For any  $d\in \N$,
there is a one-to-one correspondence between representations of the dg operad $\cG_\infty$ in $\R^d$
(concentrated in degree $0$)
and morphisms of sheaves, $\mu: \odot^2 \cT_{\R^d}\rar \cT_{\R^d}$, making the tangent sheaf, $\cT_{\R^d}$,
into a commutative and associative algebra and
satisfying the Hertling-Manin integrability condition, $R_\mu=0$, and the vanishing condition
$\mu\mid_{0\in \R^d}=0$.
}

\mip

We recall that the Hertling-Manin
tensor, $R_\mu$,  of an arbitrary commutative and associative product, $\mu:\cT_M\odot \cT_M \rar \cT_M$,
on the tangent sheaf of an arbitrary  smooth manifold
$M$ is a map \cite{HM}
$$
\Ba{rccl} R_\mu: & \ot^4 \cT_{M} & \lon & \cT_{M} \\
               & X\ot Y\ot Z\ot W       & \lon & R_\mu(X,Y,Z,W)
\Ea
$$
where
\Beqrn
R_\mu(X,Y,Z,W)&=&
  [\mu(X, Y), \mu(Z,W)] - \mu([\mu(X, Y), Z], W)
  - \mu(Z,[\mu(X,Y), W])
   \\
 && - \mu(X, [Y,\mu(Z, W)]) -  \mu[X,
 \mu(Z,W)], Y)  + \mu(X, \mu(Z, [Y,W]))\\
 &&  + \mu(X, \mu([Y,Z], W))
 + \mu([X,Z],\mu( Y, W)) + \mu([X,W],\mu(Y, Z)).
  \Eeqrn
A remarkable fact is that this map is linear not only over $\R$
but also over arbitrary smooth functions, $f\in \f_M$, on $M$, that is,
$R_\mu(fX,Y,Z,W)=fR_\mu(X,Y,Z,W)$, $R_\mu(X,fY,Z,W)=fR_\mu(X,Y,Z,W)$, etc.
One can view the Hertling-Manin integrability equation as a diffeomorphism covariant
version of the WDVV equation \cite{HM,HM2}.

\section{Applications to deformation theory}
\subsection{From minimal resolutions to $L_\infty$-algebras}
One of the advantages of knowing a  dg free  resolution, $\cP_\infty$, of a $\fG$-algebra
controlling
a mathematical structure $\cP$  is that $\cP_\infty$ paves a direct
way  to the deformation theory of $\cP$-structures.
In the heart of this approach to the deformation theory of many algebraic and geometric
structures
is  observation~\ref{3: Main Def Th} (see below) which was proven in \cite{MV} in several ways. For
its precise formulation we need the following notion.

\subsubsection{\bf Definitions}
 A $L_\infty$-algebra $(\fg, \{\mu_n:\wedge^n \fg\rar \fg[2-n]\}_{n\geq 1})$ is called
{\em filtered}\, if $\fg$ admits a non-negative decreasing
Hausdorff filtration,
$$
\fg_0=\fg \supseteq \fg_1\supseteq\ldots\supseteq \fg_i \supseteq
\ldots,
$$
such that $\Img \mu_n \subset \fg_{n}$ for all $n\geq n_0$ beginning with some
$n_0\in \N$. In this case it makes sense to define the associated set, $\cM\cC(\fg)$, of {\em Maurer-Cartan
elements}\, as a subset of $\fg$ consisting of degree 1 elements $\Ga$ satisfying the equation
$\sum_{n\geq 1}\frac{1}{n!}\mu_n(\Ga,\ldots, \Ga)=0$.
\mip

A very useful fact is that to
 every Maurer-Cartan element, $\Ga\in \cM\cC(\fg)$, of a filtered $L_\infty$-algebra
 $(\fg, \{\mu_n:\wedge^n \fg\rar \fg\}_{n\geq 1})$ there corresponds a $\Ga$-{\em twisted}\,
 $L_\infty$-algebra structure,  $\{\mu^\Ga_n:\wedge^n \fg\rar \fg\}_{n\geq 1}$, on $\fg$.
If one thinks of the original
$L_\infty$ algebra as of a dg manifold $(\cM_\fg,\eth)$ (see \S\ref{3: Main Def Th}),
then the set $\cM\cC(\fg)$ can be identified
with the zero set of the homological vector field $\eth$, and the $\Ga$-twisted
$\caL_\infty$-algebra structure on $\fg$ corresponds to
that homological vector field $\eth^\Gamma$ on $\cM_\fg$
which is obtained from $\eth$ by the translation diffeomorphism, $x\rar x+\Ga$, $\forall x\in \cM_\fg$.

\subsubsection{\bf Theorem \cite{MV}}\label{3: Main Def Th}
{\em Let $(\cF^\fG\langle E\rangle, \delta)$ be a dg free $\fG$-algebra (see Table 1)
 generated by an $\bS$-bimodule
$E$, and let $(\cQ, \delta_\cQ)$ be an arbitrary dg $\fG$-algebra. Then
\begin{itemize}
\item[(i)] the graded vector space,
$\fg:=\Hom_{\bS}(E,\cQ)[-1]$,  is canonically a filtered
$L_\infty$-algebra;
\item[(ii)] the set of all morphisms, $\{ \cF^\fG\langle E\rangle \rar \cQ\}$, of dg
$\fG$-algebras
is canonically isomorphic to the Maurer-Cartan set, $\cM\cC(\fg)$, of the  $L_\infty$-algebra
in (i).
\end{itemize}
}
\begin{proof} As an illustration we show an elementary proof of the theorem
 in the simplest
case $\fG=\fG^|$ (see Table 1), i.e.\ in the case when
$\cF^\fG\langle E\rangle$ is the free
associative algebra, $\ot^\bullet E$,
 generated by a graded vector space $E$ and $(\cQ, \delta_\cQ)$ is an arbitrary
dg associative algebra (we refer the reader to \cite{MV} for all other cases from
Table 1 except $\fG^\circlearrowright$, and to \cite{Gr2}
for the case $\fG^\circlearrowright$).
With these data
we shall associate a cohomological vector field, $\eth$, on the space
$\fg[1]=\Hom(E,\cQ)=\cQ\ot E^*$, and we shall do it in local coordinates by assuming further
 (only for simplicity of sign factors in formulae) that the graded vector spaces
 $E$ and  $\cQ$ are free modules
over some graded commutative ring, $R=\oplus_{i\in \Z} R^i$, with  {\em degree 0}\,
 generators $\{e_a\}_{a\in I}$
and, respectively, $\{e_\al\}_{\al\in J}$. Then the differentials in $\ot^\bullet E$
and $\cQ$, as well as multiplication $\circ$  in $\cQ$, have, respectively,
the following coordinate representations,
$$
\delta e_a=\sum_{k\geq 1 \atop a_1\ldots a_k\in I}
\delta_{a}^{a_1\ldots a_k} e_{a_1}\ot \ldots \ot e_{a_k},
\ \ \ \ \
\delta_\cQ e_\al=\sum_{\be\in J} Q_\al^\be e_\be,
\ \ \ \ \ e_\al\circ a_\be= \sum_{\ga\in J} \mu_{\al\be}^\ga e_\ga,
$$
for some coefficients
$\delta_a^{a_1\ldots a_k}\in R^1$, $Q_\al^\be\in R^{1}$ and  $\mu_{\al\be}^\ga\in R^0$.
The vector space of all $\cR$-linear maps, $\Hom(E,\cQ)$, is naturally graded,
$\Hom(E,\cQ)=\oplus_{i\in\Z} \Hom^i(E,\cQ)$, with $\Hom^i(E,\cQ)$ denoting the space
of all homogeneous linear maps of degree $i$. In the chosen bases a generic
element $\ga\in
 \Hom^i(E,\cQ)$ gets a coordinate representation,
$
\ga(e_a)=\sum_{\al\in J} \ga_{a (i)}^{\al} e_\al,
$
for some coefficients $\ga_{a (i)}^{\al} \in R^i$. The family of parameters
$\{ \ga_{a (i)}^{\al}\}_{a\in I,\al\in J, i\in \Z}$ provides us with a coordinate
system on the formal manifold $\cM_{\fg}\simeq \Hom(E,\cQ)$. In these coordinates
the required homological vector field on $\cM_{\fg}$, that is, a
$\caL_\infty$-structure on $\Hom(E,\cQ)[-1]$, is given explicitly by
$$
\eth =
\left(\sum_{\al,\be,a, i}Q_\be^\al \ga_{a (i)}^\be - \sum_{a, a_\bullet,\al,i}(-1)^i
\delta_a^{a_1 \ldots a_k} \ga_{a_1\ldots a_k (i)}^{\al}\right)
\frac{\p}{\p \ga_{a (i)}^{\al}}
$$
where, for $k\geq 2$,
$$
\ga_{a_1a_2\ldots a_k(i)}^\al = \sum_{\be_\bullet,\ga_\bullet \in
J\atop i_1+\ldots +i_k=i}
\mu_{\be_1\ga_1}^\al\mu_{\be_2\ga_2}^{\ga_1}\ldots \mu_{\be_{k-1}\be_k}^{\ga_{k-2}}
 \ga_{a_1(i_1)}^{\be_1}  \ga_{a_2(i_2)}^{\be_2}\ldots  \ga_{a_k(i_k)}^{\be_k}.
$$
The equation $[\eth,\eth]=0$ follows straightforwardly from the assumptions that
$\delta^2=0$, $\delta_\cQ^2=0$, as well as from the associativity of the product $\circ$ and its
compatibility with $\delta_\cQ$. This proves (i).

The Maurer-Cartan set $\cM\cC(\fg)$ is precisely  the set $\{\ga\in
\Hom^0(E,\cQ): \eth|_\ga=0\}$ and, therefore, consists of all
points in  $\Hom(E,\cQ)$ which have all the
 coordinates $\{\ga_{a (i)}^{\al}\}_{i\neq 0}$
vanishing, and the coordinate $\ga_{a (0)}^{\al}$ satisfying the equation,
$$
\sum_{\be\in J}Q_\be^\al \ga_{a (0)}^\be - \sum_{a_1,\ldots,a_k\in I}
\delta_a^{a_1 \ldots a_k} \ga_{a_1\ldots a_k (0)}^{\al}=0.
$$
Which just says that the associated to  $\ga_{a (0)}^{\al}$ map
of associative algebras, $\odot^\bullet E\rar \cQ$, commutes with the differentials
$\delta$ and
$\delta_\cQ$ defining thereby a morphism of {\em dg}\, algebras.
This proves claim (ii).
\end{proof}

\subsection{Deformation theory}
The theory of operads and props gives a universal approach to the deformation theory
of many algebraic and geometric structures and provides us with a conceptual explanation of
the well-known ``experimental" observation
that a deformation theory is controlled by  a differential graded Lie or, more generally, a $L_\infty$-algebra.
What happens is the following \cite{KS,MV, vdL}:
\Bi
\item[(I)] an algebraic or a (germ of) geometric structure, $\fs$, in
a vector space $V$ (which is an {\em object}\,
in the corresponding category, $\fS$, of algebraic or geometric structures) can
often be interpreted as a {\em representation}, $\al_\fs: \cS\rar \cE nd_V$,
of a $\fG$-algebra $\cS$ uniquely associated to the category
of $\fs$-structures;
\item[(II)] a dg resolution,
$\pi:\cS_\infty=(\cF^\fG\langle E\rangle, \delta)\rar \cS$,
 of the $\fG$-algebra $\cS$
gives rise, by Theorem~\ref{3: Main Def Th}, to a filtered $L_\infty$-algebra
on the vector space $\fg=\Hom_\bS(E,\cE nd_V)[-1])$ whose Maurer-Cartan
elements correspond to all possible representations, $\cS_\infty\rar \cE nd_V$;
in particular, our original algebraic or geometric structure $\fs$ defines
a Maurer-Cartan element  $\Ga_\fs:=\al_\fs\circ \pi$ in $\cM\cC(\fg)$;

\item[(III)] the $\Ga_\fs$-twisted $L_\infty$-algebra structure on $\fg$
is precisely the one which controls, in  Deligne's sense,
 the deformation theory of $\fs$.
\Ei
For example, if $\fs$ is the structure of associative algebra on a vector space $V$, then,
\Bi
\item[(i)] there is an operad, $\cA ss$, uniquely associated to the category
of associative algebras such that $\fs$ corresponds to a morphism,
$\al_\fs: \cA ss \rar \cE nd_V$, of operads (see \S 2.1);
\item[(ii)] there is a unique minimal resolution (see Theorem~\ref{2: Stasheff_theorem}),
$\cA ss_\infty$, of $\cA ss$ which is generated by the $\bS$-module $E=\{\K[\bS_n][n-2]\}$ and
whose representations, $\pi:\cA ss_\infty\rar \cE nd_V$, in a dg space $V$ are in one-to-one correspondence with
 Maurer-Cartan elements in the Lie algebra,
$$
 \left(\cG:=\Hom_\bS(E, \cE nd_V)[-1]= \oplus_{n\geq 1} \Hom_\K(V^{\ot n}, V)[1-n], [\ ,\ ]_G\right),
$$
 where $[\ ,\ ]_G$ is the Gerstenhaber bracket.
 \item[(iii)] the particular associative algebra structure $\fs$ on $V$ gives, therefore, rise to
 the associated Maurer-Cartan element $\ga_\fs:= \al_\fs\circ \pi$ in $\cG$; twisting $\cG$  by $\ga_\fs$ gives
 the Hochschild dg Lie algebra, $\cG_\fs=(
 \oplus_{n} \Hom_\K(V^{\ot n}, V)[1-n], [\ ,\ ]_G, d_H:=[\ga_\fs,\ ]_G)$ which indeed controls the
  deformation theory of $\fs$.
\Ei
This is a classical example illustrating how the machine works. For some new
applications
of this approach to deformation theory (e.g.\ to the proof of Deligne's conjecture
or to the deformation theory of associative bialgebras)
 we refer to \cite{KS,MV} and to many references cited there.


\section{Koszul duality theory, quantum BV manifolds and effective $BF$-actions}

\subsection{Quadratic $\fG$-algebras and their Koszul duals}
Koszul duality theory
of quadratic $\fG$-algebras is one of the most powerful theorem-proving
techniques in the theory of (wheeled) operads and properads and their applications.

\sip

What is a {\em quadratic}\, $\fG$-algebra? Every family of graphs, $\fG$, from Table 1 has a uniquely defined
subfamily, $\fG_{gen}$, of {\em generating}\, graphs, which, by definition, is the smallest subset of
$\fG$ with the defining property  that, for every $G\in \fG$ and
any $\fG$-algebra $\cP$, the associated  ``contraction" composition $\mu_\fG: G\langle\cP\rangle\rar \cP$
can be represented as an iteration (in the sense of (\ref{graph-associativity})) of compositions $\mu_{G_i}$
for some $G_i\in \fG_{gen}$, $i\in I$. For example,
$$
\fG_{gen}^\curlywedge=\left\{
\Ba{c}
\begin{xy}
<0mm,0mm>*{\bullet};
<0mm,0mm>*{};<0mm,4mm>*{}**@{-},
<0mm,0mm>*{};<-4mm,-4mm>*{}**@{-},
<0mm,0mm>*{};<-9mm,-4mm>*{}**@{-},
<0mm,0mm>*{};<4mm,-4mm>*{}**@{-},
<0mm,0mm>*{};<9mm,-4mm>*{}**@{-},
<0mm,0mm>*{};<0mm,-4mm>*{}**@{-},
<0mm,-4mm>*{\bullet};
<0mm,-4mm>*{};<-5mm,-8mm>*{}**@{-},
<0mm,-4mm>*{};<-2.7mm,-8mm>*{}**@{-},
<0mm,-4mm>*{};<2.7mm,-8mm>*{}**@{-},
<0mm,-4mm>*{};<5mm,-8mm>*{}**@{-},
<-6mm,-3.7mm>*{...};
<6mm,-3.7mm>*{...};
<0mm,-7.7mm>*{...};
\end{xy}
\Ea
\right\}, \ \
\fG^\uparrow_{gen}=\left\{
\Ba{c}
\begin{xy}
 <0mm,0mm>*{\bullet};
 <0mm,0mm>*{};<-7mm,4mm>*{}**@{-},
 <0mm,0mm>*{};<-4mm,4mm>*{}**@{-},
 <0mm,0mm>*{};<-0.5mm,3.5mm>*{...}**@{},
 <0mm,0mm>*{};<2mm,4mm>*{}**@{-},
 <0mm,0mm>*{};<13mm,5mm>*{}**@{-},
 <0mm,0mm>*{};<-7mm,-4mm>*{}**@{-},
 <0mm,0mm>*{};<-4mm,-4mm>*{}**@{-},
 <0mm,0mm>*{};<0mm,-3.5mm>*{...}**@{},
 <0mm,0mm>*{};<4mm,-4mm>*{}**@{-},
 <0mm,0mm>*{};<7mm,-4mm>*{}**@{-},
 <0mm,0mm>*{};<6mm,4mm>*{}**@{-},
 <6mm,4mm>*{};<13mm,5mm>*{}**@{-},
 <0mm,0mm>*{};<5mm,0mm>*{}**@{-},
<5mm,0mm>*{};<13mm,5mm>*{}**@{-},
 <13mm,5mm>*{\bullet}**@{},
 <13mm,5mm>*{};<6mm,9mm>*{}**@{-},
 <13mm,5mm>*{};<9mm,9mm>*{}**@{-},
 <13mm,5mm>*{};<13mm,8.5mm>*{...}**@{},
 <13mm,5mm>*{};<17mm,9mm>*{}**@{-},
 <13mm,5mm>*{};<20mm,9mm>*{}**@{-},
 <13mm,5mm>*{};<9mm,0mm>*{}**@{-},
 <13mm,5mm>*{};<13mm,0mm>*{...}**@{},
 <13mm,5mm>*{};<16.5mm,0mm>*{}**@{-},
 <13mm,5mm>*{};<20mm,0mm>*{}**@{-},
 \end{xy}
 \Ea\ , \
\Ba{c}
\begin{xy}
 <0mm,0mm>*{\bullet};
 <-0.6mm,0.44mm>*{};<-7mm,4mm>*{}**@{-},
 <-0.4mm,0.7mm>*{};<-3.5mm,4mm>*{}**@{-},
 <0mm,0mm>*{};<0mm,3mm>*{...}**@{},
 <0.4mm,0.7mm>*{};<3.5mm,4mm>*{}**@{-},
 <0.6mm,0.44mm>*{};<7mm,4mm>*{}**@{-},
 <-0.6mm,-0.44mm>*{};<-7mm,-4mm>*{}**@{-},
 <-0.4mm,-0.7mm>*{};<-3.5mm,-4mm>*{}**@{-},
 <0mm,0mm>*{};<0mm,-3mm>*{...}**@{},
 <0.4mm,-0.7mm>*{};<3.5mm,-4mm>*{}**@{-},
 <0.6mm,-0.44mm>*{};<7mm,-4mm>*{}**@{-},
\end{xy}
\
\begin{xy}
 <0mm,0mm>*{\bullet};
 <-0.6mm,0.44mm>*{};<-7mm,4mm>*{}**@{-},
 <-0.4mm,0.7mm>*{};<-3.5mm,4mm>*{}**@{-},
 <0mm,0mm>*{};<0mm,3mm>*{...}**@{},
 <0.4mm,0.7mm>*{};<3.5mm,4mm>*{}**@{-},
 <0.6mm,0.44mm>*{};<7mm,4mm>*{}**@{-},
 <-0.6mm,-0.44mm>*{};<-7mm,-4mm>*{}**@{-},
 <-0.4mm,-0.7mm>*{};<-3.5mm,-4mm>*{}**@{-},
 <0mm,0mm>*{};<0mm,-3mm>*{...}**@{},
 <0.4mm,-0.7mm>*{};<3.5mm,-4mm>*{}**@{-},
 <0.6mm,-0.44mm>*{};<7mm,-4mm>*{}**@{-},
\end{xy}
 \Ea
 \right\},
$$
and\vspace{-2mm}
$$
\fG^\circlearrowright_{c, gen}=\left\{
\Ba{c}
 \begin{xy}
 <0mm,0mm>*{\bullet};
 <-0.6mm,0.44mm>*{};<-8mm,5mm>*{}**@{-},
 <-0.4mm,0.7mm>*{};<-4.5mm,5mm>*{}**@{-},
 <0mm,0mm>*{};<0mm,5mm>*{\ldots}**@{},
 <0.4mm,0.7mm>*{};<4.5mm,5mm>*{}**@{-},
 <0.6mm,0.44mm>*{};<12.4mm,4.8mm>*{}**@{-},
 <-0.6mm,-0.44mm>*{};<-8mm,-5mm>*{}**@{-},
 <-0.4mm,-0.7mm>*{};<-4.5mm,-5mm>*{}**@{-},
 <0mm,0mm>*{};<-1mm,-5mm>*{\ldots}**@{},
 <0.4mm,-0.7mm>*{};<4.5mm,-5mm>*{}**@{-},
 <0.6mm,-0.44mm>*{};<8mm,-5mm>*{}**@{-},
 <13mm,5mm>*{};<13mm,5mm>*{\bullet}**@{},
 <12.6mm,5.44mm>*{};<5mm,10mm>*{}**@{-},
 <12.6mm,5.7mm>*{};<8.5mm,10mm>*{}**@{-},
 <13mm,5mm>*{};<13mm,10mm>*{\ldots}**@{},
 <13.4mm,5.7mm>*{};<16.5mm,10mm>*{}**@{-},
 <13.6mm,5.44mm>*{};<20mm,10mm>*{}**@{-},
 <12.4mm,4.3mm>*{};<8mm,0mm>*{}**@{-},
 <12.6mm,4.3mm>*{};<12mm,0mm>*{\ldots}**@{},
 <13.4mm,4.5mm>*{};<16.5mm,0mm>*{}**@{-},
 <13.6mm,4.8mm>*{};<20mm,0mm>*{}**@{-},
 \end{xy}\Ea \ , \
\begin{xy}
 <0mm,0mm>*{\bullet};
 <-0.6mm,0.44mm>*{};<-8mm,5mm>*{}**@{-},
 <-0.4mm,0.7mm>*{};<-4.5mm,5mm>*{}**@{-},
 <0mm,0mm>*{};<0mm,4mm>*{\ldots}**@{},
 <0.4mm,0.7mm>*{};<4.5mm,5mm>*{}**@{-},
 <0.6mm,0.44mm>*{};<6mm,4mm>*{}**@{-},
 <-0.6mm,-0.44mm>*{};<-8mm,-5mm>*{}**@{-},
 <-0.4mm,-0.7mm>*{};<-4.5mm,-5mm>*{}**@{-},
 <0mm,0mm>*{};<0mm,-4mm>*{\ldots}**@{},
 <0.4mm,-0.7mm>*{};<4.5mm,-5mm>*{}**@{-},
 <0.6mm,-0.44mm>*{};<6mm,-4mm>*{}**@{-},
(6,4)*{}
   \ar@{->}@(ur,dr) (6,-4)*{}
 \end{xy}
\right\}.
$$

\subsubsection{\bf Weight gradation}
let $\fG$ be a family of graphs from Table 1.
For any genus $q$ graph $G\in \fG$ with $p$ vertices we set $||G||:=p+q$ if the family $\fG$ contains wheels and set $||G||:=p$ otherwise. This number is called the {\em weight}\, of $G$. Thus $\fG_{gen}\subset \fG$ consists precisely of graphs
of weight $2$.

 \sip

 For an $\bS$-bimodule $E$ let
$\cF_{(\lambda)}^\fG\langle E \rangle$
 stand for a subspace of the free $\fG$-algebra $\cF^\fG\langle E \rangle$
spanned by decorated graphs of weight $\la$. Operadic compositions in  $\cF^\fG\langle E \rangle$ are homogeneous with respect to the weight gradation.

\subsubsection{\bf Definition}
A $\fG$-algebra, $\cP$,  is called {\em quadratic}\, if it is the
  quotient, $\cF^\fG\langle E\rangle / \langle\cR \rangle$, of a free $\fG$-algebra
 (generated by an $\bS$-bimodule $E$) modulo the ideal generated
by a subspace
$\cR\subset \cF_{(2)}^\fG\langle E \rangle=\bigoplus_{G\in \fG_{gen}} G\langle E \rangle$.
It comes equipped with an induced weight gradation,
$\cP=\oplus_{\la\geq 1}\cP_{(\la)}$, where $\cP_{(\la)}=
\cF_{(\lambda)}^\fG\langle E \rangle/\langle\cR\rangle$.
In particular, $\cP_{(1)}=E$ and $\cP_{(2)}=\cF_{(2)}^\fG\langle E \rangle/\cR$.

\subsection{Koszul duality}
Let $\fG_c$ be any family of {\em connected}\, graphs from Table 1. In this case one can associate
to any quadratic $\fG_c$-algebra $\cP$ its Koszul dual $\fG_c$-coalgebra $\cP^{\mbox {\scriptsize !`}}$.
We omit technical details (referring to  \cite{GK,GJ,G,V,MMS,Me-BV}) and explain just the working scheme:
\sip

(i) The notion of $\fG_c$-{\em coproperad}\, is obtained by an obvious dualization of the notion of $\fG_c$-algebra (see
\S 2.3): this is an $\bS$-bimodule $\cP=\{\cP(m,n)\}$ together with
a family of linear $\bS_m\times \bS_n$-equivariant maps,
$$
\left\{\Delta_G: \cP(m,n) \rar G\langle \cP\rangle \right\}_{G\in \fG_c(m,n), m,n\geq 0},
$$
which
satisfy the  coassociativity condition, $
\Delta_G=\Delta_H' \circ \Delta_{G/H}$,
 for any
 subgraph $H\subset G$ which belongs to the family $\fG$. Here $\Delta_H':  (G/H)\langle E\rangle \rar G\langle E \rangle$
 is the map
which equals $\Delta_H$ on the distinguished vertex of $G/H$  and which is the identity on
all other vertices of $G$.

\sip
(ii) There exists a pair of adjoint {\em exact}\, functors
$$
\Ba{rcccl}
B: &
 \mbox{the category of dg $\fG_c$-algebras} & \rightleftarrows &
 \mbox{the category of dg $\fG_c$-coalgebras}: & B^c\\
 &\cP & \lon & (B(\cP), \p_\cP) &\\
 &(B^c(\cQ),\p_\cQ) &\longleftarrow & \cQ &
 \Ea
$$
such that, for any dg $\fG_c$-algebra $\cP$ the composition $B^c(B(\cP))$ is a dg free
resolution of $\cP$. The differential $\p_\cP$ in $B(\cP)$ encodes both the differential
and all the generating contraction compositions,
$\{\mu_\fG: G\langle\cP \rangle\rar \cP\}_{G\in \fG_{gen}}$, in the $\fG_c$-algebra $\cP$
(and similarly for $\p_\cQ$).

\sip

(iii) As a vector space $B(\cP)$ is isomorphic to the free $\fG_c$-algebra, $\cF^{\fG_c}\langle
\hat{\cP}\rangle$, generated by an $\bS$-bimodule $\hat{\cP}$ which is linearly isomorphic to $\cP$
and hence comes equipped with an induced weight gradation. The subspace $\cB(\cP_{(1)}):=
\cF^{\fG_c}\langle \hat{\cP}_{(1)}\rangle$ of
$\cB(\cP)$ is obviously a sub-coproperad. On the other hand, $B(\cP)$ has its own ``outer" weight gradation, $B(\cP)=\oplus_{\mu\geq 1} B_{(\mu)}(\cP)$, induced from the weight gradation of the free algebra,
$B_{(\mu)}(\cP):=\cF^{\fG_c}_{(\mu)}\langle
\hat{\cP}\rangle$; the cobar differential $\p_{\cP}$ has weight $-1$ with respect to this outer weight gradation.

\subsubsection{\bf Definition}
Given a quadratic $\fG_c$-algebra $\cP$, the
$\fG_c$-coalgebra $\cP^{\mbox {\scriptsize !`}}=\oplus_{\mu\geq 1}\cP^{\mbox {\scriptsize !`}}_{(\mu)}$
with $\cP^{\mbox {\scriptsize !`}}_{(\mu)}:= B_{(\mu)}(\cP_{(1)}) \cap \Ker \p_{\cP}\subset B(\cP)$
is called {\em Koszul dual}\, to $\cP$.

\sip

The beauty of this notion is that $\cP^{\mbox {\scriptsize !`}}$ is again quadratic and, moreover,
 can
often be
easily  computed directly from generators and relations, $E$ and $\cR$, of $\cP$.

\subsubsection{\bf Definition }
A quadratic $\fG_c$-algebra $\cP$ is called {\em Koszul}, if the associated
inclusion of dg coproperads, $\imath: (\cP^{\mbox {\scriptsize !`}},0)  \lon (B(\cP), \p_\cP)$, is
 a quasi-isomorphism.

\mip

As the cobar construction functor $B^c$  preserves quasi-isomorphisms between
connected $\fG$-coalgebras,
 the composition
$$
\pi: \cP_\infty:=B^c(\cP^{\mbox {\scriptsize !`}}) \stackrel{B^c(\imath)}{\lon} B^c(B(\cP))
\stackrel{\mbox{\small natural projection}}{\lon} \cP
$$
is a quasi-isomorphism if and only if $\cP$ is Koszul;
in this case the dg free $\fG_c$-algebra $\cP_\infty$ gives us
 a minimal resolution of the quadratic algebra $\cP$.
 Almost all minimal resolutions listed in \S 2 have been obtained
 in this way.

 \subsection{Homotopy transfer formulae}\label{5: homotopy transfer}
  If $\cP_\infty$ is a minimal resolution of some
 $\fG$-algebra $\cP$, and $(V,d)$ is a complex carrying a $\cP$-structure, then one might expect
that the associated cohomology space, $H(V,d)$, carries an induced structure of
$\cP_\infty$-algebra. In the case when $\cP$ is an operad of associative algebras, existence
of such induced $\cA ss_\infty$-structures was proven by Kadeishvili in \cite{Ka}
and the first explicit formulae have been shown in \cite{Me0}. Later Kontsevich and Soibelman
\cite{KS} have nicely rewritten these homotopy transfer formulae in terms of sums of decorated graphs.
In fact, it is a general phenomenon that the homotopy transfer formulae
can be represented as sums of graphs. The required graphs are precisely the ones
which describe the image of the natural inclusion
 $\imath: (\cP^{\mbox {\scriptsize !`}},0)  \lon (B(\cP), \p_\cP)$, and apply to any quadratic
 $\fG$-algebra, not necessarily the Koszul one \cite{Me-BV}.

\subsection{Example: unimodular Lie 1-bialgebras versus quantum BV manifolds}
The wheeled prop, $\cU\LB$, of {\em unimodular}\, Lie 1-bialgebras was defined in \cite{Me-BV}
(cf.\ \S\ref{2: subsect LB and Poisson}) as the quotient,
$\cF^\circlearrowright_c\langle B\rangle/\langle \cR    \rangle$, of the free wheeled
properad generated by $\bS$-bimodule (\ref{2: generators LB}) modulo the ideal generated by
relations (\ref{2: jacobi equations}), (\ref{2: Lie 1bi relations}) and the following ones,
$
\Ba{c}
\begin{xy}
 <0mm,-0.55mm>*{};<0mm,-2.5mm>*{}**@{-},
 <0.5mm,0.5mm>*{};<2.2mm,2.2mm>*{}**@{-},
 <-0.48mm,0.48mm>*{};<-2.9mm,3.2mm>*{}**@{-},
 <0mm,0mm>*{\bullet};<0mm,0mm>*{}**@{},
(2.2,2.2)*{}
   \ar@{->}@(ur,d) (0,-1)*{}
 \end{xy}=0,$
 $\begin{xy}
 <0mm,0.66mm>*{};<0mm,3mm>*{}**@{-},
 <0.39mm,-0.39mm>*{};<2.2mm,-2.2mm>*{}**@{-},
 <-0.35mm,-0.35mm>*{};<-2.8mm,-3.2mm>*{}**@{-},
 <0mm,0mm>*{\bullet};<0mm,0mm>*{}**@{},
(0.0,1.0)*{}
   \ar@{->}@(u,dr) (2.2,-2.2)*{}
\end{xy}=0,
\Ea
$
expressing unimodularity of both  binary operations. This is a quadratic wheeled properad
so  that one can apply the above general machinery to compute its Koszul dual coproperad,
$\cU\LB^{\mbox {\scriptsize !`}}$, and then the dg  properad
$\cP_\infty:=B^c(\cU\LB^{\mbox {\scriptsize !`}})$ which turns out to be a free
wheeled properad, $\cF^\circlearrowright_c\langle Z\rangle$, generated by an $\bS$-bimodule,
\vspace{-2mm}
$$
Z(m,n):= \bigoplus_{a\geq 0}^\infty
\sgn_m\ot \id_n[m-2-2a]=\mbox{span}\left\langle
\begin{xy}
 <0mm,0mm>*{\mbox{$\xy *=<3mm,3mm>
\txt{{{$_a$}}}*\frm{-}\endxy$}};
<-1.5mm,1.5mm>*{};<-5mm,5mm>*{}**@{-},
<-1mm,1.5mm>*{};<-2.7mm,5mm>*{}**@{-},
<1.5mm,1.5mm>*{};<5mm,5mm>*{}**@{-},
<1mm,1.5mm>*{};<2.7mm,5mm>*{}**@{-},
<0mm,4.5mm>*{...};
<-1.5mm,-1.5mm>*{};<-5mm,-5mm>*{}**@{-},
<-1mm,-1.5mm>*{};<-2.7mm,-5mm>*{}**@{-},
<1.5mm,-1.5mm>*{};<5mm,-5mm>*{}**@{-},
<1mm,-1.5mm>*{};<2.7mm,-5mm>*{}**@{-},
<0mm,-4.5mm>*{...};
<0mm,0mm>*{};<-5.8mm,6.3mm>*{^1}**@{},
   <0mm,0mm>*{};<-3mm,6.3mm>*{^2}**@{},
   <0mm,0mm>*{};<2.6mm,6.3mm>*{^{m\hspace{-0.5mm}-\hspace{-0.5mm}1}}**@{},
   <0mm,0mm>*{};<7.0mm,6.3mm>*{^m}**@{},
<0mm,0mm>*{};<-5.8mm,-6.3mm>*{_1}**@{},
   <0mm,0mm>*{};<-3mm,-6.3mm>*{_2}**@{},
   <0mm,0mm>*{};<3mm,-6.3mm>*{_{n\hspace{-0.5mm}-\hspace{-0.5mm}1}}**@{},
   <0mm,0mm>*{};<7.0mm,-6.3mm>*{_n}**@{},
 \end{xy}
\right\rangle_{m+n+2a\geq 3\atop m+a\geq 1, n+a\geq 1}.
$$
and equipped with the following differential\vspace{-2mm}
$$\label{3:Differential in ULieB_infty}
\delta
\begin{xy}
 <0mm,0mm>*{\mbox{$\xy *=<3mm,3mm>
\txt{{{$_a$}}}*\frm{-}\endxy$}};
<-1.5mm,1.5mm>*{};<-5mm,5mm>*{}**@{-},
<-1mm,1.5mm>*{};<-2.7mm,5mm>*{}**@{-},
<1.5mm,1.5mm>*{};<5mm,5mm>*{}**@{-},
<1mm,1.5mm>*{};<2.7mm,5mm>*{}**@{-},
<0mm,4.5mm>*{...};
<-1.5mm,-1.5mm>*{};<-5mm,-5mm>*{}**@{-},
<-1mm,-1.5mm>*{};<-2.7mm,-5mm>*{}**@{-},
<1.5mm,-1.5mm>*{};<5mm,-5mm>*{}**@{-},
<1mm,-1.5mm>*{};<2.7mm,-5mm>*{}**@{-},
<0mm,-4.5mm>*{...};
<0mm,0mm>*{};<-5.8mm,6.3mm>*{^1}**@{},
   <0mm,0mm>*{};<-3mm,6.3mm>*{^2}**@{},
   <0mm,0mm>*{};<2.6mm,6.3mm>*{^{m\hspace{-0.5mm}-\hspace{-0.5mm}1}}**@{},
   <0mm,0mm>*{};<7.0mm,6.3mm>*{^m}**@{},
<0mm,0mm>*{};<-5.8mm,-6.3mm>*{_1}**@{},
   <0mm,0mm>*{};<-3mm,-6.3mm>*{_2}**@{},
   <0mm,0mm>*{};<3mm,-6.3mm>*{_{n\hspace{-0.5mm}-\hspace{-0.5mm}1}}**@{},
   <0mm,0mm>*{};<7.0mm,-6.3mm>*{_n}**@{},
 \end{xy}
 =(-1)^{m-1}
\begin{xy}
 <0mm,0mm>*{\mbox{$\xy *=<6mm,3mm>
\txt{{{$_{a-1}$}}}*\frm{-}\endxy$}};
<-1.5mm,1.5mm>*{};<-5mm,5mm>*{}**@{-},
<-1mm,1.5mm>*{};<-2.7mm,5mm>*{}**@{-},
<1.5mm,1.5mm>*{};<4mm,4mm>*{}**@{-},
<1mm,1.5mm>*{};<2.7mm,5mm>*{}**@{-},
<0mm,4.5mm>*{...};
<-1.5mm,-1.5mm>*{};<-5mm,-5mm>*{}**@{-},
<-1mm,-1.5mm>*{};<-2.7mm,-5mm>*{}**@{-},
<1.5mm,-1.5mm>*{};<4mm,-4mm>*{}**@{-},
<1mm,-1.5mm>*{};<2.7mm,-5mm>*{}**@{-},
<0mm,-4.5mm>*{...};
<0mm,0mm>*{};<-5.8mm,6.3mm>*{^1}**@{},
   <0mm,0mm>*{};<-3mm,6.3mm>*{^2}**@{},
   <0mm,0mm>*{};<3.0mm,6.3mm>*{^m}**@{},
<0mm,0mm>*{};<-5.8mm,-6.3mm>*{_1}**@{},
   <0mm,0mm>*{};<-3mm,-6.3mm>*{_2}**@{},
   <0mm,0mm>*{};<3.0mm,-6.3mm>*{_n}**@{},
   \ar@{->}@(ur,dr) (4.0,4.0)*{};(4.0,-4.0)*{},
 \end{xy}
 \ + \sum_{a=b+c\atop b,c\geq 0}\sum_{m=I'\sqcup I''\atop
 [n]=J'\sqcup J''}
 (-1)^{\sigma(I_1\sqcup I_2) + |I_1|(|I_2|+1)}
 \Ba{c}
 \begin{xy}
 <0mm,0mm>*{\mbox{$\xy *=<3mm,3mm>
\txt{{{$_b$}}}*\frm{-}\endxy$}};
<-1.5mm,1.5mm>*{};<-5mm,5mm>*{}**@{-},
<-1mm,1.5mm>*{};<-2.7mm,5mm>*{}**@{-},
<1.5mm,1.5mm>*{};<8mm,8mm>*{}**@{-},
<1mm,1.5mm>*{};<2.7mm,5mm>*{}**@{-},
<0mm,4.5mm>*{...};
<-1.5mm,-1.5mm>*{};<-5mm,-5mm>*{}**@{-},
<-1mm,-1.5mm>*{};<-2.7mm,-5mm>*{}**@{-},
<1.5mm,-1.5mm>*{};<5mm,-5mm>*{}**@{-},
<1mm,-1.5mm>*{};<2.7mm,-5mm>*{}**@{-},
<0mm,-4.5mm>*{...};
 <9.5mm,9.5mm>*{\mbox{$\xy *=<3mm,3mm>
\txt{{{$_c$}}}*\frm{-}\endxy$}};
<8mm,11mm>*{};<5mm,14mm>*{}**@{-},
<8.5mm,11mm>*{};<7mm,14mm>*{}**@{-},
<10.5mm,11mm>*{};<12mm,14mm>*{}**@{-},
<11mm,11mm>*{};<14mm,14mm>*{}**@{-},
<9.5mm,13.5mm>*{...};
<8.5mm,8mm>*{};<7mm,4.5mm>*{}**@{-},
<10.5mm,8mm>*{};<12mm,4.5mm>*{}**@{-},
<11mm,8mm>*{};<14mm,4.5mm>*{}**@{-},
<9.5mm,5mm>*{...};
<10.5mm,1.8mm>*{\underbrace{\ \ \ \ \ \ \  }_{J''}};
<0mm,-8mm>*{\underbrace{\ \ \ \ \ \ \ \ \  }_{J'}};
<-1mm,8.2mm>*{\overbrace{\ \ \ \ \ \ \  }^{I'}};
<9.5mm,17mm>*{\overbrace{\ \ \ \ \ \ \  \ }^{I''}};
 \end{xy}
 \Ea
$$
It is not known at present whether or not $\cU\LB$ is Koszul, i.e.\ whether or not
the above free properad is a (minimal) resolution of the latter. In any case, $\cU\LB_\infty$
gives us an approximation to that minimal resolution, and has, in fact, a geometrically
meaningful set, $\{ \cU\LB_\infty\rar \cE nd_V\}$, of all possible representations. To
describe this set let us recall a few notions from the Schwarz model \cite{Sc} of the Batalin-Vilkovisky
quantization formalism \cite{BV}.

\subsection{Formal quantum BV manifolds} Let
$\{x^a, \psi_a, \hbar\}_{1\leq a \leq n}$, $n\in \N$, be a set of formal
homogeneous variables of degrees
$|x^a| + |\psi_a|=1$ and $|\hbar|=2$, and let $\f_{x,\psi}^\hbar:=\K[[x^a,\psi_a,\hbar]]$ be
the associated
free graded commutative ring which we view from now on as a $\K[[\hbar]]$-algebra.
The degree $-1$ Lie bracket,
$$
\{f\bullet g\}:= (-1)^{|f|}\Delta(fg) - (-1)^{|f|}\Delta(f)g - f\Delta(g),\ \forall \ f,g\in
\f^\hbar_{x,\psi}
$$
make $\f_{x,\psi}^\hbar$
 into a Gerstenhaber $\K[[\hbar]]$-algebra (see \S\ref{3: Gerst algebras and HM}). Here and elsewhere
$\Delta:=\sum_{a=1}^n(-1)^{|x^a|}\frac{\p^2}{\p x^a\p \psi_a}$. A {\em quantum master function}\,
is, by definition, a degree 2 element $\Gamma\in \f_{x,\psi}^\hbar$ satisfying a
so called {\em quantum master equation}
\Beq\label{5: QM equation}
\hbar\Delta \Gamma + \frac{1}{2}\{\Gamma\bullet \Gamma\}=0.
\Eeq
Such an element makes the  $\K[[\hbar]]$-module  $\f_{x,\psi}^\hbar$
{\em differential}\, with the differential $\Delta_\Ga:= \hbar \Delta + \{\Ga\bullet\ \ \}$. Note that this differential does {\em not}\, respect the algebra structure in  $\f_{x,\psi}^\hbar$ but respects the Poisson brackets.

Consider a group of $\K[[\hbar]]$-algebra automorphisms,
$F: \f_{x,\psi}^\hbar\rar \f_{x,\psi}^\hbar$, preserving the Lie brackets, $F(\{f\bullet g\})=
\{F(f)\bullet F(g)\}$ (but not necessarily the operator $\Delta$); this group is uniquely determined
 by a collection, $\cN:=\{|x^a|, |\psi_a|\}_{1\leq a\leq n}$, of
$2n$ integers and is denoted by
$Symp_\cN$. It is often called a group of {\em symplectomorphsims} of the Gerstenhaber
algebra $(\f_{x,\psi}^\hbar, \{\ \bullet\ \})$.
A remarkable fact \cite{Kh} is that $Symp_\cN$  acts
on the  set
of quantum master functions by the formula,
\Beq\label{Semidensity_Gamma_trasnformation}
e^{\frac{F(\Gamma)}{\hbar}}:
=
\left[\mbox{Ber}\left(\Ba{cc}
\frac{\p F(x^a)}{\p x^b} & \frac{\p F(x^a)}{\p \psi_b}\\
\frac{\p F(\psi_a)}{\p x^b} & \frac{\p F(\psi_a)}{\p \psi_b}
\Ea
\right)
\right]^{-\frac{1}{2}}
e^{\frac{\Gamma({x}, {\psi}, \hbar)}{\hbar}}.
\Eeq

\vspace{0mm}

\subsubsection{\bf Definition} An equivalence class of
 pairs, $(\f_{x,\psi}^\hbar, \Ga)$, under the action of the group
 $Symp_\cN$ is called a formal {\em quantum BV manifold}\, $\cM$ of dimension $\cN$. A particular
 representative, $(\f_{x,\psi}^\hbar, \Ga)$, of $\cM$ is called
 a {\em Darboux coordinate chart}\, on $\cM$.

\sip

In geometric terms, $\cM$ is a formal odd symplectic manifold equipped with a special type
semidensity \cite{Kh,Sc}. We need an extra structure on $\cM$ which we again define with the
help of
a Darboux coordinate chart. Notice that the ideals, $I_x$ and $I_\psi$, in the
$\K[[\hbar]]$-algebra $\f_{x,\psi}^\hbar$ generated, respectively, by $\{x^a\}_{1\leq a\leq n}$
and $\{\psi_a\}_{1\leq a\leq n}$, are also Lie ideals; geometrically, they define a pair
of transversally intersecting Lagrangian submanifolds of $\cM$. A quantum BV manifold
$\cM$ is said to have {\em split quasi-classical limit}\, (or, slightly shorter,
 $\cM$ is  {\em quasi-classically split})
 if it admits a Darboux coordinate chart
in which the master function, $\Ga(x,\psi,\hbar)=\sum_{n\geq 0} \Ga_n(x,\psi)\hbar^n$, satisfies
the following two boundary conditions,
$$
\Ga_0\in I_xI_y, \ \ \  \Ga_1\in I_x+ I_y.
$$
In plain terms, these conditions mean that $\Ga(x,\psi,\hbar)$ is given by a formal power series
of the form,
\Beq\label{5: Ga_coordinate_decomposition}
\Gamma(x, \psi, \hbar)=
\underbrace{\sum_{a,b}\Gamma_{(0)\, b}^{\ \ \ a}x^b \psi_a}_{\Gamma_0}
+ \underbrace{
\sum_{p+q+2n\geq 3\atop
p+n\geq 1, q+n\geq 1}\frac{1}{p!q!}
\Gamma_{(n)\, a_1\ldots a_p}^{\ \ \ b_1\ldots b_q} x^{a_1}\ldots  x^{a_p}\psi_{b_1}\ldots
\psi_{b_q} \hbar^n}_{\bGa}
\Eeq
for some  $\Gamma_{(n)\, a_1\ldots a_p}^{\ \ \ b_1\ldots b_q}\in \K$. Quantum master equation
(\ref{5: QM equation}) immediately implies that $\{\Ga_0, \Ga_0\}=0$ so that $\eth:=\{\Gamma_0\bullet\ \
\}$ is a differential in the Gerstenhaber algebra $\f_{x,\psi}^\hbar$. Then master equation
(\ref{5: QM equation})
for a quasi-classically split quantum master function  can be equivalently rewritten in the form,
$$
\eth \bGa + \hbar \Delta \bGa + \frac{1}{2}\{\bGa\bullet \bGa\}=0,
$$
where $\bGa$ is an element of $\f_{x,\psi}^\hbar$ of polynomial order at least $3$
(here we set, by definition, the polynomial order of generators $x$ and $\psi$  equal to $1$
and the polynomial order of $\hbar$  equal to 2). The differential $\eth$ induces a differential
on the tangent space, $\cT_*\cM$, to $\cM$ at the distinguished point; we denote it by the same letter
$\eth$. Such a quantum BV manifold is called {\em minimal}\, if $\eth=0$ and {\em contractible}\, if there exists
a Darboux coordinate chart in which $\bGa=0$ (i.e.\ $\Ga=\eth$) and the tangent complex $(\cT_*\cM,\eth)$ is acyclic
(cf.\ \S \ref{2: subsub on Lie_infty}).

An important class of so called $BF$ field theories (see, e.g., \cite{CR} and references cited there)
have associated quantum BV manifolds which do satisfy
the split quasi-classical limit condition.

\subsubsection{\bf Proposition}\label{5: Prop on repr of ULB and quantum BV}
{\em For any dg vector space $V$, there is a one-to-one correspondence
between representations, $\cU\LB_\infty\rar \cE nd_V$,  and structures of formal
quasi-classically split quantum BV manifold
 on $\cM_{V\oplus V^*[1]}$, the formal manifold associated to
$V\oplus V^*[1]$.}

\subsection{Morphisms of quantum BV manifolds \cite{Me-BV}} The above Proposition together with
the Koszul duality theory approach to the homotopy transfer outlined in \S\ref{5: homotopy transfer}
provide us with highly non-trivial formulae for constructing quantum BV manifold structures
out of dg unimodular Lie 1-bialgebras. We would like to have
a {\em category}\, of quantum BV manifolds in which such homotopy transfer formulae
can be interpreted as morphisms. This can be achieved via the following

\subsubsection{\bf Definitions} (i) A {\em morphism}\, of quasi-classically split quantum BV manifolds, $F:\cM\rar \hat{\cM}$,
is, by definition, a morphism of dg $\K[[\hbar]]$-modules,
$$
F: \left(\f_{\hat{\cM}}\simeq \f_{\hat{x},\hat{\psi}}^\hbar, \Delta_{\hat{\Ga}}\right)\lon
 \left(\f_{\cM}\simeq\f_{x,\psi}^\hbar, \Delta_\Ga\right),
$$
inducing in the classical limit $\hbar\rar 0$ a morphism of algebras, $F|_{\hbar=0}:
\f_{\hat{x},\hat{\psi}}\rar \f_{{x},{\psi}}$ which preserves the ideals, $F|_{\hbar=0}(\langle\hat{x}\rangle)\subset \langle x \rangle$ and $F|_{\hbar=0}(\langle \hat{\psi} \rangle)\subset \langle\psi \rangle$,
 of the distinguished Lagrangian submanifolds in $\hat{\cM}|_{\hbar=0}$ and ${\cM}|_{\hbar=0}$.
\sip

 (ii) If $F:\cM\rar \hat{\cM}$ is a morphism of quantum BV manifolds, then $dF|_{\hbar=0}$ induces in fact
 a morphism of dg vector spaces, $(\cT_*\cM, \eth) \rar (\cT_*\hat{\cM}, \hat{\eth})$;
 $F$ is called a {\em quasi-isomorphism}\, if the latter map induces an isomorphism of the associated
 cohomology groups.

\subsubsection{\bf Theorem \cite{Me-BV}} {\em Every quantum quasi-classically split  BV manifold is isomorphic to
the product of a minimal one and of a contractible one. In particular, every such  manifold is quasi-isomorphic to a minimal one}.

\subsection{Homotopy transfer of quantum BV-structures via Feynman integral} Homotopy transfer formulae
of $\cP_\infty$-structures given by Koszul duality theory are given by  sums
of decorated graphs which resemble Feynman diagrams in quantum field theory. This resemblance
was made a rigorous fact in \cite{Mn} for the case of the wheeled operad of unimodular Lie algebras (see
\S\ref{3: subsect unimod Lie}).
\sip

Given any complex $V$ and a dg Lie 1-bialgebra
structure on $V$ with degree $0$ Lie cobrackets $\Delta^{\CoLie}: V\rar \wedge^2 V$
and degree $1$ Lie brackets  $[\ \bullet\ ]: \odot^2 V\rar V[1]$, the associated by
Koszul duality theory
homotopy  formulae transfer this rather trivial quantum BV manifold structure on $V$ to a highly
non-trivial quantum master function on its cohomology $H(V)$; the same formulae
can also be described \cite{Me-BV}
by a standard Batalin-Vilkovisky quantization \cite{BV} of a $BF$-type field theory on
the space $V\oplus V^*[1]$ with the action given by,
$$
\Ba{rccc}
S:& V\oplus V^*[1]  & \lon & \K\vspace{2mm}\\
& p\oplus \om & \lon & S(p,\om):= \langle p, d\om \rangle
 \ + \ \frac{1}{2}\langle p, [\om , \om]\rangle +
\frac{1}{2}\langle [p \bullet p], \om\rangle,
\Ea
$$
where $\langle\ ,\ \rangle$ stands for the natural pairing, and
$[\ ,\ ]: \odot^2(V^*[1])\rar V^*[2]$ for the dualization of $\Delta^{\CoLie}$.
Thus at least in some cases the Koszul duality technique for homotopy transfer of $\infty$-structures
is identical to the Feynman diagram technique in theoretical physics. The beauty of the latter
lies in its combinatorial simplicity (due to the Wick theorem), while the power of the former
lies in its generality: the Koszul duality theory applies to any (non-commutative case including)
quadratic $\fG$-algebras.


\begin{thebibliography}{10}

\bibitem[BaVi]{BV} I.\ Batalin and G.\ Vilkovisky, {\em Gauge algebra and quantization},
Phys.\ Lett.\ B {\bf 102} (1981), 27.

\bibitem[CaRo]{CR} A.\ Cattaneo and C.\ Rossi, {\em Higher-dimensional $BF$
theories in the Batalin-Vilkovisky formalism: the BV action and generalized Wilson loops},
Commun.\ Math.\ Phys. {\bf 221} (2001), 591-657.


\bibitem[ChLi]{CL}  F.\ Chapoton and  M.\ Livernet,
 {\em  Pre-Lie algebras and the rooted trees operad},
Internat. Math. Res. Notices, {\bf 8} (2001), 395–408.





\bibitem[DoKh]{DK}
V.\ Dotsenko and A.\ Khoroshkin,
{\em Character formulas for the operad of two compatible brackets
and for the bihamiltonian operad}. Funktsional. Anal. i Prilozhen., {\bf 41}(1) (2007) 1–17.


\bibitem[Dr]{Dr}
V.\ Drinfeld,
{\em On some unsolved problems in quantum group theory}.
In: Lecture Notes in Math., Springer,  {\bf 1510} (1992), 1-8.



\bibitem[Ga]{G} W.L.\ Gan, {\em Koszul duality for dioperads}, Math.\ Res.\ Lett.
{\bf 10} (2003), 109-124.

\bibitem[GeJo]{GJ} E.\ Getzler and J.D.S.\ Jones, {\em Operads, homotopy algebra, and
 iterated integrals for double loop spaces},
 preprint  hep-th/9403055.


\bibitem[GiKa]{GK} V.\ Ginzburg and M.\ Kapranov, {\em Koszul duality for operads},
{Duke Math. J.} {\bf 76} (1994) {203--272}.

\bibitem[Gr1]{Gr1} J.\ Gran\aa ker, {\em Unimodular $L_\infty$-algebras},
preprint arXiv:0803.1763 (2008).

\bibitem[Gr2]{Gr2} J.\ Gran\aa ker, {\em Quantum BV manifolds and Lie quasi-bialgebras},
to appear in J.\ Diff.\ Geom. and Applications.

\bibitem[HeMa]{HM} C.\ Hertling and Yu.I.\ Manin, {\em Weak Frobenius manifolds},
 Intern.\ Math.\ Res.\ Notices {\bf 6} (1999) 277-286.

\bibitem[HMT]{HM2}
  C.\ Hertling, Yu.\ Manin and C.\ Teleman, {\em
An update on semisimple quantum cohomology and F-manifolds}, preprint arXiv:0803.2769 (2008).

\bibitem[Ho]{Ho}
E.\  Hoffbeck, {\em A Poincar´e-Birkhoff-Witt criterion for Koszul operads}. Prerint arXiv:0709.2286, 2008.

 \bibitem[Ka]{Ka} T.V.\ Kadeishvili,   {\em The algebraic structure in the cohomology of
 $A(\infty)$-algebras}, Soobshch.\ Akad.\ Nauk Gruzin.\ SSR {\bf 108} (1982),  249-252.

\bibitem[Kh]{Kh} H.\ Khudaverdian, {\em Semidensities on odd symplectic supermanifolds},
 Commun.\ Math.\ Phys. {\bf 247} (2004), 353–390.


\bibitem[Ko1]{Ko0} M.\ Kontsevich, unpublished (but see \cite{MaVo}).


 \bibitem[Ko2]{Ko} M.\ Kontsevich, {\em Deformation quantization
 of Poisson manifolds}, Lett.\ Math.\ Phys. {\bf 66} (2003), 157-216.

\bibitem[KoSo]{KS} M.\ Kontsevich and Y.\ Soibelman, {\em
Deformations of algebras over operads and the {D}eligne
              conjecture}, {Conf\'erence Mosh\'e Flato 1999, Vol. I (Dijon)}, pp. {255--307},
              Kluwer Acad. Publ. (Dordrech),  2000.


\bibitem[MMS]{MMS}  M.\ Markl, S.\ Merkulov and S.\ Shadrin, {\em Wheeled props and the master
equation}, preprint math.AG/0610683, J.\ Pure and Appl.\ Algebra {\bf 213} (2009), 496-535.


 \bibitem[MaVo]{MaVo} M.\ Markl and A.A.\ Voronov,
{\em PROPped up graph cohomology}.
In: ``Algebra, Arithmetic
and Geometry - Manin Festschrift" (eds. Yu.\ Tschinkel and Yu.\ Zarhin),
Progress in Mathematics, Birkha\"user (2008).


\bibitem[May]{May}  J.P. May. The Geometry of Iterated Loop Spaces, volume 271 of Lecture Notes in Mathematics. Springer-
Verlag, New York, 1972.

\bibitem[Mc]{Mc} S.\ McLane, {\em Categorical algebra},
Bull.\ Amer.\ Math.\ Soc. {\bf 71} (1965), 40-106.

\bibitem[Me1]{Me0} S.A.~Merkulov,~{\em Strong homotopy algebras of a K\"ahler manifold},
Intern.\ Math.\ Res.\ Notices {\bf } (1999), 153-164.


\bibitem[Me2]{Me-F}  S.A.~Merkulov,~{\em Operads, deformation theory and $F$-manifolds}.
In:  Frobenius manifolds ( Festschrift for Yu.I. Manin),
 213--251, Aspects Math., E36, Vieweg, Wiesbaden, 2004.



\bibitem[Me3]{Me1} S.A.\ Merkulov, {\em Prop profile of Poisson geometry},
 Commun.\ Math.\ Phys. {\bf 262} (2006),
117-135.

\bibitem[Me4]{Me2} S.A.\ Merkulov, {\em Nijenhuis infinity and contractible
dg manifolds}, Compositio Mathematica {\bf 141} (2005), 1238-1254.

\bibitem[Me5]{Me3} S.A.\ Merkulov, {\em Graph complexes with loops and wheels}.
In: ``Algebra, Arithmetic
and Geometry - Manin Festschrift" (eds. Yu.\ Tschinkel and Yu.\ Zarhin),
Progress in Mathematics, Birkha\"user (2008)  .


\bibitem[Me6]{Me4} S.A.\ Merkulov, {\em Lectures on props,
Poisson geometry and deformation quantization}.
In: { Poisson Geometry in Mathematics and Physics}, Contemporary Mathematics vol.\ 450,
 AMS (2008), 223-257.



\bibitem[Me7]{Me-BV} S.A.\ Merkulov, {\em Wheeled pro(p)file of the Batalin-Vilkovisky
formalism},
preprint  arXiv:0804.2481 (2008).




\bibitem[MeVa]{MV} S.A.\ Merkulov and B.\ Vallette,
\newblock{\em Deformation theory of representations of prop(erad)s I}. J.\ Reine Angew.\ Math. {\bf 634}
 (2009), 51–106.


\bibitem[Mn]{Mn} P.\ Mnev,
\newblock{\em On simplicial BF theory},
Dokl.\ Akad.\ Nauk {\bf 418} (2008), 308-312.

\bibitem[Pe]{Pe} R.\ Penrose, {\em
Applications of negative dimensional tensors}. In: Combinatorial Mathematics and its Applications
(Proc. Conf., Oxford, 1969) pp. 221--244 Academic Press, London.

\bibitem[Sc]{Sc} A.\ Schwarz, {\em Geometry of Batalin-Vilkovisky quantization},
Commun.\ Math.\ Phys.\
{\bf 155} (1993), 249-260.

\bibitem[Sta]{St} J.D.\ Stasheff, {\em On the homotopy
  associativity of $H$-spaces, I \@ II}, Trans.\ Amer.\
  Math.\ Soc.\ {\bf 108} (1963), 272-292 \& 293-312.

\bibitem[Str1]{Str} H.\ Strohmayer, {\em Operad profiles of Nijenhuis structures}. To appear in
{J.\ Diff.\ Geom.\ and Applications}.

\bibitem[Str2]{Str2} H.\ Strohmayer, {\em Prop profile of bi-Hamiltonian structures}.
To appear in { J.\ of Nocommutative Geom.}.






\bibitem[Va1]{V}
B.\ Vallette,  {\em A Koszul duality for
props}, Trans.\ AMS, {\bf 359} (2007), 4865-4943.

\bibitem[Va2]{V2}
B.\ Vallette, {\em Homology of generalized partition posets}, J. Pure Appl. Algebra,
208(2):699–725, 2007.


\bibitem[vdL]{vdL} P.\  van der Laan,
    {\em Operads up to homotopy and deformations of operad maps},
    math.QA/0208041.


\bibitem[We]{We} A.\ Weinstein, {\em
The modular automorphism group of a
Poisson Manifold}, J.\  Geom. Phys. {\bf 23} (1997), 379-394.





  \end{thebibliography}
  \end{document}